\documentclass[a4paper]{article}

\usepackage{bbm}
 
\usepackage[all]{xy}\usepackage[latin1]{inputenc}  
\usepackage[dvips]{graphics,graphicx}
\usepackage{amsfonts,amssymb,amsmath,color,mathrsfs, amstext,bbm}
\usepackage{amsbsy, amsopn, amscd, amsxtra, amsthm,authblk, enumerate}
\usepackage{upref}
\usepackage{geometry}
\usepackage{bm}
\usepackage{dsfont}
\usepackage{todonotes}
\usepackage{mathtools}
\usepackage{appendix}
\usepackage{float}

\usepackage[all]{xy}\usepackage[latin1]{inputenc}        
\usepackage[dvips]{graphics,graphicx}
\usepackage{amsfonts,amssymb,amsmath,xcolor,mathrsfs, amstext}
\usepackage{amsbsy, amsopn, amscd, amsxtra, amsthm,authblk, enumerate}
\usepackage{upref}
\usepackage{geometry}
\usepackage{float}

\usepackage[colorlinks,
            linkcolor=blue,
            anchorcolor=green,
            citecolor=blue
            ]{hyperref}

\numberwithin{equation}{section}

\def\e{\varepsilon}
\def\epsilon{\varepsilon}
\def\eps{\varepsilon}

\newcommand{\ol}{\overline}
\newcommand{\wt}{\widetilde}
\def\alb#1\ale{\begin{align*}#1\end{align*}}
\newcommand{\eqb}{\begin{equation}}
\newcommand{\eqe}{\end{equation}}

\DeclareMathOperator{\dist}{dist}

\DeclareMathOperator{\loc}{loc}

\newcommand{\inner}[2]{\langle {#1}, {#2}\rangle}
\newcommand{\bbC}{\mathbb{C}}
\newcommand{\bbD}{\mathbb{D}}
\newcommand{\bbE}{\mathbb{E}}
\newcommand{\bbH}{\mathbb{H}}

\newcommand{\bbR}{\mathbb{R}}
\newcommand{\bbP}{\mathbb{P}}

\newcommand{\cD}{\mathcal{D}}

\newcommand{\cL}{\mathcal{L}}

\newcommand{\cW}{\mathcal{W}}
\newcommand{\cO}{\mathcal{O}}

\newcommand{\QD}{\mathrm{QD}}
\newcommand{\GQD}{\mathrm{GQD}}

\newcommand{\LF}{\mathrm{LF}}
\newcommand{\LP}{\mathrm{LP}}
\newcommand{\SLE}{\mathrm{SLE}}
\newcommand{\mSLE}{\mathrm{mSLE}}

\newcommand{\Wd}{\mathrm{Weld}}
\newcommand{\Md}{{\mathcal{M}}^\mathrm{disk}}
\newcommand{\Mfd}{{\mathcal{M}}^\mathrm{f.d.}}

\newtheorem{theorem}{Theorem}[section]
\newtheorem{lemma}[theorem]{Lemma}
\newtheorem{proposition}[theorem]{Proposition}
\newtheorem*{proposition*}{Proposition}
\newtheorem{corollary}[theorem]{Corollary}
\newtheorem*{corollary*}{Corollary}
\newtheorem{remark}[theorem]{Remark}

\newtheorem{definition}[theorem]{Definition}
\newtheorem*{definitions*}{Definitions}

\newtheorem*{example*}{\bf Example}

\numberwithin{equation}{section}

\title{Conformal welding of quantum disks and multiple SLE:\\ the non-simple case}
\author{Morris Ang\thanks{Department of Mathematics, UC San Diego} \qquad Nina Holden\thanks{Courant Institute of Mathematical Sciences, New York University} \qquad Xin Sun\thanks{Beijing International Center for Mathematical Research, Peking University.}  \qquad Pu Yu\thanks{Courant Institute of Mathematical Sciences, New York University}}
\date{\today}

\begin{document}

\maketitle

\begin{abstract}

    Two-pointed quantum disks with a weight parameter $W>0$ is a canonical family of finite-volume random surfaces in Liouville quantum gravity. We prove that the conformal welding of the forested variant of this disk gives a two-pointed quantum disk with an independent SLE$_\kappa$ for $\kappa\in(4,8)$. Furthermore, we show that the conformal welding of multiple forested quantum disks gives a surface arising in Liouville conformal field theory decorated by multiple SLE$_\kappa$ for $\kappa\in(4,8)$, such that the random conformal modulus contains the SLE partition function as a multiplicative factor. In partiuclar, this gives a construction of the multiple SLE$_\kappa$ associated with any given link pattern. As a corollary, for $\kappa\in(4,8)$, we prove the existence of the multiple SLE partition functions, which are smooth functions satisfying a system of PDEs and conformal covariance. This was open for $\kappa \in  (6,8)$ and $N\ge  3$ prior to our work.
\end{abstract}




\section{Introduction}
In the last two decades, two dimensional random conformal geometry has been an active area of research in probability theory. This article concerns the connections between the three central topics in this area: the Schramm-Loewner evolution ($\SLE_\kappa$), Liouville quantum gravity (LQG) and Liouville conformal field theory (LCFT). $\SLE_\kappa$  is an important 
family of random non-self-crossing curves introduced by Schramm~\cite{Sch00},
which are natural candidates to describe the scaling limits of various two-dimensional lattice models at criticality, e.g.~\cite{smirnov2001critical,lawler2011conformal,schramm2009contour,chelkak2014convergence}. LQG is introduced by Polyakov in his seminal work~\cite{polyakov1981quantum} with parameter $\gamma\in(0,2]$, and has been shown to describe the scaling limits of a large class of random planar maps, see e.g.~\cite{LeGall13, BM17,HS19,gwynne2021convergence}. As the fundamental building block of LQG, LCFT is the 2D quantum field theory which is made rigorous by~\cite{DKRV16} and later works. 
See~\cite{Law08,vargas-dozz,GHS19,BP21,gwynne2020random,She22} for more background on these topics.

One of the deepest results in random planar geometry is the \emph{conformal welding} of random surfaces.  {R}oughly speaking, when we glue two $\gamma$-LQG surfaces together, we get another  {$\gamma$-}LQG surface decorated by an $\SLE_\kappa$ curve with $\kappa=\gamma^2$. This type of result was first proved in Sheffield's quantum zipper paper~\cite{She16a} and extended to a broader class of infinite area surfaces in~\cite{DMS14}. In~\cite{AHS23}, similar results were proved for a  class of canonical finite area LQG surfaces called (two-pointed) quantum disks. When $\gamma\in(\sqrt{2},2)$ and $\kappa = 16/\gamma^2$, it is  shown in~\cite{DMS14,MSW22CLE} that certain   LQG surfaces with non-simple boundaries, or \emph{generalized LQG surfaces} can be conformally welded together with the interface being $\SLE_\kappa$ and $\mathrm{CLE}_\kappa$ curves. Our first result is the conformal welding of the generalized quantum disks {in} this setting, which extends~\cite{AHS23} to the $\kappa\in(4,8)$ regime.

The convergence of interfaces in statistical physics model{s} to $\SLE_\kappa$ can {often} be extended to the  {case}  of multiple curves~\cite{Izy17,KS17scaling,peltola2019global}, which give{s} rise to the notion of multiple $\SLE_\kappa$~\cite{bauer2005multiple}. For $\kappa\in(0,4]$, the multiple SLE is well-studied in both simply and multiply connected domains~\cite{dubedat2007commutation,graham2007multiple,kytola2016pure,kozdron2006configurational,lawler2009partition,jahangoshahi2018multiple}, and admits a natural partition function. On the other hand, the $\kappa\in(4,8)$ regime is far less understood. \cite{wu2020hypergeometric,peltola2019toward} give a probabilistic construction of \emph{global multiple} $\SLE_\kappa$ for $\kappa\in(4,6]$, and it is proved that the Loewner equations driven by the multiple SLE partition functions generate the local multiple SLEs. In a recent work~\cite{Zhan23}, Zhan gave a construction of the multiple SLE as a $\sigma$-finite measure for $\kappa\in(4,8)$ based on similar ideas, and proved the uniqueness of the measure. For $\kappa\in(6,8)$, the finiteness of the multiple $N$-$\SLE$   measure when $N\ge3$ remained an open problem.

Our second main result is the construction of the global multiple SLE measure for $\kappa\in(4,8)$ from the conformal welding of multiple generalized quantum disks. In a concurrent work~\cite{SY23} by the third and fourth authors, it is shown that the conformal welding of a certain collection of LQG surfaces can be described by LCFT decorated with multiple SLE curves, and the density of the random moduli  is given by the LCFT partition function 
times the SLE partition function for the interfaces. In Theorem~\ref{thm:main}, we prove an analogous result in the $\kappa\in(4,8)$ setting, and the multiple SLE measure there agrees with the ones in~\cite{wu2020hypergeometric,peltola2019toward,Zhan23}. We then further infer that the partition function is finite, which completes the existence and uniqueness of the multiple SLE for $\kappa\in(4,8)$. Moreover, we will show that as probability
measures on curve segments, the global multiple SLE for $\kappa\in(4,8)$ agrees with the  {local multiple $\SLE_\kappa$} driven by {the global multiple SLE}  partition functions. 


 
\subsection{Multiple SLE and partition functions}\label{subsec:intro-sle}
The chordal $\SLE_\kappa$ in the upper half plane $\bbH$ is a probability measure $\mu_\bbH(0,\infty)$ on non-crossing curves from 0 to $\infty$ which is scale invariant and satisfies the domain Markov property. The $\SLE_\kappa$ curves are simple when $\kappa\in(0,4]$, non-simple and non-space-filling for $\kappa\in(4,8)$, and space-filling when $\kappa\ge8$. 
By conformal invariance, for a simply connected domain $D$ and $z,w\in\partial D$ distinct, one can define the $\SLE_\kappa$ probability measure $\mu_D(z,w)^\#$ on $D$ by taking conformal maps $f:\bbH\to D$ where $f(0)=z$, $f(\infty)=w$. For  {$\rho_-,\rho_+>-2$, $\SLE_\kappa(\rho_-;\rho_+)$}  {is a variant of $\SLE_\kappa$ 
studied in numerous works}, e.g.~\cite{LSW03,Dub05,MS16a}. 

For $N>0$, consider $N$ disjoint simple curves in $\overline{\bbH}$ connecting $1,2,...,2N\in\partial\bbH$. Topologically, these $N$ curves form a planar pair
partition, which we call a link pattern and denote by  $\alpha=\{\{i_1,j_1\},...,\{i_N,j_N\}\}$. The pairs $\{i,j\}$ in $\alpha$ are called links, and the set of link patterns with $N$ links is denoted by $\mathrm{LP}_N$.

Let $(D;x_1,...,x_{2N})$ be a topological polygon, in the sense that $D\subset\bbC$ is a simply connected domain and $x_1, ..., x_{2N} \in\partial D$ are 
{$2N$} distinct boundary points appearing in counterclockwise
order on boundary segments. 
In this paper, we work  {with} polygons where $\partial D$ is smooth near each of $x_j$ for $j=1,...,2N$. Consider a link pattern $\alpha= \{\{i_1,j_1\},...,\{i_N,j_N\}\}\in\mathrm{LP}_N$. Let $X_\alpha(D;x_1,...,x_{2N})$ be the space of $N$ non-crossing continuous curves $(\eta_1,...,\eta_N)$ in $\ol D$ such that  for each $1\le k\le N$, $\eta_k$ starts  {at} $x_{i_k}$, ends at $x_{j_k}$, and does not partition $D$ into components where some {pair of boundary points} corresponding to
a link in $\alpha$  belong to different components.

 Now we introduce the following definition of the global multiple SLE.
\begin{definition}[\cite{beffara2021uniqueness,Zhan23}]\label{def:msle}
    Let $\kappa\in(0,8)$. Let $N\ge1$ and fix a link pattern $\alpha\in\LP_N$. We call a probability measure on the families of curves $(\eta_1,...,\eta_N)\in X_\alpha(D;x_1,...,x_{2N})$ an \emph{$N$-global $\SLE_\kappa$ associated with $\alpha$}, if for each $1\le k\le N$, the conditional law of the curve $\eta_k$ given $\eta_1,...,\eta_{k-1},\eta_{k+1},...,\eta_N$ is the chordal $\SLE_\kappa$ connecting $x_{i_k}$ and $x_{j_k}$ in the connected component of the domain $D\backslash\cup_{k'\neq k}\eta_{k'}$ containing the endpoints $x_{i_k}$ and $x_{j_k}$ of $\eta_k$ on its boundary.
\end{definition}

\begin{theorem}\label{thm:existence-uniqueness}
    Let $\kappa\in(0,8)$, $N\ge1$ and $\alpha\in\LP_N$. Let $(D;x_1,...,x_{2N})$ be a topological polygon. Then there exists a unique $N$-global $\SLE_\kappa$ associated with $\alpha$, which we denote by $\mSLE_{\kappa,\alpha}(D;x_1,...,x_{2N})^\#$.
\end{theorem}

 Theorem~\ref{thm:existence-uniqueness} is already known when $\kappa\in(0,4]$~\cite{kozdron2006configurational,peltola2019global,beffara2021uniqueness} or $N=2$~\cite{miller2018connection}. {For $\kappa \in (4,8)$,}  uniqueness is proved in~\cite{Zhan23}, and for $\kappa\in(4,6]$,  existence is shown in~\cite{wu2020hypergeometric,peltola2019toward}. For $\kappa\ge8$, the existence  is trivial while the uniqueness fails. The remaining part is the existence for $\kappa\in(6,8)$, which shall be proved via Theorem~\ref{thm:main}.

One naturally associated object is the \emph{multiple SLE  partition function}. For $\kappa\in(0,8)$, let $$b = \frac{6-\kappa}{2\kappa}$$ 
be the conformal weight. Fix $N>0$. Let $\mathfrak{X}_{2N} = \{(x_1,...,x_{2N})\in\bbR^{2N}:x_1<...<x_{2N}\}$ be the configuration space. Following~\cite{dubedat2007commutation}, a  multiple $\SLE_\kappa$ partition function is a positive smooth function $\mathcal{Z}:\mathfrak{X}_{2N}\to \bbR_{+}$  satisfying the following two properties:
\begin{itemize}
    \item[(PDE)]\label{it:msle-pdf} \emph{Partial differential equations of second order}: We have
\begin{equation}\label{eq:msle-pde}
    \bigg[\frac{\kappa}{2}\partial_i^2+\sum_{j\neq i}\big(\frac{2}{x_j-x_i}\partial_j-\frac{2b}{(x_j-x_i)^2} \big)   \bigg]\mathcal{Z}(\bbH;x_1,...,x_{2N}) = 0\ \ \ \text{for\ } i=1,{\dots},2N.
\end{equation}
\item[(COV)]\label{it:msle-conformal-conf} \emph{M\"{o}bius covariance}: For any conformal map $f:\bbH\to\bbH$ with $f(x_1)<...<f(x_{2N})$,
\begin{equation}\label{eq:msle-conformal-conf}
    \mathcal{Z}(\bbH;x_1,...,x_{2N}) = \prod_{i=1}^{2N}f'(x_i)^b\times \mathcal{Z}(\bbH;f(x_1),...,f(x_{2N})).
\end{equation}
\end{itemize}


By~\eqref{eq:msle-conformal-conf}, the notion of partition function can be extended to other simply connected domains. Let $(D;x_1,...,x_{2N})$ be a topological polygon. Then  for any conformal map $f:D\to f(D)$, one has
\begin{equation}\label{eq:partition-conformal-conf}
    \mathcal{Z}(D;x_1,...,x_{2N}) = \prod_{i=1}^{2N}f'(x_i)^b\times \mathcal{Z}(D;f(x_1),...,f(x_{2N})).
\end{equation}

The {multiple SLE partition functions} are  related to another approach to construct the multiple $\SLE_\kappa$, namely the local $N$-$\SLE_\kappa$. One generate{s} several $\SLE_\kappa$ curves by describing their time evolution via Loewner chains. It is shown in~\cite{dubedat2007commutation} that the {local $N$-$\SLE_\kappa$} can be classified by the partition function $\mathcal{Z}$ in terms of Loewner driving functions, while~\cite[Theorem 1.3]{peltola2019global} proved that the global $N$-$\SLE_\kappa$ agree with local  $N$-$\SLE_\kappa$ when $\mathcal{Z}=\mathcal{Z}_\alpha$ and $\kappa\in(0,4]$. See Section~\ref{subsec:local-msle} for a detailed discussion.

The multiple SLE partition functions are constructed explicitly in~\cite{peltola2019global} for $\kappa\in(0,4]$, ~\cite{wu2020hypergeometric} for $\kappa\in(4,6]$, and ~\cite{KP16partitionfunc} (relaxing the positivity constraint) for  $\kappa\in(0,8)\backslash\mathbb{Q}$. 
 For $\kappa\in(6,8)$, we have the following result.
\begin{theorem}\label{thm:partition-func}
    Let $\kappa\in(6,8)$ and $N\ge1$. Then for each link pattern $\alpha\in\mathrm{LP}_N$, there exists an associated positive smooth function $\mathcal{Z}_\alpha:\mathfrak{X}_{2N}\to \bbR_{+}$ satisfying (PDE) and (COV).    Moreover, the local  $N$-SLE$_\kappa$ driven by $\mathcal{Z}_\alpha$ agrees with the initial segments of the global $N$-SLE$_\kappa$.
\end{theorem}

Prior to Theorem~\ref{thm:partition-func}, the existence of the multiple SLE partition function for $\kappa\in(6,8)$ was unknown for $N\ge3$. One major  difficulty is that $b<0$ in this range, and the current technical estimates are insufficient for building the partition function directly as in~\cite{peltola2019global,wu2020hypergeometric}. 

For a conformal map $\varphi:D\to \tilde D$ and a measure $\mu(D;x,y)$ on continuous curves from $x$ to $y$ in $\overline{D}$, we write $\varphi\circ\mu(D,x,y)$ for the law of $\varphi\circ\eta$ when $\eta$ is sampled from $\mu(D;x,y)$. Given Theorems~\ref{thm:existence-uniqueness} and~\ref{thm:partition-func}, we define the measure 
{
\eqb\label{eq-mSLE}
\mSLE_{\kappa, \alpha}(D;x_1,...,x_{2N}) = \mathcal{Z}_\alpha(D;x_1,...,x_{2N})\times \mSLE_{\kappa, \alpha}(D;x_1,...,x_{2N})^\#.
\eqe
}
Then we have the following conformal covariance
\begin{equation}\label{eq:msle-conformal-covariance}
     f\circ \mSLE_{\kappa,\alpha}(D;x_1,...,x_{2N}) = \prod_{i=1}^{2N}f'(x_i)^b\times\mSLE_{\kappa,\alpha}(f(D);f(x_1),...,f(x_{2N}))
\end{equation}
whenever the boundaries of  $D$ and $f(D)$ are smooth near the marked points.

In Theorems~\ref{thm:existence-uniqueness} and~\ref{thm:partition-func}, the measure $ \mSLE_{\kappa,\alpha}(D;x_1,...,x_{2N})^\#$ and the   partition function $\mathcal{Z}_\alpha(D;x_1,...,x_{2N})$ will be defined via a cascade relation as in~\cite{wu2020hypergeometric,Zhan23}; see Section~\ref{subsec:pre-msle}. Under this inductive definition, we will first prove Theorem~\ref{thm:main} {below}, and infer that the partition function $\mathcal{Z}_\alpha(D;x_1,...,x_{2N})$ is finite, which further completes the induction for the proof of  Theorems~\ref{thm:existence-uniqueness} and~\ref{thm:partition-func}. Moreover, in Section~\ref{subsec:local-msle}, we will show that the local multiple $\SLE_\kappa$ driven by the  partition function $\mathcal{Z}_\alpha$ agrees with the global multiple $\SLE_\kappa$ associated to $\alpha$.

\subsection{Liouville quantum gravity surfaces and conformal welding}\label{subsec:intro-lqg}
Let $D\subset\bbC$ be a simply connected domain. The Gaussian Free Field (GFF) on $D$ is the centered Gaussian process on $D$ whose covariance kernel is the Green's function~\cite{She07}. For $\gamma\in(0,2)$ and $\phi$ a variant of the GFF, the $\gamma$-LQG area measure in $D$ and length measure on $\partial D$ is roughly defined by $\mu_\phi(dz)=e^{\gamma\phi(z)}dz$ and $\nu_\phi(dx) = e^{\frac{\gamma}{2}\phi(x)}dx$, and are made rigorous by regularization and renormalization~\cite{DS11}. Two pairs $(D,h)$  and $(D',h')$ represent the same quantum surface if there is a conformal map between $D$ and $D'$ preserving the geometry{; see the discussion around \eqref{eq:lqg-changecoord}}.

For $W>0$, the two-pointed quantum disk of weight $W$, whose law is denoted by $\Md_2(W)$, is a quantum surface with two boundary marked points introduced in~\cite{DMS14,AHS23}, which has finite quantum area and length.  The surface is simply connected when $W\ge\frac{\gamma^2}{2}$, and consists of a chain of countably many weight $\gamma^2-W$ quantum disks when $W\in(0,\frac{\gamma^2}{2})$. For the special case $W=2$, the two boundary marked points are \emph{quantum typical}  {with respect to} the LQG boundary length measure~\cite[Proposition A.8]{DMS14}. 

As shown in~\cite{cercle2021unit,AHS21}, the quantum disks can be alternative{ly} described in terms of LCFT. The \emph{Liouville field} $\LF_\bbH$ is an infinite measure on the space of generalized functions on $\bbH$ obtained by an additive perturbation of the GFF. For $i=1,...,m$ and $(\beta_i,s_i)\in\bbR\times\partial\bbH$, we can make sense of the measure $\LF_\bbH^{(\beta_i,s_i)_i}(d\phi) = \prod_i e^{\frac{\beta_i}{2}\phi(s_i)}\LF_\bbH(d\phi)$ via  regularization and renormalization, which leads to the notion of \emph{Liouville fields with boundary insertions}. See Definition~\ref{def-lf-H-bdry} and Lemma~\ref{lm:lf-insertion-bdry}. 

Next we briefly recall the \emph{generalized quantum surfaces} (or forested quantum surfaces) from~\cite{DMS14,MSW22CLE}, which is based on the construction of the loop-trees    in~\cite{CK13looptree}. Fix $\gamma\in(\sqrt{2}, 2)$. Given an $\alpha$-stable L\'{e}vy process $(X_t)_{t>0}$ with no negative jumps and $\alpha = \frac{4}{\gamma^2}\in(1,2)$, we first determine the tree structure as in~\cite{CK13looptree}.  This can be done by gluing a Poisson point process of loop-trees to one side of the half real line.   Then we assign a conformal structure to each of the loops using a (standard) independent quantum disk. This defines  the \emph{forested line} as in~\cite[Definition 1.14]{DMS14}. Points on the line corresponds to the running infimum of $(X_t)_{t>0}$, which we parameterize by the LQG length measure; the boundaries of the loop trees are parameterized by \emph{generalized quantum length}~\cite{MSW22CLE}. The generalized quantum surfaces are then constructed by truncating and gluing independent forested lines to the boundary arcs of classical quantum surfaces.   
See Section~\ref{subsec:gqd} for more details.

For $W>0$, we write $\Mfd_2(W)$ for the law of the generalized quantum surface obtained by truncating and gluing independent forested lines to both sides of its boundary according to the LQG boundary length. We call a sample from $\Mfd_2(W)$ a \emph{weight $W$ forested quantum disk.} For the special weight $W = \gamma^2-2$, the two marked points are quantum typical {with respect to} the generalized LQG length measure (see Proposition~\ref{prop:def-gqd}). By sampling additional marked points from the generalized quantum length measure, we obtain a multiply marked generalized quantum disk. For $m\ge1$, we write $\GQD_{m}$ for the law of the generalized quantum disk with $m$ marked points on the boundary sampled from the generalized quantum length measure; see Definition~\ref{def:GQD} for a precise description.

Given a pair of classical quantum surfaces, following ~\cite{She16a,DMS14}, there exists a way to \emph{conformally weld} them together according to the length measure provided that the interface lengths agree; see e.g.~\cite[Section 4.1]{AHS21} and~\cite[Section 4.1]{ASY22} for more explanation.  In~\cite[Theorem 1.15]{DMS14}, it is proved that for $\kappa=\frac{16}{\gamma^2}$,  by drawing an $\SLE_\kappa(\frac{\kappa}{2}-4;\frac{\kappa}{2}-4)$ curve $\eta$ on an independent a weight $2-\frac{\gamma^2}{2}$ quantum wedge $\mathcal{W}$, one cuts the wedge into two independent forested lines $\cL_\pm$ whose boundaries are identified via the generalized quantum length. Moreover, $(\mathcal{W}, \eta)$ is measurable with respect to $\cL_\pm$. {Note in particular that in \cite{DMS14}, $\cL_\pm$ determine $(\mathcal{W}, \eta)$ in a weaker sense than what we have in conformal welding. A somewhat stronger notion of uniqueness was proven in \cite{McEnteggart-Miller-Qian}. Finally, in light of the recent work~\cite{kavvadias2023conformal} on conformal removability of non-simple SLEs for $\kappa\in(4,\kappa_0)$, where $\kappa_0$ is some constant in $(4,8)$, it is possible to identify the recovery of $(\mathcal{W}, \eta)$ from $\cL_\pm$ as actual conformal welding as in the $\kappa\in(0,4)$ case.} 

In Proposition~\ref{prop:weld:segment}, we prove the analog of~\cite[Theorem 1.15]{DMS14} for $\Md_2(2-\frac{\gamma^2}{2})$ and segments of forested lines. Following this weaker notion of uniqueness, we define the conformal welding of forested lines or segments of forested lines to be this procedure of gluing them together to get a quantum wedge or quantum disk decorated by $\SLE_\kappa(\frac{\kappa}{2}-4;\frac{\kappa}{2}-4)$ curves. 

The conformal welding operation discussed as above naturally extends to generalized quantum surfaces. Let $\mathcal{M}^1, \mathcal{M}^2$ be measures on the space of generalized quantum surfaces with boundary marked points. For $i=1,2$, fix some boundary arcs $e_1, e_2$ such that $e_i$ are boundary arcs of finite generalized quantum length on samples from $\mathcal{M}^i$, and define the measure $\mathcal{M}^i(\ell_i)$ via the disintegration 
\eqb\label{eq:def-weld-general-qs}\mathcal{M}^i = \int_0^\infty \mathcal{M}^i(\ell_i)d\ell_i\eqe
over the generalized quantum lengths of $e_i$. For $\ell>0$, given a pair of surfaces sampled from the product measure $\mathcal{M}^1(\ell)\times \mathcal{M}^2(\ell)$, we can first weld the forested line segments of $e_1$ and $e_2$ together according to the generalized  quantum length to get an $\SLE_\kappa(\frac{\kappa}{2}-4;\frac{\kappa}{2}-4)$-decorated quantum disk $\cD$, and then weld $\cD$ to the remaining parts of $\mathcal{M}^1$ and $\mathcal{M}^2$ by the conformal welding of classical quantum surfaces. This yields a single surface decorated with an interface from the gluing.  
We write $\text{Weld}( \mathcal{M}^1(\ell), \mathcal{M}^2(\ell))$ for the law of the resulting curve-decorated surface, and let $$\text{Weld}(\mathcal{M}^1, \mathcal{M}^2):=\int_{\bbR}\, \text{Weld}(  \mathcal{M}^1(\ell), \mathcal{M}^2(\ell))\,d\ell$$ be the  welding of $\mathcal{M}^1,  \mathcal{M}^2$ along the boundary arcs $e_1$ and $e_2$. The case where both $e_1$ and $e_2$ have infinite generalized quantum length can be treated analogously.  {
 {We note that the curve-decorated quantum surface obtained from the above procedure is measurable with respect to the two generalized quantum surfaces being welded,}
thanks to the measurability for gluing forested line segments along with the uniqueness of the conformal welding for classical quantum surfaces.}  {In our setting, the latter uniqueness follows from the conformal welding of classical quantum disks~\cite{AHS23,AHS21}; the general case follows from the conformal removability of $\SLE_\kappa$ curves for $\kappa<4$ as discussed in~\cite[Section 3.5]{DMS14}.} By induction, this definition extends to the welding of multiple generalized quantum surfaces, where we first specify some pairs of boundary arcs on the quantum surfaces, and then identify each pair of arcs according to the generalized quantum length. 



\begin{figure}[htb]
    \centering
    \includegraphics[scale = 0.65]{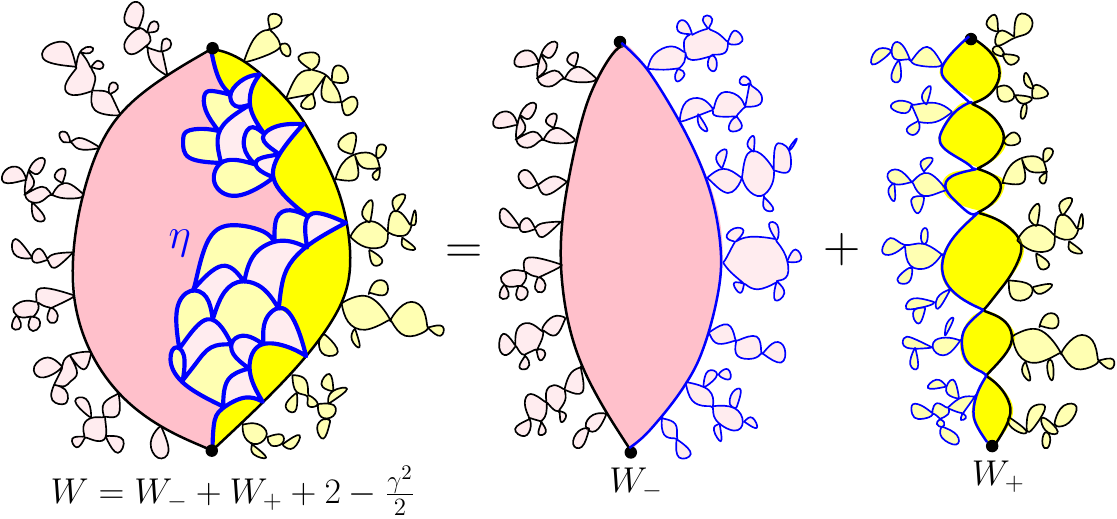}
    \caption{An illustration of Theorem~\ref{thm:fd-disk} in the case $W_-\ge\frac{\gamma^2}{2}$ and $W_+\in(0,\frac{\gamma^2}{2})$. If we draw an independent $\SLE_\kappa(\rho_-;\rho_+)$ curve on top of a generalized quantum disk of weight $W_-+W_+$, then the curve-decorated quantum surface is equal to the welding of a pair of  {weight} $W_1$ and $W_2$ generalized quantum disk{s} conditioned on having the same generalized quantum length for the interface. If $W<\frac{\gamma^2}{2}$, then the interface $\eta$ is understood as the concatenation of $\SLE_\kappa(\rho_-;\rho_+)$ curves in each bead of the weight $W$ forested quantum disk.} 
    \label{fig:fdweld}
\end{figure}

Now we state our result on the  welding of generalized quantum disks. See Figure~\ref{fig:fdweld} for an illustration. 

\begin{theorem}\label{thm:fd-disk}
     Let $\gamma\in(\sqrt{2},2)$ and $\kappa= \frac{16}{\gamma^2}$. Let {$W_-,W_+>0$} and $\rho_\pm = \frac{4}{\gamma^2}(2-\gamma^2+W_\pm)$. Let $W = W_++W_-+2-\frac{\gamma^2}{2}$. Then for some constant $c\in(0,\infty)$,
     \begin{equation}\label{eq:thm-fd-disk}
     \begin{split}
         \Mfd_2(W)\otimes\SLE_\kappa(\rho_-;\rho_+) = c\int_0^\infty \Wd(\Mfd_2(W_-;\ell),\Mfd_2(W_+;\ell))d\ell.
         \end{split}
     \end{equation}
  \end{theorem}

One immediate consequence of Theorem~\ref{thm:fd-disk} is the reversibility of $\SLE_\kappa(\rho_-;\rho_+)$ with $\kappa\in(4,8)$ and $\rho_\pm>\frac{\kappa}{2}-4$.   {We work on the case when $\rho_-=\rho_+=0$; as explained in~\cite[Lemma 3.1]{IGIII}, the reversibility of general $\SLE_\kappa(\rho_-;\rho_+)$ is equivalent to the reversibility of $\SLE_\kappa$ thanks to~\cite[Theorem 1.1]{MS16b} and~\cite[Proposition 7.30]{MS16a}.}  {In Proposition~\ref{prop:typical}, we will prove that the forested lines can be viewed as a Poisson point process 
 {of} generalized quantum disks, which further implies the reversibility of forested line segments. This gives the symmetry of forested quantum disks about their two marked points. 
Consider the welding on the right hand side of Figure~\ref{fig:fdweld} with $W_+=W_-=\gamma^2-2$ and let $(S,\eta)$ be the curve-decorated quantum surface, which has law $\Mfd_2(\frac{3\gamma^2}{2}-2)\otimes\SLE_\kappa$. On one hand, 
as we view the welding on the right hand side upside down, we get a  curve-decorated quantum surface $(S,\wt\eta)$, which also has law $\Mfd_2(\frac{3\gamma^2}{2}-2)\otimes\SLE_\kappa$, and $\wt\eta$ has opposite direction as $\eta$. We embed $S$ and $\wt S$ on $\bbH$. Since we are welding the same quantum surfaces, this  generates a homeomorphism $\varphi$ of $\bbH$ sending $\eta$ to $\wt\eta$ whose reflection is conformal in each connected component of $\bbH\backslash\eta$. Now by~\cite[Theorem 1.1]{McEnteggart-Miller-Qian}, the reflection of $\varphi$ must be a conformal automorphism of $\bbH$, indicating that $\wt\eta$ agrees with the time reversal of $\eta$ up to a conformal automorphism of $\bbH$. This verifies the reversibility of $\SLE_\kappa$.} 
Previously, the only known approach for the reversibility of $\SLE_\kappa$ when $\kappa\in(4,8)$ is through the imaginary geometry~\cite{IGIII}, while our conformal welding result provides a new perspective.  {[*Removed the whole plane $\SLE_\kappa(\rho)$ part since essentially this is just reversibility of $\SLE_\kappa(\kappa/2-4;\kappa/2-4)$*]}

\subsection{Multiple SLEs from   welding of quantum surfaces}

\subsubsection{Statement of the   welding result}
We start with the definition of the welding of generalized quantum disks according to link patterns. Fix $N\ge1$. Let $(D;x_1,...,x_{2N})$ be a topological polygon.  We draw $N$ disjoint simple curves $\tilde\eta_1,...,\tilde\eta_N$ according to  {a link pattern} $\alpha$, dividing $D$ into $N+1$ connected components $S_1,...,S_{N+1}$. For $1\le k\le N+1$, let $n_k$ be the number of  {points} on the boundary of $S_k$, and let $\tilde\eta_{k,1},...,\tilde\eta_{k,m_k}$ be the interfaces which are part of the boundary of $S_k$. Then for each $S_k$, we assign a generalized quantum disk with $n_k$ marked points on the boundary from $\mathrm{GQD}_{n_k}$, and consider the disintegration $\mathrm{GQD}_{n_k} =  {\int_{\bbR_+^{m_k}}}\mathrm{GQD}_{n_k}(\ell_{k,1},...,\ell_{k,m_k})\,d\ell_{k,1}...d\ell_{k,m_k}$ over the generalized quantum length of the boundary arcs  corresponding to the interfaces. For $(\ell_1,...,\ell_N)\in\bbR_+^N$, let $\ell_{k,j} = \ell_i$ if the interface $\tilde\eta_{k,j}=\tilde\eta_i$. We sample $N+1$ quantum surfaces from $\prod_{k=1}^{N+1}\mathrm{GQD}_{n_k}(\ell_{k,1},...,\ell_{k,m_k})$ and conformally weld them together by generalized LQG boundary length according to the link pattern $\alpha$, and write $\Wd_{\alpha}(\mathrm{GQD}^{N+1})(\ell_1,...,\ell_N)$ for the law of the resulting quantum surface decorated by $N$ interfaces. 
Define $\Wd_{\alpha}(\mathrm{GQD}^{N+1})$ by
$$ \Wd_{\alpha}(\mathrm{GQD}^{N+1}) = \int_{\bbR_+^N}\Wd_{\alpha}(\mathrm{GQD}^{N+1})(\ell_1,...,\ell_N)\,d\ell_1...d\ell_N. $$

 {Note that in the above definition, we have implicitly used the fact that the welding map is associative. Indeed, we may first apply Proposition~\ref{prop:weld:segment} to weld together all the forested line segments to get  weight $2-\frac{\gamma^2}{2}$ quantum disks. The associativity for generalized quantum surfaces is immediate from its counterpart for classical quantum surfaces.}
\begin{figure}[htb]
    \centering
    \begin{tabular}{cc} 
		\includegraphics[width=0.5\textwidth]{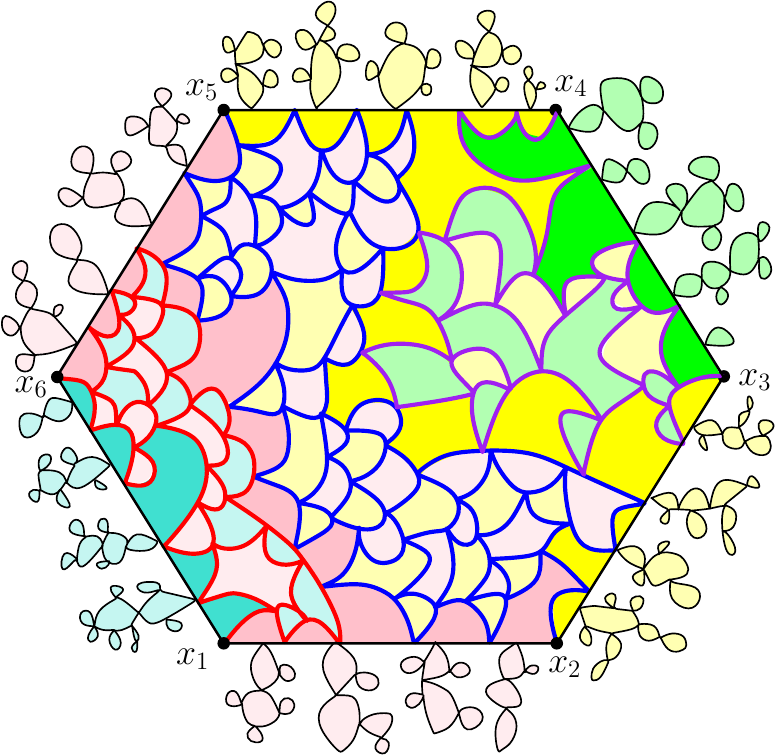}
		& 
		\includegraphics[width=0.5\textwidth]{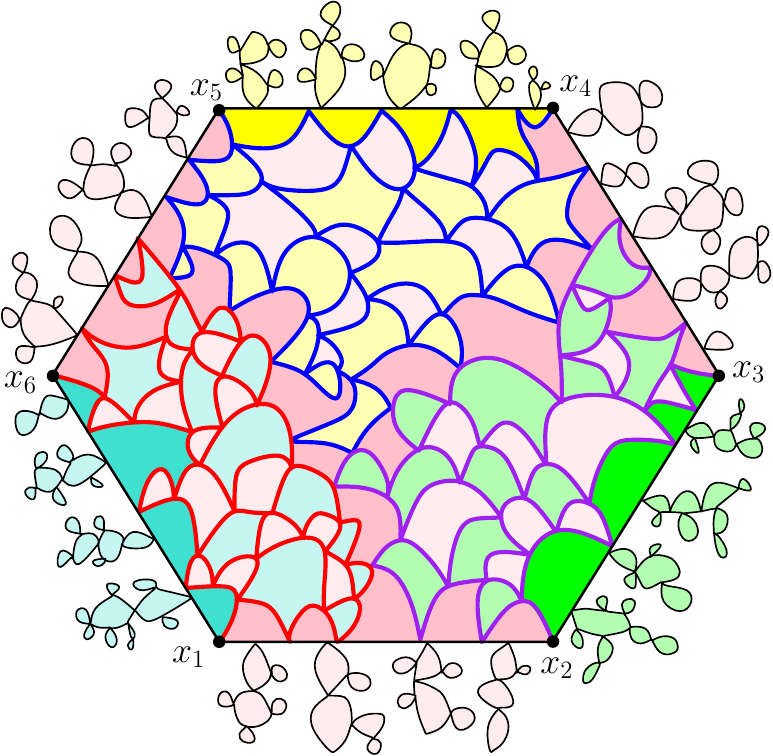}
	\end{tabular}
 \caption{An illustration of Theorem~\ref{thm:main}. \textbf{Left:} Under the link pattern $\alpha=\{\{1,6\},\{2,5\},\{3,4\} \}$, we are welding two samples from $\GQD_{2}$ (drawn in green and turquoise) with two samples from $\GQD_{4}$ (drawn in pink and yellow), restricted to the event $E$ where the welding output has the structure of a simply connected quantum surface glued to forested lines. {Each generalized disk is composed of a countable number of (regular) disks, and the (regular) disks that are used to connect two of the marked boundary points of the generalized disk are shown in dark color, while the other disks are shown in light color.} If we let $\ell_1,\ell_2,\ell_3$ be the interface lengths ordered from the left to the right, then the precise welding equation can be written as $\int_{\bbR_+^3}\Wd(\GQD_{2}(\ell_1),\GQD_{4}(\ell_1,\ell_2),\GQD_{4}(\ell_2,\ell_3),\GQD_{2}(\ell_3))\Big|_Ed\ell_1\,d\ell_2\,d\ell_3 $. \textbf{Right:} A similar setting where the link pattern $\alpha = \{\{1,6\},\{2,3\},\{4,5\} \}$, we are welding three samples from $\GQD_{2}$ (drawn in green, {turquoise} and yellow) with one sample from $\GQD_{6}$ (drawn in pink). The corresponding welding equation is given by 
 $\int_{\bbR_+^3}\Wd(\GQD_{2}(\ell_1),\GQD_{2}(\ell_2),{\GQD_{2}}(\ell_3),\GQD_{6}(\ell_1,\ell_2,\ell_3))\Big|_Ed\ell_1\,d\ell_2\,d\ell_3 $.
  The forested line part of each generalized quantum disk is drawn in a lighter shade.}\label{fig:main}
 \end{figure}


When $m\ge2$, for a sample  from $\GQD_{m}$, 
{consider} 
the paths of disks connecting the marked points. 
{Shrink each path of disks connecting each pair of marked points into a curve. The resulting set of curves form a tree $\mathcal{T}$  with at most $2m-3$ edges and the number of leaf nodes equal to $m$.}   Then $\Wd_\alpha(\GQD^{N+1})$ induces a topological gluing of trees, and, different from the simple case, the marked points are not in the same connected component if and only if in this gluing of trees, there exists an edge in some tree $\mathcal{T}$ not glued to any other edges. 

\begin{figure}[htb]
    \centering
    \begin{tabular}{cc} 
		\includegraphics[scale=0.5]{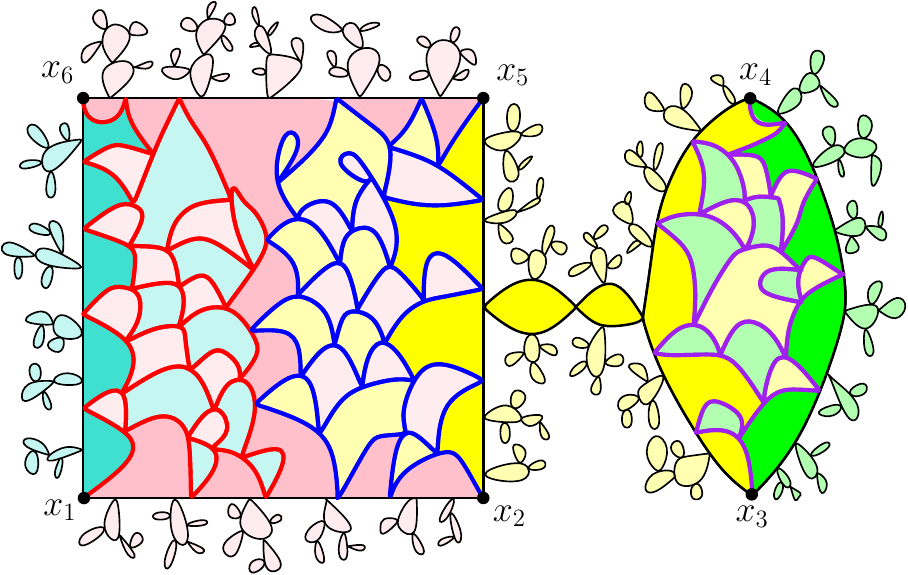} & \includegraphics[scale=0.67]{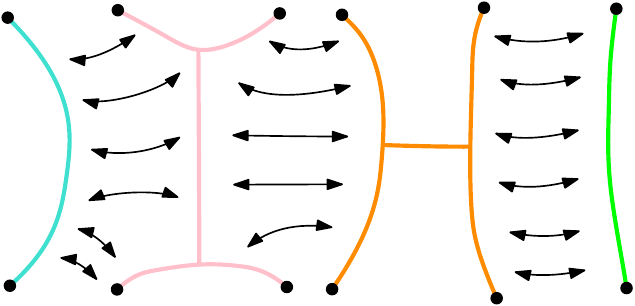}
		\\
		\includegraphics[scale=0.5]{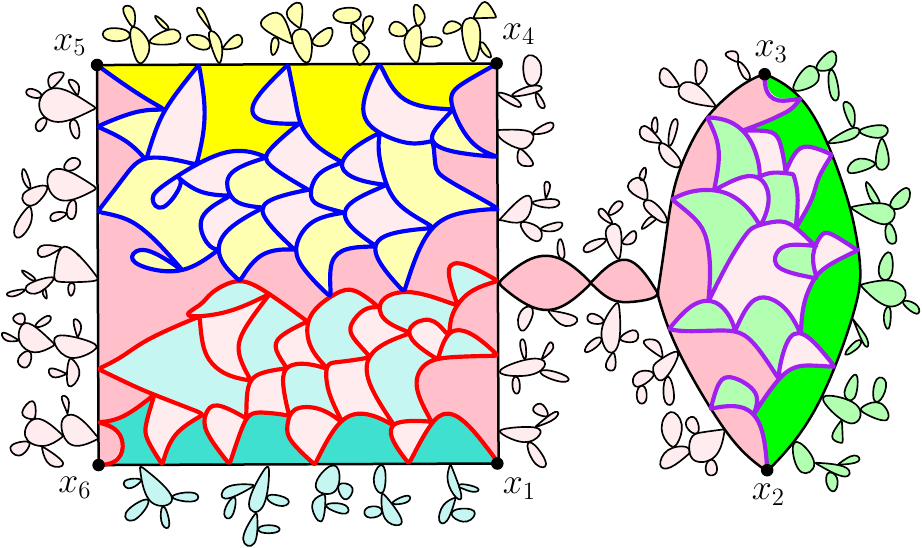} & \includegraphics[scale = 0.8]{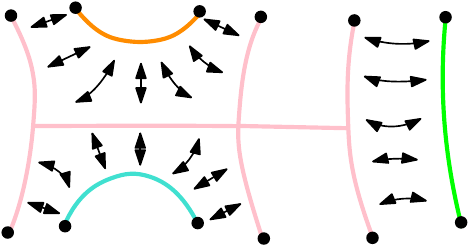}
	\end{tabular}
 \caption{An illustration of the case where the event $E$ in Theorem~\ref{thm:main} fails to happen. The left two panels illustrate weldings of  generalized quantum disks following the same link patterns as in Figure~\ref{fig:main}. In these examples the output surface  {cannot be described by forested lines glued to a single simply connected surface}. 
 The right two panels indicate  the corresponding topological gluing of trees, where in both cases there exists an edge which is not glued to any other edge.
 }\label{fig:non-connect}
 \end{figure}


We are now ready to state our result for non-simple multiple SLEs and conformal welding of generalized quantum disks. See Figure~\ref{fig:main} for an illustration.   {Recall the multiple SLE measure $\mSLE_{\kappa, \alpha}$ defined in~\eqref{eq-mSLE}.}

 \begin{theorem}\label{thm:main}
Let $\gamma\in(\sqrt{2},2)$, $\kappa=16/\gamma^2$ and $\beta = \frac{4}{\gamma}-\frac{\gamma}{2}$. Let $N\ge 2$ and $\alpha\in\mathrm{LP}_N$ be a link pattern. Let $c$ be the constant in Theorem~\ref{thm:fd-disk} for $W_-=W_+=\gamma^2-2$, and 
$c_N = \frac{\gamma^{2N-1} }{2^{2N+1}c^N(Q-\beta)^{2N}}$. Consider the curve-decorated quantum surface $(\bbH,\phi_N,0,y_1,...,y_{2N-3},1,\infty,\eta_1,...,\eta_N)/{\sim_\gamma}$ where $(\phi_N,y_1,...,y_{2N-3},\eta_1,...,\eta_N)$ has law
\begin{equation}\label{eq:thm-main}
\begin{split}
    c_N\int_{0<y_1<...<y_{2N-3}<1}\bigg[&\LF_\bbH^{(\beta,0),(\beta,1),(\beta,\infty),(\beta,y_1),...,(\beta,y_{2N-3})}(d\phi_N)\\&\times \mathrm{mSLE}_{\kappa,\alpha}(\bbH,0,y_1,...,y_{2N-3},1,\infty)(d\eta_1\,...\,d\eta_N)\bigg]\,dy_1\,...\,dy_{2N-3}.
    \end{split}
\end{equation}
If we truncate and glue an independent forested line to the quantum surface described above, then we obtain the conformal welding $\Wd_{\alpha}(\mathrm{GQD}^{N+1})$ restricted to the event
$E$  {where the interfaces all lie in a single region homeomorphic to the disk}.  

\end{theorem}

Theorem~\ref{thm:main} is the analog of~\cite[Theorem 1.1]{SY23}, and gives a concrete coupling between Liouville CFT and multiple SLE for $\kappa\in(4,8)$. One key observation is that the conformal weight  of the multiple SLE ($b = \frac{6-\kappa}{2\kappa}$) matches with that of Liouville CFT ($\Delta_\beta = \frac{\beta}{2}(\frac{2}{\gamma}+\frac{\gamma-\beta}{2})$) in the sense $b+\Delta_\beta = 1$.

\subsubsection{Overview of the proof}
Different from the existing works~\cite{kytola2016pure,peltola2019global,wu2020hypergeometric}, our proof of Theorems~\ref{thm:existence-uniqueness} and~\ref{thm:partition-func} rely on the novel coupling between Liouville CFT and multiple SLE in Theorem~\ref{thm:main}. Theorem~\ref{thm:fd-disk} has {independent} interest and further applications as well. We first prove  Theorem~\ref{thm:fd-disk}, and then prove Theorems~\ref{thm:existence-uniqueness},~\ref{thm:partition-func} and~\ref{thm:main} {via induction}. 

To show Theorem~\ref{thm:fd-disk},  by~\cite{AHS23}, it would suffice to work on the case of forested line segments (Proposition~\ref{prop:weld:segment}), which could be thought as the ``$W_\pm=0$" case. We begin with the welding of forested lines as in~\cite[Theorem 1.15]{DMS14}, and encode the locations of the cut points of the two independent forested lines by the zeros of a pair of independent squared Bessel processes $Z_\pm$ of dimension $2-\frac{\gamma^2}{2}$. Then by the additivity of squared Bessel processes, the locations of the cut points of the weight $2-\frac{\gamma^2}{2}$ quantum wedge are encoded by a squared Bessel process of dimension $4-\gamma^2\in(0,2)$. Then we use the Poisson point process description of quantum wedges to ``cut" a weight $2-\frac{\gamma^2}{2}$ quantum disk off the quantum wedge. The location of the cut points then will be encoded by  squared Bessel bridges, from which we infer the welding equation~\eqref{eq:weld:segment}. Finally, the measurablity {result follows} 
from~\cite[Theorem 1.15]{DMS14} by a local comparison. 

The proofs of the theorems related to multiple SLE are based on an induction. For $N=2$, Theorems~\ref{thm:existence-uniqueness} and ~\ref{thm:partition-func} are 
{known from previous works}~\cite{miller2018connection}. To prove Theorem~\ref{thm:main}, we begin with Theorem~\ref{thm:fd-disk} with $W_-=W_+=\gamma^2-2$. By definition, if we sample two marked points $x,y$ from the generalized quantum length measure and glue a sample from $\GQD_{2}$ along the boundary arc between $x,y$, then we obtain the desired conformal welding picture, and the $N=2$ case for Theorem~\ref{thm:main} follows from the techniques in~\cite{AHS21,SY23}. Now assume Theorems~\ref{thm:existence-uniqueness},~\ref{thm:partition-func} and~\ref{thm:main} hold for $2,\dots,N$. For $\alpha\in\LP_{N+1}$, we first construct the measure $\mSLE_{\kappa,\alpha}$ using the cascade relation in~\cite{wu2020hypergeometric,Zhan23}; see Section~\ref{subsec:pre-msle}. Then we prove Theorem~\ref{thm:main} for $N+1$ with the measure $\mSLE_{\kappa,\alpha}$.   To prove that the measure $\mSLE_{\kappa,\alpha}$ is finite, we 
consider the event $F_0$ where the generalized quantum lengths of all the boundary segments in the welding described in  Theorem~\ref{thm:main} 
are between 1 and 2. On the one hand, the fact that the two marked points on $\Mfd_2(\gamma^2-2)$ are typical to the generalized quantum length measure allows us to derive the joint law of boundary lengths of $\GQD_{n}$ (see Proposition~\ref{prop:gqd-bdry-length}), and it follows from a disintegration  {that} $\Wd_\alpha(\GQD^{N+2})[F_0]<\infty$. On the other hand,  in the expression~\eqref{eq:thm-main}, the event $F_0$  is not depending on     $\mSLE_{\kappa,\alpha}$. Therefore we infer that $\big|\mathrm{mSLE}_{\kappa,\alpha}(\bbH,0,y_1,...,y_{2N-1},1,\infty)\big|$ is locally integrable in terms of $(y_1,...,y_{2N-1})$, and is in fact smooth following a  hypoellipticity argument in~\cite{dubedat2015sle,peltola2019global}. This proves Theorem~\ref{thm:existence-uniqueness} and Theorem~\ref{thm:partition-func} for $N+1$ and completes the induction.

\subsection{Outlook and perspectives}
In this section, we discuss some ongoing works and future directions regarding the multiple SLE and the conformal welding of generalized quantum surfaces.

\begin{itemize}
    \item Based on conformal welding of classical quantum surfaces, with Remy the first three authors have proved a number of results on integrability of SLE and LCFT~\cite{AHS21,ARS21,AS21,ARS22}. Using the  conformal welding of generalized quantum surfaces, it is possible to extend the integrability results from~\cite{AS21,ARS22} to non-simple CLE. Moreover, in a forthcoming work by the first, third and fourth authors with Zhuang, we will give an explicit formula for the boundary touching probability of non-simple CLE, where Theorem~\ref{thm:fd-disk} is a critical input.

    \item In this paper, we focus on the chordal multiple SLE for $\kappa\in(4,8)$. For $\kappa\in(0,4]$, the multiple SLE on multiply connected domains is constructed in~\cite{jahangoshahi2018multiple}, and taking a limit yields the multiple radial SLE~\cite{healey2021n}. The existence of multiple SLE in general planar domains or radial setting remains open for $\kappa\in(4,8)$, and we believe that the conformal welding of generalized quantum surfaces can be applied to settle these problems. 
     {An interesting problem is to determine the SLE partition function in these settings.}
    \item Our construction of the partition function $\mathcal{Z}_\alpha$ is coherent with~\cite{wu2020hypergeometric,peltola2019toward} when $\kappa\in(4,6]$. In this range, $\{\mathcal{Z}_\alpha:\alpha\in\LP_N\}$ satisfy a certain asymptotic property {as two of the marked boundary points collide} and a strong power law bound as in~\cite[Theorem 1.1]{peltola2019global}, which uniquely specify the 
    partition functions by~\cite{FK15PDEa,FK15PDEb}. If one can prove the asymptotic property and the power law bound for $\kappa\in(6,8)$ as well, then $\{\mathcal{Z}_\alpha:\alpha\in\LP_N\}$ in Theorem~\ref{thm:partition-func} gives the 
    partition function for the range $\kappa\in(6,8)$ as in~\cite{peltola2019global}.

    {\item Following the SLE duality~\cite{zhan2008duality}, one may also consider the mixed multiple SLE with both $\SLE_{\tilde\kappa}$ and $\SLE_\kappa$ curves where $\tilde\kappa = \frac{16}{\kappa}\in(0,4)$, i.e., replace some of the curves in $\mSLE_{\kappa,\alpha}$ by $\SLE_{\tilde\kappa}$ curves.  It would be interesting to consider the partition functions in this setting, and their relations with the conformal field theory.} 
    
    
\end{itemize}

{We refer to the outlook of \cite{SY23} for further future directions related multiple-SLE and conformal welding.}

\vskip20pt

\noindent {\bf Acknowledgements}. We thank Eveliina Peltola and Dapeng Zhan  for helpful discussions. We thank two anonymous referees for their careful reading and many helpful comments. M.A.\ was supported by the Simons Foundation as a Junior Fellow at the Simons Society of Fellows. N.H.\ was supported by grant DMS-2246820 of the National Science Foundation. X.S.\ was partially supported by the NSF Career award 2046514, a start-up grant from the University of Pennsylvania, and a fellowship from the Institute for Advanced Study (IAS) at Princeton.  P.Y.\  was partially supported by NSF grant DMS-1712862. 

\section{Preliminaries}


In this paper we work with non-probability measures and extend the terminology of ordinary probability to this setting. For a finite or $\sigma$-finite  measure space $(\Omega, \mathcal{F}, M)$, we say $X$ is a random variable if $X$ is an $\mathcal{F}$-measurable function with its \textit{law} defined via the push-forward measure $M_X=X_*M$. In this case, we say $X$ is \textit{sampled} from $M_X$ and write $M_X[f]$ for $\int f(x)M_X(dx)$. \textit{Weighting} the law of $X$ by $f(X)$ corresponds to working with the measure $d\wt{M}_X$ with Radon-Nikodym derivative $\frac{d\wt{M}_X}{dM_X} = f$, and \textit{conditioning} on some event $E\in\mathcal{F}$ (with $0<M[E]<\infty$) refers to the probability measure $\frac{M[E\cap \cdot]}{M[E]} $  over the space $(E, \mathcal{F}_E)$ with $\mathcal{F}_E = \{A\cap E: A\in\mathcal{F}\}$.  If $M$ is finite, we write $|M| = M(\Omega)$ and $M^\# = \frac{M}{|M|}$ for its normalization. We also fix the notation $|z|_+:=\max\{|z|,1\}$ for $z\in\bbC$.

We also extend the terminology to the setting of more than one random variable sampled from non-probability measures. By saying ``we first sample $X_1$ from $M_1$ and then sample $X_2$ from $M_2$", we refer to a sample $(X_1,X_2)$ from $M_1\times M_2$.  In this setting, weighting the law of $X_2$ by $f(X_2)$ corresponds to working with the measure $d\wt{M}_X$ with Radon-Nikodym derivative $\frac{d\wt{M}_X}{dM_1\times dM_2}(x_1,x_2) = f(x_2)$. In the case where $M_2$ is a probability measure, we say that the marginal law of $X_1$ is $M_1$.

For a 
{M\"{o}bius} 
transform $f:\bbH\to\bbH$ and $s\in\bbR$, if $f(s)=\infty$, then we define $f'(s) = (-\frac{1}{f(w)})'|_{w=s}$. Likewise, if $f(\infty)=s$, then we set $f'(\infty) = ((f^{-1})'(s))^{-1}$. In particular, if $f(z) = a+\frac{\lambda}{x-z}$, then $f'(x) = \lambda^{-1}$ and $f'(\infty) = \lambda$. If $f(z) = a+rz$ with $a\in\bbR,r>0$, then we write $f'(\infty)=r^{-1}$.  These align with the conventions in~\cite{lawler2009partition}.

\subsection{The Gaussian free field and Liouville quantum gravity surfaces}
Let $m$ be the uniform measure on the unit  semicircle  {$\bbH\cap\partial\bbD$}. Define the Dirichlet inner product $\langle f,g\rangle_\nabla = (2\pi)^{-1}\int_X \nabla f\cdot\nabla g $ on the space $\{f\in C^\infty(\bbH):\int_\bbH|\nabla f|^2<\infty; \  \int f(z)m(dz)=0\},$ and let $H(\bbH)$ be the closure of this space {with respect to} the inner product $\langle f,g\rangle_\nabla$. Let $(f_n)_{n\ge1}$ be an orthonormal basis of $H(\bbH)$, and $(\alpha_n)_{n\ge1}$ be a collection of independent standard Gaussian variables. Then the summation
$$h = \sum_{n=1}^\infty \alpha_nf_n$$
a.s.\ converges in the space of distributions on $\bbH$, and $h$ is the \emph{Gaussian Free Field} on $\bbH$ normalized such that $\int h(z)m(dz) = 0$. Let $P_\bbH$ be the law of $h$. See~\cite[Section 4.1.4]{DMS14} for more details.

For $z,w\in\overline{\bbH}$, we define
$$
    G_\bbH(z,w)=-\log|z-w|-\log|z-\bar{w}|+2\log|z|_++2\log|w|_+; \qquad G_\bbH(z,\infty) = 2\log|z|_+.
$$
Then the GFF $h$ is the centered Gaussian field on $\bbH$ with covariance structure $\bbE [h(z)h(w)] = G_\bbH(z,w)$. As pointed out in~\cite[Remark 2.3]{AHS21}, if $\phi = h+f$ where $f$ is a  function continuous everywhere except for finitely many log-singularities, then $\phi$ is a.s.\ in the dual space $H^{-1}(\bbH)$ of $H(\bbH)$.   

Now let $\gamma\in(0,2)$ and $Q=\frac{2}{\gamma}+\frac{\gamma}{2}$. Consider the space of pairs $(D,h)$, where $D\subseteq \mathbb C$ is a planar domain and $h$ is a distribution on $D$ (often some variant of the GFF). For a conformal map $g:D\to\tilde D$ and a generalized function $h$ on $D$, define the generalized function $g\bullet_\gamma h$ on $\tilde{D}$ by setting 
\begin{equation}\label{eq:lqg-changecoord}
    g\bullet_\gamma h:=h\circ g^{-1}+Q\log|(g^{-1})'|.
\end{equation}
 Define the equivalence relation $\sim_\gamma$ {as follows. We say that} $(D, h)\sim_\gamma(\wt{D}, \wt{h})$ if there is a {conformal map $g:D\to\wt{D}$ such that $\tilde h = g\bullet_\gamma h$}.
A \textit{quantum surface} $S$ is an equivalence class of pairs $(D,h)$ under the {equivalence} relation $\sim_\gamma$, and we say {that} $(D,h)$ is an \emph{embedding} of $S$ if $S = (D,h)/\mathord\sim_\gamma$. {Abusing notation, we will sometimes call $(D,h)$ as a quantum surface, and we are then referring to the equivalence class $(D,h)/\sim_\gamma$ that it defines.} Likewise, a \emph{quantum surface with} $k$ \emph{marked points} is an equivalence class of {tuples} 
of the form
	$(D, h, x_1,\dots,x_k)$, where $(D,h)$ is a quantum surface, the points  $x_i\in {\overline{D}}$, and with the further
	requirement that marked points (and their ordering) are preserved by the conformal map $\varphi$ in \eqref{eq:lqg-changecoord}. A \emph{curve-decorated quantum surface} is an equivalence class of tuples $(D, h, \eta_1, ..., \eta_k)$,
	where $(D,h)$ is a quantum surface, $\eta_1, ..., \eta_k$ are curves in $\overline D$, and with the further
	requirement that $\eta$ is preserved by the conformal map $g$ in \eqref{eq:lqg-changecoord}. Similarly, we can
	define a curve-decorated quantum surface with $k$ marked points. Throughout this paper, the curves $\eta_1, ..., \eta_k$ are $\mathrm{SLE}_\kappa$ type curves (which have conformal invariance properties) sampled independently of the surface $(D, h)$.  

 For a $\gamma$-quantum surface $(D, h, z_1, ..., z_m)$, its \textit{quantum area measure} $\mu_h$ is defined by taking the weak limit $\epsilon\to 0$ of $\mu_{h_\epsilon}:=\epsilon^{\frac{\gamma^2}{2}}e^{\gamma h_\epsilon(z)}d^2z$, where $d^2z$ is the Lebesgue area {measure} and $h_\epsilon(z)$ is the circle average of $h$ over $\partial B(z, \epsilon)$. When $D=\mathbb{H}$, we can also define the  \textit{quantum boundary length measure} $\nu_h:=\lim_{\epsilon\to 0}\epsilon^{\frac{\gamma^2}{4}}e^{\frac{\gamma}{2} h_\epsilon(x)}dx$ where $h_\epsilon (x)$ is the average of $h$ over the semicircle $\{x+\epsilon e^{i\theta}:\theta\in(0,\pi)\}$. It has been shown in \cite{DS11, SW16} that all these weak limits are well-defined  for the GFF and its variants we are considering in this paper, {and that} $\mu_{h}$ and $\nu_{h}$ {can} be conformally extended to other domains using the relation $\bullet_\gamma$.

Consider a pair $(D,h)$ where $D$ is now a closed set (not necessarily homeomorphic to a {closed} disk)
such that each component of its interior together with its prime-end boundary is homeomorphic to
the closed disk, and $h$ is only defined as a distribution on each of these components.  We extend the equivalence relation $\sim_\gamma$ described after~\eqref{eq:lqg-changecoord}, such that $g$ is now allowed to be any homeomorphism from $D$ to $\tilde D$ that is conformal on each component of the interior of $D$.  A \emph{beaded quantum surface} $S$ is an equivalence class of pairs $(D,h)$ under the equivalence relation $\sim_\gamma$ as described above, and we say $(D,h)$ is an embedding of $S$ if $S = (D,h)/{\sim_\gamma}$. Beaded quantum surfaces with marked points and curve-decorated beaded quantum surfaces can be defined analogously. 

 As argued in \cite[Section 4.1]{DMS14},  we have the decomposition $H(\mathbb{H}) = H_1(\mathbb{H})\oplus H_2(\mathbb{H})$, where $H_1(\mathbb{H})$ is the subspace of radially symmetric functions, and $H_2(\mathbb{H})$ is the subspace of functions having mean 0 about all semicircles $\{|z|=r,\ \text{Im}\ z>0\}$. As a consequence, for the GFF $h$ sampled from $P_\bbH$, we can decompose $h=h_1+h_2$, where $h_1$ and $h_2$ are independent distributions given by the projection of $h$ onto $H_1({\mathbb{H}})$ and $ H_2({\mathbb{H}})$, respectively. 

 We now turn to the definition of \emph{quantum disks}, which is split in two different cases: \emph{thick quantum disks} and \emph{thin quantum disks}. These surfaces can also be equivalently constructed via methods in Liouville conformal field theory (LCFT) as we shall briefly discuss in the next subsection; see e.g.~\cite{DKRV16, HRV-disk} for these constructions and see \cite{AHS17,cercle2021unit,AHS21} for proofs of equivalence with the surfaces defined above. 
	
	\begin{definition}[Thick quantum disk]\label{def-quantum-disk}
		Fix $\gamma\in(0,2)$ and let $(B_s)_{s\ge0}$ and $(\wt{B}_s)_{s\ge0}$ be independent standard one-dimensional Brownian motions.  Fix a weight parameter $W\ge\frac{\gamma^2}{2}$ and let $\beta = \gamma+ \frac{2-W}{\gamma}\le Q$. Let $\mathbf{c}$ be sampled from the infinite measure $\frac{\gamma}{2}e^{(\beta-Q)c}dc$ on $\bbR$ independently from $(B_s)_{s\ge0}$ and $(\wt{B}_s)_{s\ge0}$.
			Let 	
			\begin{equation*}
				Y_t=\left\{ \begin{array}{rcl} 
					B_{2t}+\beta t+\mathbf{c} & \mbox{for} & t\ge 0,\\
					\wt{B}_{-2t} +(2Q-\beta) t+\mathbf{c} & \mbox{for} & t<0,
				\end{array} 
				\right.
			\end{equation*}
			conditioned on  $B_{2t}-(Q-\beta)t<0$ and  $ \wt{B}_{2t} - (Q-\beta)t<0$ for all $t>0$. Let $h$ be a free boundary  GFF on $\mathbb{H}$ independent of $(Y_t)_{t\in\bbR}$ with projection onto $H_2(\mathbb{H})$ given by $h_2$. Consider the random distribution
			\begin{equation*}
				\psi(\cdot)=X_{-\log|\cdot|} + h_2(\cdot) \, .
		\end{equation*}
		Let $\mathcal{M}_2^{\mathrm{disk}}(W)$ be the infinite measure describing the law of $({\mathbb{H}}, \psi,0,\infty)/\mathord\sim_\gamma $. 
		We call a sample from $\mathcal{M}_2^{\textup{disk}}(W)$ a \emph{quantum disk} of weight $W$ with two marked points.
		
		We call $\nu_\psi((-\infty,0))$ and $\nu_\psi((0,\infty))$ the left and right{, respectively,} quantum   {boundary} length of the quantum disk  $(\mathbb{H}, \psi, 0, \infty)$.
	\end{definition}
	
	When $0<W<\frac{\gamma^2}{2}$, we define the \emph{thin quantum disk} as the concatenation of weight $\gamma^2-W$ thick disks with two marked points as in \cite[Section 2]{AHS23}.

	\begin{definition}[Thin quantum disk]\label{def-thin-disk}
		Fix $\gamma\in(0,2)$. For $W\in(0, \frac{\gamma^2}{2})$, the infinite measure $\mathcal{M}_2^{\textup{disk}}(W)$ on doubly marked beaded quantum surfaces is defined as follows. First sample a random variable $T$ from the infinite measure $(1-\frac{2}{\gamma^2}W)^{-2}\textup{Leb}_{\mathbb{R}_+}$; then sample a Poisson point process $\{(u, \mathcal{D}_u)\}$ from the intensity measure $\mathds{1}_{t\in [0,T]}dt\times \mathcal{M}_2^{\textup{disk}}(\gamma^2-W)$; and finally consider the ordered (according to the order induced by $u$) collection of doubly-marked thick quantum disks $\{\mathcal{D}_u\}$, called a \emph{thin quantum disk} of weight $W$.  {The number $T$ is referred to as the cut point measure of the quantum disk.}
		
		Let $\mathcal{M}_2^{\textup{disk}}(W)$ be the infinite measure describing the law of this ordered collection of doubly-marked quantum disks $\{\mathcal{D}_u\}$.
		The left and right{, respectively,} boundary length of a sample from $\mathcal{M}_2^{\textup{disk}}(W)$ is set to be equal to the sum of the left and right boundary lengths of the quantum disks $\{\mathcal{D}_u\}$.
	\end{definition}

 For $W>0$, one can {disintegrate} the measure $\Md_2(W)$ according to its the quantum length of the left and right boundary arc, i.e., 
 \begin{equation}
     \Md_2(W) = \int_0^\infty\int_0^\infty \Md_2(W;\ell_1,\ell_2)d\ell_1\,d\ell_2,
     \label{eq:disint}
 \end{equation}
where $\Md_2(W;\ell_1,\ell_2)$ is supported on the set of doubly-marked quantum surfaces with left and right boundary arcs having quantum lengths $\ell_1$ and $\ell_2$, respectively. One can also define the measure $\Md_2(W;\ell) := \int_0^\infty \Md_2(W;\ell,\ell')d\ell'$, i.e., the disintegration over the quantum length of the left 
boundary arc. Then we have
\begin{lemma}[Lemma 2.16 and Lemma 2.18 of~\cite{AHS23}]\label{lem-disk-length}
    Let $W\in (0,2+\frac{\gamma^2}{2})$. There exists some constant $c$ depending on $W$ and $\gamma$, such that $$|\Md_2(W;\ell)| = c\ell^{-\frac{2W}{\gamma^2}}.$$
     {The same is true for the disintegration of $\Md_2(W)$ over the right boundary length or the total boundary length.}
\end{lemma}

Finally the weight 2 quantum disk is special in the sense that its two marked points are typical with
respect to the quantum boundary length measure~\cite[Proposition A.8]{DMS14}. Based on this we can define the family of quantum disks marked with multiple quantum typical points.

\begin{definition}\label{def:QD}
    Let $(\mathbb{H},\phi, 0, \infty)$ be  the embedding of a sample from $\Md_2(2)$ as in Definition~\ref{def-quantum-disk}. Let $L=\nu_\phi(\partial\bbH)$, and $\mathrm{QD}$ be the law of $(\bbH,\phi)$ under the reweighted measure $L^{-2}\Md_2(2)$. For $n\ge0$, let  $(\bbH,\phi)$ be a sample {from $\frac{1}{n!}L^n\mathrm{QD}$ and then sample $s_1,...,s_n$ on $\partial\bbH$ according to the probability measure $n!\cdot1_{s_1<...<s_n}\nu_\phi^\#(ds_1)...\nu_\phi^\#(ds_n)$}. Let $\mathrm{QD}_{n}$ be the law of $(\bbH,\phi,s_1,...,s_n)/\sim_\gamma$, and we call a sample from $\QD_{n}$ \emph{a quantum disk with $n$ boundary marked points}.
\end{definition}

\subsection{Liouville conformal field theory on the upper half plane}
Recall that $P_\bbH$ is the law of the free boundary GFF on  $\bbH$ normalized to have average zero on $\partial\bbD\cap\bbH$.
\begin{definition}
Let $(h, \mathbf c)$ be sampled from $P_\bbH \times [e^{-Qc} dc]$ and $\phi = h-2Q\log|z|_+ + \mathbf c$. We call $\phi$ the Liouville field on $\bbH$, and we write
$\LF_\bbH$ for the law of $\phi$.
\end{definition}

\begin{definition}[Liouville field with boundary insertions]\label{def-lf-H-bdry}
Write $\partial\bbH = \bbR\cup\{\infty\}.$	Let $\beta_i\in\mathbb{R}$  and $s_i\in \partial\mathbb{H}$ for $i = 1, ..., m,$ where $m\ge 1$ and all the $s_i$'s are distinct. Also assume $s_i\neq\infty$ for $i\ge 2$. We say  {that} $\phi$ is a \textup{Liouville Field on $\mathbb{H}$ with insertions $\{(\beta_i, s_i)\}_{1\le i\le m}$} if $\phi$ can be produced as follows by  first sampling $(h, \mathbf{c})$ from $C_{\mathbb{H}}^{(\beta_i, s_i)_i}P_\mathbb{H}\times [e^{(\frac{1}{2}\sum_{i=1}^m\beta_i - Q)c}dc]$ with
	$$C_{\mathbb{H}}^{(\beta_i, s_i)_i} =
	\left\{ \begin{array}{rcl} 
\prod_{i=1}^m  |s_i|_+^{-\beta_i(Q-\frac{\beta_i}{2})} \exp(\frac{1}{4}\sum_{j=i+1}^{m}\beta_i\beta_j G_\mathbb{H}(s_i, s_j)) & \mbox{if} & s_1\neq \infty\\
	\prod_{i=2}^m  |s_i|_+^{-\beta_i(Q-\frac{\beta_i}{2}-\frac{\beta_1}{2})}\exp(\frac{1}{4}\sum_{j=i+1}^{m}\beta_i\beta_j G_\mathbb{H}(s_i, s_j)) & \mbox{if} & s_1= \infty
	\end{array} 
	\right.
	 $$
	and then taking
	\begin{equation}\label{eqn-def-lf-H}
\phi(z) = h(z) - 2Q\log|z|_++\frac{1}{2}\sum_{i=1}^m\beta_i G_\mathbb{H}(s_i, z)+\mathbf{c}
	\end{equation}
with the convention $G_\mathbb{H}(\infty, z) = 2\log|z|_+$. We write $\textup{LF}_{\mathbb{H}}^{(\beta_i, s_i)_i}$ for the law of $\phi$.
\end{definition}

The following lemma explains that adding a $\beta$-insertion point at $s\in\partial\mathbb{H}$  is  {equivalent} to weighting the law of Liouville field $\phi$ by $e^{\frac{\beta}{2}\phi(s)}$.

\begin{lemma}[{Lemma 2.8 of~\cite{SY23}}]\label{lm:lf-insertion-bdry}
	For $\beta, s\in\mathbb{R}$ such that $s\notin\{s_1, ..., s_m\}$, in the sense of vague convergence of measures,
	\begin{equation}
	\lim_{\epsilon\to 0}\epsilon^{\frac{\beta^2}{4}}e^{\frac{\beta}{2}\phi_\epsilon(s)}\textup{LF}_{\mathbb{H}}^{(\beta_i, s_i)_i} = \textup{LF}_{\mathbb{H}}^{(\beta_i, s_i)_i, (\beta, s)}.
	\end{equation}
 Similarly, for $a,s_1,...,s_m\in\bbR$, we have
 \begin{equation}
\lim_{R\to\infty}R^{Q\beta-\frac{\beta^2}{4}}e^{\frac{\beta}{2}\phi_R(a)}    \LF_{\bbH}^{(\beta_i,s_i)_i} =  \LF_{\bbH}^{(\beta_i,s_i)_i,(\beta,\infty)}.
\end{equation}
\end{lemma}
The Liouville fields have a nice compatibility with the notion of quantum surfaces. To be more precise, for a measure $M$ on the space of distributions on a domain $D$ and a conformal map $\psi:D\to \tilde{D}$, let $\psi_* M$ be the push-forward of $M$ under the mapping $\phi\mapsto\psi\bullet_\gamma\phi$. Then we have the following conformal covariance of the Liouville field {due to \cite[Proposition 3.7]{HRV-disk} when none of the boundary points are $\infty$; we state a slight generalization by \cite{SY23}.} For $\beta\in\bbR$, we use the shorthand 
$$\Delta_\beta:=\frac{\beta}{2}(Q-\frac{\beta}{2}).$$ 
\begin{lemma}[Lemma 2.9 of~\cite{SY23}]\label{lm:lcft-H-conf}
	Fix $(\beta_i, s_i)\in\mathbb{R}\times\partial\bbH$ 
 for $i=1, ..., m$ with $s_i$'s being distinct.  Suppose $\psi:\mathbb{H}\to\mathbb{H}$ is conformal map. Then $\textup{LF}_{\mathbb{H}} = \psi_*\textup{LF}_{\mathbb{H}}$, and 
	\begin{equation}
	\textup{LF}_{\mathbb{H}}^{(\beta_i, \psi(s_i))_i} = \prod_{i=1}^m|\psi'(s_i)|^{-\Delta_{\beta_i}}\psi_*\textup{LF}_{\mathbb{H}}^{(\beta_i, s_i)_i}.
	\end{equation} 
\end{lemma}

The next lemma shows that sampling points from the Liouville field according to the LQG length measure corresponds to adding $\gamma$-insertions to the field.

\begin{lemma}[{Lemma 2.13 of~\cite{SY23}}]\label{lm:gamma-insertion}
    Let $m\ge 2,n\ge1$ and $(\beta_i,s_i)\in \bbR\times\partial \bbH$ with $\infty\ge s_1>s_2>...>s_m>-\infty$. Let $1\le k\le m-1$ and $g$ be a non-negative measurable function supported on $[s_{k+1},s_k]$. Then as measures we have the identity
\begin{equation}\label{eq:gamma-insertion}
\begin{split}
    &\int\int_{x_1,...,x_n\in[{s_{k+1}},{s_k}]}\, g(x_1,...,x_n)\nu_\phi(dx_1)...\nu_\phi(dx_n) \LF_\bbH^{(\beta_i,s_i)_i}(d\phi) \\&= \int_{x_1,...,x_n\in[{s_{k+1}},{s_k}]}\int\, \LF_\bbH^{(\beta_i,s_i)_i,(\gamma,x_1),...,(\gamma,x_n)}(d\phi)\, g(x_1,...,x_n)dx_1...dx_n.
    \end{split}
\end{equation}
\end{lemma}


Next we recall the relations between marked quantum disks and Liouville fields. The statements in~\cite{AHS21} are involving Liouville fields on the strip $\mathcal{S}:=\bbR\times(0,\pi)$, yet we can use the map $z\mapsto e^z$ to transfer to the upper half plane.

\begin{definition}\label{def:three-pointed-disk}
	Let $W>0$. First sample  a quantum disk from $\Md_2(W)$ and weight its law by the quantum length of its left boundary arc. Then sample a marked point on the left boundary according to the probability measure proportional to the LQG length. We denote the law of the triply marked quantum surface  by $\mathcal{M}_{2, \bullet}^{\textup{disk}}(W)$, where this newly added point is referred as the third marked point.
\end{definition}

\begin{proposition}[Proposition 2.18 of \cite{AHS21}]\label{prop:m2dot}
	For $W>\frac{\gamma^2}{2}$ and $\beta = \gamma+\frac{2-W}{\gamma}$, let $\phi$ be sampled from $\frac{\gamma}{2(Q-\beta)^2}\textup{LF}_{\mathbb{H}}^{(\beta,\infty), (\beta,0), (\gamma, 1)}$. Then $(\bbH, \phi, 0, \infty, 1)/{\sim_\gamma}$ has the same law as $\mathcal{M}_{2, \bullet}^{\textup{disk}}(W)$.
\end{proposition}

The proposition above gives rise to the quantum disks with general third insertion points, which could be defined via three-pointed Liouville fields.

\begin{definition}\label{def:m2dot-alpha}
	Fix $W>\frac{\gamma^2}{2},\beta=\gamma+\frac{2-W}{\gamma}$ and let $\beta_3\in\mathbb{R}$. Set $\mathcal{M}_{2, \bullet}^{\textup{disk}}(W;\beta_3)$ to be the law of $(\bbH, \phi, 0,\infty, 1)/\sim_\gamma$ with $\phi$ sampled from $\frac{\gamma}{2(Q-\beta)^2}\textup{LF}_{\bbH}^{(\beta,0),(\beta,\infty), (\beta_3, 1)}$. We call the boundary arc between the two $\beta$-singularities with (resp.\ not containing) the $\beta_3$-singularity the marked (resp.\ unmarked) boundary arc.
\end{definition}

Next we turn to the $W\in(0,\frac{\gamma^2}{2})$ case. Recall the following fact from~\cite{AHS23}.
\begin{lemma}[Proposition 4.4 of \cite{AHS23}]\label{lem:m2dot-thin}
    For $W\in(0,\frac{\gamma^2}{2})$ we have
    \begin{equation*}
        \Md_{2,\bullet}(W) = (1-\frac{2}{\gamma^2}W)^2\Md_2(W)\times\Md_{2,\bullet}(\gamma^2-W)\times\Md_2(W),
    \end{equation*}
    where the right hand side is the infinite measure on ordered collection of quantum surfaces obtained
by concatenating samples from the three measures.
\end{lemma}

\begin{definition}\label{def:m2dot-thin}
    Let $W\in(0,\frac{\gamma^2}{2})$ and $\beta\in\bbR$. Given a sample $(S_1, S_2, S_3)$ from
    \begin{equation*}
        (1-\frac{2}{\gamma^2}W)^2\Md_2(W)\times\Md_{2,\bullet}(\gamma^2-W;\beta)\times\Md_2(W),
    \end{equation*}
    let $S$ be their concatenation in the sense of Lemma~\ref{lem:m2dot-thin} with $\beta$ in place of $\gamma$. We define the infinite measure $\Md_{2,\bullet}(W;\beta)$ to be the law of $S$. 
\end{definition}

\subsection{The Schramm-Loewner evolution}
Fix $\kappa>0$. We start with the $\SLE_\kappa$ process on the upper half plane $\bbH$. Let $(B_t)_{t\ge0}$ be the standard Brownian motion. The $\SLE_\kappa$ is the probability measure on continuously growing {curves $\eta$ in $\overline{\bbH}$}, {whose mapping out function $(g_t)_{t\ge0}$ (i.e., the unique conformal transformation from the unbounded component of $\mathbb{H}\backslash \eta([0,t])$ to $\mathbb{H}$ such that $\lim_{|z|\to\infty}|g_t(z)-z|=0$) can be described by}
\begin{equation}\label{eq:def-sle}
g_t(z) = z+\int_0^t \frac{2}{g_s(z)-W_s}ds, \ z\in\mathbb{H},
\end{equation} 
where $W_t=\sqrt{\kappa}B_t$ is the Loewner driving function. For weights $\rho^-,\rho^+>-2$, the $\SLE_\kappa(\rho^-;\rho^+)$ process is the probability measure on {curves $\eta$  in $\overline{\bbH}$} 
 {such that \eqref{eq:def-sle} is still satisfied}, except that
the Loewner driving function $(W_t)_{t\ge 0}$ is now characterized by 
\begin{equation}\label{eq:def-sle-rho}
\begin{split}
&W_t = \sqrt{\kappa}B_t+\sum_{q\in\{+,-\}}\int_0^t \frac{\rho^{q}}{W_s-V_s^{q}}ds; \\
& V_t^{\pm} = 0^{\pm}+\int_0^t \frac{2}{V_s^{\pm}-W_s}ds, \ q\in\{L,R\}.
\end{split}
\end{equation}
It has been proved in \cite{MS16a} that the SLE$_\kappa(\rho^-;\rho^+)$ process a.s.\ exists, is unique and generates a continuous curve. The curve is simple for $\kappa\in[0,4]$, has self-touchings for $\kappa\in(4,8)$ and is space-filling when $\kappa\ge8$. 

The $\SLE_\kappa$, as a probability measure, can be defined on other domains by conformal maps. To be more precise, let $\mu_\bbH(0,\infty)$ be the $\SLE_\kappa$ on $\bbH$ from 0 to $\infty$, $D$ be a simply connected domain, and $f:\bbH\to D$ be a conformal map with $f(0)=x,f(\infty)=y$. Then we can define a probability measure $\mu_D({x,y})^\# = f\circ\mu_\bbH(0,\infty)$. Let $$b=\frac{6-\kappa}{2\kappa}$$ be the \emph{boundary scaling exponent}, and recall that for $x,y\in\partial D$ such that $\partial D$ is smooth near $x,y$, the boundary Poisson kernel is defined by $H_D(x,y)=\varphi'(x)\varphi'(y)(\varphi(x)-\varphi(y))^{-2}$ where $\varphi:D\to\bbH$ is a conformal map. Then as in~\cite{lawler2009partition}, one can define the $\SLE_\kappa$ in $(D,x,y)$ as a non-probability measure by setting $\mu_D({x,y}) = H_D(x,y)^b\cdot\mu_D({x,y})^\#$, which satisfies the conformal covariance
$$f\circ\mu_{D}(x,y) = |f'(x)|^b|f'(y)|^b\mu_{f(D)}(f(x),f(y))  $$
for any conformal map $f:D\to f(D)$.

\subsection{The multiple SLE and its partition function}\label{subsec:pre-msle}
In this section, we review the background and some basic properties of the multiple chordal SLE as established in e.g.~\cite{lawler2009partition,peltola2019global,peltola2019toward}. In particular, we shall focus on the probabilistic construction of the partition function and multiple SLE for $\kappa\in(0,6]$ in~\cite{wu2020hypergeometric,peltola2019global}, which will be the base of our results for  {the} $\kappa\in(6,8)$ regime. 

Recall that for $N=1$, the partition function for the $\SLE_\kappa$ in $(D,x,y)$ is $H_D(x,y)^b$. For $N=2$, it is shown that Theorem~\ref{thm:existence-uniqueness} holds for $\kappa\in(0,8)$~\cite{miller2018connection}. Moreover, the solutions to (PDE) and (COV) have explicit expressions, which have the form of the hypergeometric functions. This gives the partition function $\mathcal{Z}_\alpha(\bbH;x_1,x_2,x_3,x_4)$, which could be extended to other simply connected domains via~\eqref{eq:partition-conformal-conf}.

For $N\ge2$, the multiple SLE for $\kappa\in(0,6]$ ($\kappa\in(0,8)$ for $N=2$) and the partition function can be defined via the following induction. Let $(D;x_1,...,x_{2N})$ be a polygon and $\alpha\in\LP_N$.  We fix the following notations:
\begin{itemize}
    \item Let $\{i,j\}\in\alpha$ with $i<j$ be a link and let $\hat\alpha\in\LP_{N-1}$ be the link pattern obtained by removing $\{i,j\}$ from $\alpha$;
    \item  For a continuous curve $\eta$ in $\ol D$, let $\mathcal{E}_\eta$ be the event where $\eta$ does not partition $D$ into components where some variables corresponding to a link in  {$\hat\alpha$} would belong to different components;
    \item On the event $\mathcal{E}_\eta$, let $\hat D_{\eta}$ be the union of connected components of $D\backslash\eta$ with some of the points $\{x_1,...,x_{2N}\}$ on the boundary, i.e., 
    \begin{equation}\label{eq:def-msle-a}
        \hat D_{\eta} = \mathop{\bigcup_{
            \tilde D \text{ c.c. of } D\backslash\eta \atop  {\partial\tilde D\cap\{x_1,...,x_{2N}\}\backslash\{x_i,x_j\}}\neq\emptyset} \tilde D
        }
    \end{equation}
    \item On the event $\mathcal{E}_\eta$, define  
    \begin{equation}\label{eq:def-msle-b}
        \mathcal{Z}_{\hat\alpha}(\hat D_{\eta};x_1,...,x_{i-1},x_{i+1},...,x_{j-1},x_{j+1},...,x_{2N}) = \mathop{\prod_{
            \tilde D \text{ c.c. of } D\backslash\eta \atop  {\partial\tilde D\cap\{x_1,...,x_{2N}\}\backslash\{x_i,x_j\} {\neq}\emptyset}}} \mathcal{Z}_{\alpha_{\tilde D}}(\tilde D;...)
    \end{equation}
    where for each $\tilde D$, the ellipses ``..." stand for those variables among $\{x_1,..., x_{2N}\} \backslash \{x_i, x_j\}$ which belong to $\partial \tilde D$, and $\alpha_{\tilde D}$ stands for the sub-link patterns of $\hat\alpha$ associated to the components $\tilde D\subset D\backslash\eta$.
\end{itemize}
The measure $\mSLE_{\kappa,\alpha}(D;x_1,...,x_{2N})$ is defined as follows:
\begin{enumerate}[(i)]
    \item Sample $\eta_1$ as an $\SLE_\kappa$ in $(D;x_i,x_j)$ and weight its law by
    $$ \mathds{1}_{{\mathcal{E}_{\eta_1}}}H_D(x_i,x_j)^b\times \mathcal{Z}_{\hat\alpha}(\hat D_{\eta_1};...);$$
    \item Sample $(\eta_2,...,\eta_N)$ from the probability measure
    \begin{equation*}
        \mathop{\prod_{
            \tilde D \text{ c.c. of } D\backslash\eta \atop  {\partial\tilde D\cap\{x_1,...,x_{2N}\}}\backslash\{x_i,x_j\}=\emptyset}} \mSLE_{\kappa,\alpha_{\tilde D}}(\tilde D;...)^\#;
    \end{equation*}
    \item Output $(\eta_1,...,\eta_N)$ and let $\mSLE_{\kappa,\alpha}(D;x_1,...,x_{2N})$  {denote} its law, while $\mathcal{Z}_\alpha(D;x_1,...,x_{2N})$ is the  {total measure} of $\mSLE_{\kappa,\alpha}(D;x_1,...,x_{2N})$. 
\end{enumerate}
We remark that the above induction is well-defined since in each of $\mSLE_{\kappa,\alpha_{\tilde D}}(\tilde D;...)^\#$ the number of marked points is {strictly} less than $2N$, and the resulting measure $\mathrm{mSLE}_{\kappa, \alpha}(D; x_1, \dots, x_{2N})$ does not depend on the choice of the link $\{i,j\}\in\alpha$~\cite[Proposition B.1]{peltola2019toward}. It is then shown in~\cite{wu2020hypergeometric,peltola2019global} that the partition function above is well-defined and satisfies the following power law bound
\begin{equation}\label{eq:powerlawbound}
    \mathcal{Z}_\alpha(\bbH;x_1,...,x_{2N}) \le \prod_{\{i,j\}\in\alpha}|x_j-x_i|^{-2b}.
\end{equation}
Moreover, it is easy to verify from definition that the probability measure $\mSLE_{\kappa,\alpha}(D;x_1,...,x_{2N})^\#$ satisfies the resampling property in Theorem~\ref{thm:existence-uniqueness}.

Finally, we comment that for $\kappa\in(6,8)$, if  the partition function $\mathcal{Z}_{\hat\alpha}$ is finite for any $\hat\alpha\in \bigsqcup_{k=1}^{N-1}\LP_k$, then the measure $\mSLE_{\kappa,\alpha}(D;x_1,...,x_{2N})$ is well-defined for any $\alpha\in\LP_N$ (which is not necessarily finite at this moment) and does not depend on the choice of the link $\{i,j\}\in\alpha$. Indeed, using the  symmetry in the exchange of the two curves in the 2-$\SLE_\kappa$, 
this follows from exactly the same argument as~\cite[Proposition B.1]{peltola2019toward}, where we may first sample the 2-$\SLE_\kappa$ and weight by the product of the  partition functions in the subdomains cut out by the two curves, and then sample the rest of the curves from the multiple SLE probability measure.

\section{Conformal welding of forested quantum disks for $\kappa\in(4,8)$}
In this section, we work on the conformal welding of forested quantum surfaces. We start with the definition of forested lines and forested quantum disks, and use a pinching argument to prove Theorem~\ref{thm:fd-disk}. Then we will give decomposition theorems for forested quantum disks with marked points, and show the resampling property of weight $\gamma^2-2$ forested quantum disks. 

\subsection{Generalized quantum surfaces}\label{subsec:gqd}

\begin{figure}
    \centering
    \begin{tabular}{cc} 
		\includegraphics[scale=0.49]{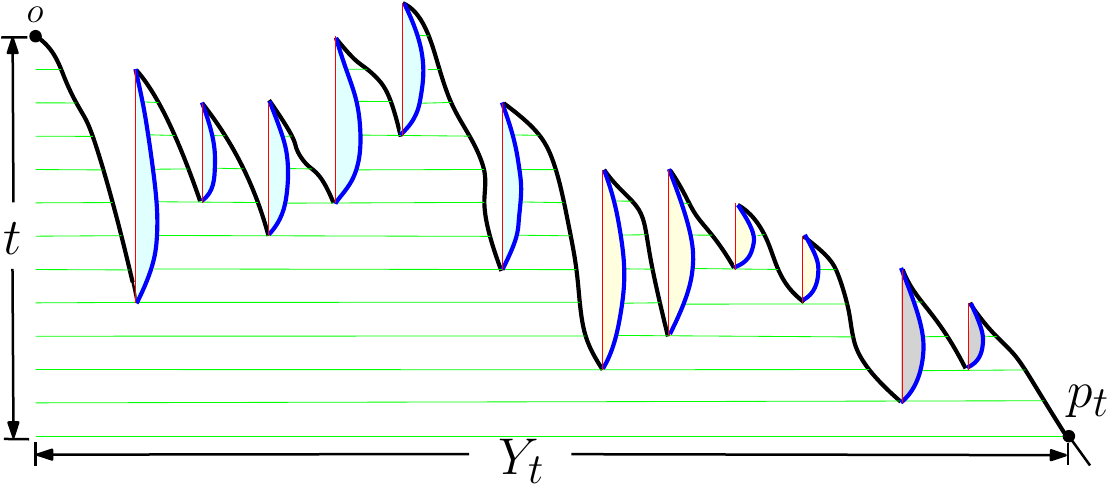}
		&
	   \includegraphics[scale=0.65]{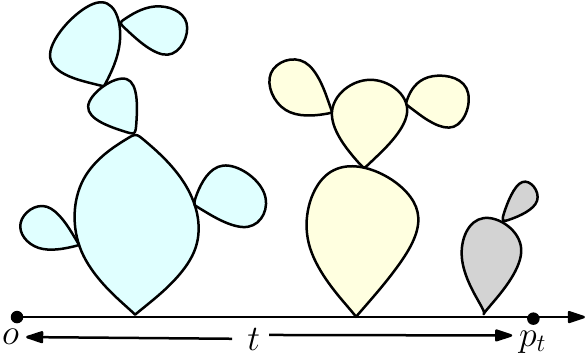}
	\end{tabular}
 \caption{\textbf{Left:} The graph of the L\'{e}vy process $(X_t)_{t>0}$ with only upward jumps. We draw the blue curves for each of the jump, and identify the points that are on the same green horizontal line.   \textbf{Right:} The L\'{e}vy tree of disks obtained from the left panel. For each topological disk we assign a quantum disk $\QD$ conditioned on having the same boundary length as the size of the jump, with the points on the red line in the left panel shrinked to a single point.  The quantum length of the line segment between the root $o$ and the point $p_t$ is $t$, while the segment along the forested boundary between $o$ and $p_t$ has generalized quantum length $Y_t = \inf\{s>0:X_s\le -t\}$.  }\label{fig:forestline-def}
 \end{figure}

We start by recalling the notion of forested lines from~\cite{DMS14}. Recall that for the measure $\QD$ on quantum surfaces {from Definition~\ref{def:QD}}, we can define the disintegration $\QD = \int_0^\infty \QD(\ell)d\ell$ over its boundary length, and  {$|\QD(\ell)| = c\ell^{-\frac{4}{\gamma^2}-2}$}.  Let $(X_t)_{t\ge0}$ be a stable L\'{e}vy process of index $\frac{\kappa}{4}\in(1,2)$ with only upward jumps, {so} $X_t\overset{d}{=}t^{\frac{4}{\kappa}}X_1$ for any $t>0$. On the graph of $X$, we draw two curves {for each time $t$ at which $X$ jumps: One curve is a straight vertical line segment connecting the points $(t,X_t)$ and $(t,X_{t^-})$, while the other curve is to the right of this line segment connecting its two end-points.  {The vertical line segment will eventually be collapsed to a single point, creating a disk whose boundary is the second curve.} The precise form of the second curve does not matter as long as it intersects each horizontal line at most once, it stays
 below the graph of $X_t$, and it does not intersect the vertical line segment except at its end-points.} 
Then we draw horizontal lines and identify the points lying on the same segment which does not go above the graph of $(X_t)_{t>0}$. We also identify   $(t,X_t)$ with $(t,X_{t^-})$ for every $t$ with $X_t\neq X_{t^-}$. This gives a tree of topological disks. For each jump of size $L$, we then independently sample a quantum disk of boundary
length $L$ from  $\QD(L)^\#$ and topologically identify the boundary of each quantum disk with its corresponding  loop  {in a clockwise length-preserving way, with rotation chosen uniformly at random (i.e., for any boundary point of the quantum disk, the point on the loop it is identified with is chosen randomly)}.  The unique point corresponding to 
 {$(0,0)$ }
on the graph of $X$ is called the \emph{root}. The closure of the collection of the points on the boundaries of the quantum disks is referred as the \emph{forested boundary arc}, while the set of  {the points corresponding to the running infimum of $(X_t)_{t\ge0}$} 
is called the \emph{line boundary arc}. Since $X$ only has positive jumps, the quantum disks are lying on the same side of the line boundary arc, whose points correspond to the running infimum of $X$.  
See Figure~\ref{fig:forestline-def} for an illustration.


As discussed in~\cite{DMS14}, forested lines are beaded quantum surfaces. Since the sum of the jumps which
occur in any non-empty open interval of time of an $\alpha$-stable L\'{e}vy process is infinite, the sum of the lengths of the loops in a looptree is infinite. As such, any path connecting two points  on the boundary of different quantum disks has infinite quantum length, and the forested boundary arc shall {instead} be parametrized by \emph{generalized quantum length}. 
\begin{definition}[Forested line]\label{def:forested-line}
For $\gamma \in (\sqrt2,2)$, let  {$(X_s)_{s\geq 0}$} be a stable L\'evy process of index $\frac{4}{\gamma^2}>1$ with only positive jumps {satisfying $X_0=0$ a.s.}. For $t>0$, let $Y_t=\inf\{s>0:X_s\le -t\}$, and {fix the multiplicative constant of $X$ such that} $\bbE[e^{-Y_1}] = e^{-1}$. Define the forested line as {described above}.

 The line boundary arc is parametrized by quantum length.  {That is, for the point $p_t$ on the line boundary arc corresponding to the time $Y_t$, we assign $t$ units of quantum length for the line segment from the root to $p_t$.} The forested boundary arc is parametrized by generalized quantum length; that is, the length of the corresponding interval of $(X_t)_{t\geq 0}$.  
 For a point $p_t$ on the line boundary arc with LQG distance $t$ to the root, the segment of the forested boundary arc between $p_t$ and the root has generalized quantum length $Y_t$.  
\end{definition}


Next we gather some lemmas  {about} the law of the generalized quantum length of forested lines.

\begin{lemma}\label{lem:law-forested-line}
For $\gamma \in (\sqrt2, 2)$, sample a forested line and for $t > 0$, let $p_t$ be the point on the line boundary arc at quantum length $t$ from the root, and let $Y_t$ be the generalized quantum length of the forested arc from $p_t$ to the root. 
Then $(Y_t)_{t\geq 0}$ is a stable subordinator of index $\frac{\gamma^2}4 \in (0,1)$. 
\end{lemma}
\begin{proof}
Recall the definition of a forested line in terms of a L\'evy process $(X_s)_{s \geq 0}$ with index $\frac4{\gamma^2}>1$. Since record minima of $(X_s)$ correspond to points on the line, we have
\[Y_t = \inf \{ s>0 \: : \: X_s \leq -t \}.\]
 {The process} $(Y_t){_{t\geq 0}}$ is increasing, 
 {the random variable $Y_t$ is infinitely divisible for each fixed $t>0$ (e.g. by the Markov property of $X$ we have $Y_t \stackrel d= Y_{t/2} + \wt Y_{t/2}$ where $\wt Y$ is an independent copy of $Y$)},
and satisfies a scaling relation $(Y_{kt}){_{t\geq 0}} \stackrel d= (k^{\frac4{\gamma^2}} Y_t){_{t\geq 0}}$ (inherited from $(X_{k^{4/\gamma^2}s}){_{s\geq 0}} \stackrel d= (k X_s){_{s\geq 0}}$). Therefore $(Y_t){_{t\geq 0}}$ is a stable subordinator with index $\frac{\gamma^2}4$. 
\end{proof}

\begin{lemma}[L\'evy process moments]\label{lem-levy-moment}
For $p < \frac{\gamma^2}{4}$,
\[\bbE[Y_1^p] = \frac4{\gamma^2} 
\frac{\Gamma(-\frac4{\gamma^2}p)}{\Gamma(-p)}. \]
Conversely, for $p \geq \frac{\gamma^2}4$, we have $\bbE[Y_1^p] = \infty$.
\end{lemma}
\begin{proof}
  By our normalization,   $\bbE[e^{-\lambda Y_1}] = e^{- \lambda^{\gamma^2/4}}$ for all $\lambda >0$.
For $p < 0$, we have
\begin{equation*}
    \begin{split}
\Gamma(-p) \bbE[Y_1^{p}] = \bbE\left[ \int_0^\infty e^{-\lambda Y_1} \lambda^{-p-1} d\lambda \right] &= \int_0^\infty\bbE\left[  e^{-\lambda Y_1}  \right]\lambda^{-p-1} d\lambda \\&= \int_0^\infty  e^{-\lambda^{\gamma^2/4}} \lambda^{-p-1} d\lambda = \Gamma(\frac4{\gamma^2} p) \frac4{\gamma^2}.
\end{split}
\end{equation*}

Similarly, for $p \in (0, \frac{\gamma^2}4)$,
\begin{equation*}
\begin{split}
    \Gamma(1-p) \bbE[Y_1^p] = \bbE \left[ \int_0^\infty Y_1 e^{-\lambda Y_1} \lambda^{-p} d\lambda\right] &= \int_0^\infty \bbE \left[  Y_1 e^{-\lambda Y_1} \right] \lambda^{-p} d\lambda \\&=\frac{\gamma^2}4 \int_0^\infty \lambda^{\gamma^2/4 - 1-p} e^{- \lambda^{\gamma^2/4}} d\lambda = \Gamma(1-\frac4{\gamma^2} p),
    \end{split}
\end{equation*}
and applying the identity $z\Gamma(z) = \Gamma(z+1)$ yields the desired formula. 
Finally, for $p \geq \frac{\gamma^2}4$, the integral in the previous equation does not converge. 
\end{proof}

Since this paper focuses on finite volume surfaces, we define the following \emph{truncation} operation. For $t>0$ and a forested line $\cL$ with root $o$, mark the point $p_t$ on the line boundary arc with quantum length $t$ from $o$. By \emph{truncation of $\cL$ at quantum length $t$}, we refer to the surface $\cL_t$ which is the union of the line boundary arc and the quantum disks on the forested boundary arc between $o$ and $p_t$. In other words, $\cL_t$ is the surface generated by $(X_s)_{0\le s\le Y_t}$ in the same way as Definition~\ref{def:forested-line}, and the generalized quantum length of the forested boundary arc of $\cL_t$ is $Y_t$. 
{The beaded quantum surface $\cL_t$ is called a forested line segment.}

\begin{definition}\label{def:line-segment}
    Fix $\gamma\in(\sqrt{2},2)$. Define $\mathcal{M}_2^\mathrm{f.l.}$ as the law of the surface obtained by first sampling $\mathbf{t}\sim \mathrm{Leb}_{\bbR_+}$ and truncating an independent forested line at quantum length $\mathbf{t}$. 
\end{definition}

\begin{lemma}[Law of forested segment length]\label{lem:fs-len}
Fix $q \in \bbR$. 
Suppose we sample $\mathbf t \sim 1_{t > 0} t^{-q} dt$ and independently sample a forested line $\cL$.  
For $q<2$,  the law of  
$Y_{\mathbf t}$ is $\ C_q  \cdot1_{L>0} L^{-\frac{\gamma^2}4q + \frac{\gamma^2}4 - 1} dL.$ where $C_q := \frac{\gamma^2}4 \bbE[Y_1^{\frac{\gamma^2}4 (q-1)}]<\infty$. If $q\geq2$, then for any $0<a<b$, the event $\{Y_{\mathbf t}\in[a,b]\}$ has infinite measure. 
\end{lemma}
\begin{proof}
Write $\beta = \frac4{\gamma^2}$ and $M$ for the reference measure describing the law of $((X_t)_{t\ge0},\mathbf t)$. 
Then for $0<a<b$,
\begin{equation*}
\begin{split}
    M[Y_{\mathbf t}\in [a, b]] &= \int_0^\infty t^{-q} \bbP[Y_t \in [a, b]]dt 
= \int_0^\infty t^{-q}\int  1_{t^\beta Y_1 \in [a,b] }\,d\bbP dt\\&
= \frac{1}{\beta}\int \int_a^b s^{\frac{1-q}{\beta}-1}Y_1^{\frac{q-1}{\beta}}\, dsd\bbP = \frac{1}{\beta}\bbE\big[Y_1^{\frac{q-1}{\beta}}\big] \int_a^b s^{\frac{1-q}{\beta}-1}ds.
\end{split}
\end{equation*}
where we applied Fubini's theorem and the change of variable $s = t^\beta Y_1$. Then the claim follows as the finiteness/infiniteness of $\bbE[Y_1^p]<\infty$ is given by Lemma~\ref{lem-levy-moment}.
\end{proof}

Now we introduce the formal definition of generalized quantum surfaces. Let $n\ge1$,  {$m\ge0$,} and $(D,\phi, {w_1,...,w_m},z_1,...,z_n)$ be an embedding of a connected possibly beaded quantum surface $S$ of finite volume, with $z_1,...,z_n\in\partial D$  ordered clockwise  {and $w_1,...,w_m\in D$}. 
    We sample independent forested lines $\cL^1,...,\cL^n$, truncate them such that their quantum lengths match the length of boundary segments $[z_1,z_2],...,[z_n,z_1]$ and glue them to $\partial D$ correspondingly. 
    Let $S^f$ be the output beaded quantum surface. 

 

\begin{definition}\label{def:f.s.}
   We call   a beaded quantum surface $S^f$ as above a (finite volume) \emph{generalized quantum surface}, 
   {and   $S^f$ together with its marked points a \emph{generalized quantum surface with marked points}.} We call this procedure \emph{foresting the boundary} of $S$, and  
   say $S$  \emph{is the spine of}  $S^f$. 
\end{definition}

When $n\ge2$, the spine $S$ can be recovered from $S^f$ in the following way. Let $(D^f, \phi^f, z_1,..., z_n)$ be an embedding of $S^f$. Let $D\subset D^f$ be the domain with $z_1,...,z_n$ on the boundary, such that the boundary arc of $D$ from $z_k$ to $z_{k+1}$ consists of the ``shortest" (i.e., intersection of all possible) clockwise path from $z_k$ to $z_{k+1}$ within $\partial D^f$. Then $S = (D, \phi^f, z_1,..., z_n)/{\sim_\gamma}$.  {Note that in this case, for each thin quantum disk of the spine, there {is} 
a marked point on each of its ends that is not a cut point of $\partial D^f$.}

\begin{definition}\label{def:f.d.}
For any $W > 0$, write ${\mathcal M}_2^{\mathrm{f.d.}}(W)$ for the infinite measure on generalized quantum surfaces obtained by first taking a quantum disk from ${\mathcal M}_2^{\mathrm{disk}}(W)$, then foresting its two boundary arcs. A sample from $\Mfd_2(W)$ is called a forested quantum disk of weight $W$. 
\end{definition}

 {Our definition of the generalized quantum surface is based on the forested lines in~\cite{DMS14}. The generalized quantum surfaces are also defined in~\cite{MSW22CLE,HL22CLE} using a different formalism. In Section~\ref{subsec:fdp.p.p.} we will see that their definitions are equivalent with ours;}   {see Remark~\ref{rmk:QD-def}.}


\begin{lemma}\label{lem:fd-exponent}
For $W \in (0, \gamma^2)$, the law of the 
generalized quantum length of the left (resp.\ whole) boundary of a forested quantum disk from ${\mathcal M}_2^{\mathrm{f.d.}}(W)$ is given by $1_{L>0} cL^{-1-\frac W2 + \frac{\gamma^2}4}dL$ for some constant $c$. 

When $W \geq \gamma^2$, the mass of forested quantum disks with generalized quantum length in $[1,2]$ is infinite.
\end{lemma}
\begin{proof}
This follows immediately from Lemma~\ref{lem:fs-len} {and Lemma~\ref{lem-disk-length}}.
\end{proof}

{Recall the disintegration \eqref{eq:disint} of the quantum disk measure. By disintegrating over the values of $Y_t$, we can similarly define a disintegration of $\mathcal M_2^\mathrm{f.l.}$: 
\[
\mathcal{M}_2^\mathrm{f.l.} = \int_{\bbR_+^2}\mathcal{M}_2^\mathrm{f.l.}(t;\ell)\,dt\,d\ell.
\]
where $\mathcal{M}_2^\mathrm{f.l.}(t;\ell)$ is the measure on forested line segments with quantum length $t$ for the line boundary arc and generalized quantum length $\ell$ for the forested boundary arc. We write $\mathcal{M}_2^\mathrm{f.l.}(\ell):=\int_0^\infty\mathcal{M}_2^\mathrm{f.l.}(t;\ell)dt$, i.e., the law of forested line segments whose forested boundary arc has generalized quantum length $\ell$.  A similar disintegration holds for the forested quantum disk, namely,
\begin{equation}
     \mathcal{M}^\mathrm{f.d.}_2(W) = \int_0^\infty\int_0^\infty \mathcal{M}^\mathrm{f.d.}_2(W;\ell_1,\ell_2)\,d\ell_1\,d\ell_2.
 \end{equation}
Indeed, this follows by defining the measure $\mathcal{M}^\mathrm{f.d.}_2(W;\ell_1,\ell_2)$ via
$$\int_{\bbR_+^2}\mathcal{M}_2^\mathrm{f.l.}(t_1;\ell_1)\times \Md_2(W;t_1,t_2)\times \mathcal{M}_2^\mathrm{f.l.}(t_2;\ell_2)\,dt_1dt_2.$$ 
}

\subsection{Forested line as a Poisson point process of forested disks}\label{subsec:fdp.p.p.}

In this section we study the Poissonian structure of forested lines and prove the resampling property of the weight $\gamma^2-2$ forested quantum disks. The results are implicitly stated in~\cite[Section 2.2]{MSW22CLE} and can be proved following ideas from~\cite{CK13looptree}. Here we include the precise statements and proofs for completeness. 

\begin{definition}\label{def:GQD1}
Let $\GQD_{2} := {\mathcal M}_2^{\mathrm{f.d.}}(\gamma^2-2)$ be the infinite measure on generalized quantum surfaces, and 
let $\GQD_{1}$ denote the corresponding measure when we forget one of the marked points and unweight by the generalized quantum length of the forested boundary. 
\end{definition}

\begin{lemma}\label{lem:excursion-jump}
Let $(\wt X_t)_{t \geq 0}$ be a stable L\'evy process with index $\beta = \frac4{\gamma^2} \in (1,2)$ with only downward jumps, and let $\wt I_{t^-} := \inf_{s < t} \wt X_s$ be its (infinitesimally lagged) infimum process. The process $(\wt X_t - \wt I_{t^-})$ is comprised of an ordered collection of excursions, each starting from 0 and ending at the first time it takes a negative value (this corresponds to  times when $\wt X_t$ jumps past its previous infimum). Let $\wt M$ be the excursion measure of $(\wt X_t - \wt I_{t^-})$.
For an excursion $e \sim \wt M$, let $u = -e(\tau)$ and $v = e(\tau^-)$, where $\tau$ is the duration of $e$. Then the $\wt M$-law of $(u,v)$ is $c (u+v)^{-1-\beta}$ for some constant $c>0$. 
\end{lemma}
\begin{proof}
   The lemma follows from \cite[Example 7]{levy-overshoot}. There, they fix $x>0$ and consider an index-$\beta$ stable L\'evy $(X_t)_{t \geq 0}$ process started at 0 and run until it first time it hits or exceeds $x$. We specialize to the case where it only has upward jumps, so the positivity parameter is $\rho = 1-\frac1\beta$ (see e.g.~\cite[Chapter VIII, above Lemma 1]{bertoin-levy}), so \cite[Example 7]{levy-overshoot} gives
\[\bbP[X_{\tau_x^+} - x \in du, x - X_{\tau_x^+-} \in dv, x - \ol X_{\tau_x^+-} \in dy] = \mathrm{const}\ \cdot\ \mathds{1}_{u,v,y>0} \cdot \frac{(x-y)^{\beta - 2}}{(u+v)^{1+\beta}} \, dy \, du \, dv,\]
where $X_{t-} :={} \lim_{ s \uparrow t} X_s$, $\ol X_{t-} := \sup_{s < t} X_s$,
and $\tau_x^+$ is the first time $t$ that $X_t \geq x$. Integrating out $y$ then sending $x$ to 0 yields the analogous result for excursions with only upward jumps. Flipping the sign then yields the lemma. 
\end{proof}

\begin{figure}
    \centering
    \begin{tabular}{cc} 
		\includegraphics[scale=0.49]{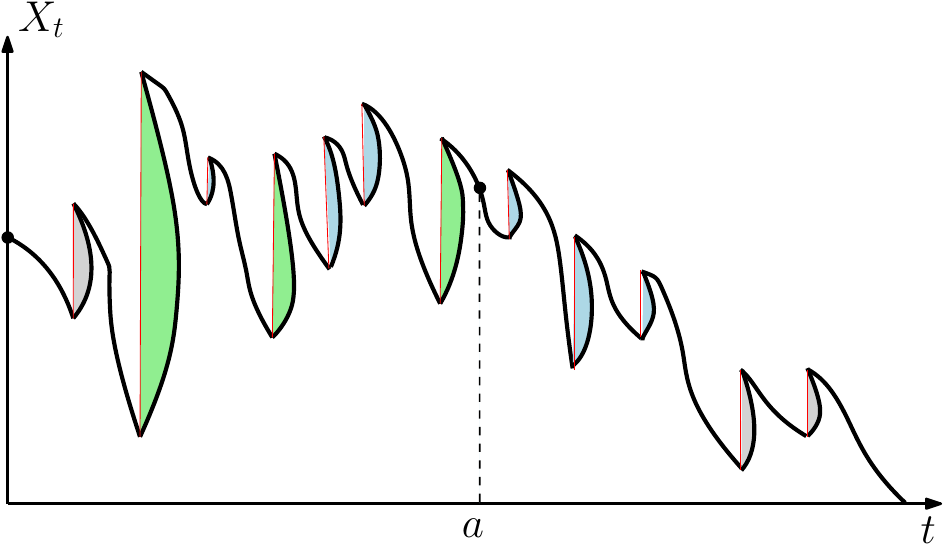}
		&
	   \includegraphics[scale=0.65]{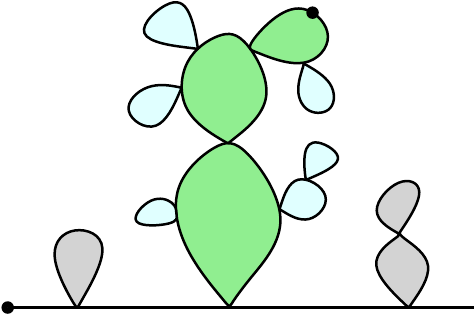}
	\end{tabular}
 \caption{An illustration of the proof of Proposition~\ref{prop:typical}. \textbf{Left:} We sample a time  {$a>0$} from the Lebesgue measure.    \textbf{Right:} The L\'{e}vy tree of disks obtained from the left panel with the marked point. The collection of green disks, which we shall prove to have law $c\Md_2(\gamma^2-2)$, correspond to the jumps of $(X_{u-t})_{t \in [0,u]}$ hitting running infimum. }\label{fig:forestline-poisson}
 \end{figure}

\begin{proposition}\label{prop:typical}
Sample a forested line, and consider the collection of pairs $(u, \mathcal D_u^f)$ such that $\mathcal D_u^f$ is a generalized quantum surface attached to the line boundary arc (with its root defined to be the attachment point) and $u$ is the quantum length from the root of the forested line  to the root of $\mathcal D_u^f$. Then the law of this collection is a Poisson point process with intensity measure $c\mathrm{Leb}_{\bbR_+} \times \GQD_{1}$ for some constant $c>0$.
\end{proposition}

 {We note that the above proposition immediately gives the reversibility of forested lines in the following sense. For any fixed $t>0$, consider a sample $\cL_t\sim \mathcal{M}_2^\mathrm{f.l.}(t)$, i.e.,  the truncation of a forested line at quantum length $\ell$. Then  $\cL_t$, when viewed from the reverse direction,  has the same law $\mathcal{M}_2^\mathrm{f.l.}(t)$.
}

\begin{proof}
Let $\beta = \frac4{\gamma^2}$, and recall that the forested line is defined by a stable L\'evy process $(X_t)_{t \geq 0}$ (with exponent $\beta$) with only upward jumps. Let $I_t := \inf_{s \leq t} X_s$, so the process $(X_t - I_t)$ decomposes as an ordered collection of excursions. The process $-I_t$ is a local time at 0 for $(X_t - I_t)$ \cite[Chapter VII, Theorem 1]{bertoin-levy}. Thus, if $M$ is the excursion measure of $(X_t - I_t)$, the set of pairs $(u, e_u)$, where $e_u$ is an excursion of $(X_t - I_t)$ and  {$u$ is the value of $-I_t$ during this excursion}, is a Poisson point process with intensity measure $\mathrm{Leb}_{\bbR_+} \times M$. 
For a sample $e \sim M$, as in Section~\ref{subsec:gqd}, one can construct a generalized quantum surface by sampling an independent quantum disk for each jump of $e$; let $M'$ be its law. 
By the construction of the forested line, the collection $(u, \cD_u^f)$ is a Poisson point process with intensity measure $\mathrm{Leb}_{\bbR_+} \times M'$, so we need to show that $M' = c \mathrm{GQD}_1$ for some constant $c$.

Sample a point from the (infinite) generalized quantum length on the forested line and let $\cD^f$ be the generalized quantum surface attached to the line boundary arc which contains this point. Note that $\cD^f$ has two marked points: the sampled point, and the root of $\cD^f$. Sampling a point from generalized quantum length measure corresponds to picking a time $a \in \bbR_+$ from the Lebesgue measure and looking at the excursion $e_u$  of 
$(X_t - I_t)$ containing $a$. The law of $e_u$ is $M$ weighted by excursion duration, so the law of $\cD^f$ is $M'$ weighted by generalized quantum length and with an additional marked point sampled according to the probability measure proportional to generalized quantum length. Therefore it suffices to show that the law of $\cD^f$ is $c{\mathcal M}_2^{\mathrm{f.d.}}(\gamma^2-2)$ for some constant $c$. 


{Let $(\wt X_t)_{t \in [0,a]}$ be the modification of $(X_{a-t} - X_a)_{t \in [0,a]}$ which is right continuous with left limits, so the marginal law of $a$ is Lebesgue measure on $\mathbb R_+$ and the conditional law of $(\wt X_t)_{t \in [0,a]}$ given $a$ is an index-$\beta$ stable L\'evy process with only downward jumps stopped at time $a$.} 
Let $\Lambda = \{t>0:\inf_{s\le t}\wt X_s = \wt X_t \}$. Using $(\wt X_t)_{t \in [0,a]}$
 we define a beaded quantum surface $\tilde\cD^f$ as follows:
\begin{enumerate}[(i)]
    \item For each time $\tau\in\Lambda$, 
    let $u_\tau = \wt X_{\tau^-} - \inf_{s<\tau}\wt X_s$, $v_\tau =\inf_{s<\tau}\wt X_s-\wt X_\tau $, and $\ell_\tau = u_\tau+v_\tau$. Independently sample a quantum disk $\tilde\cD_\tau$ from $\QD(\ell_\tau)^\#$, uniformly choose a marked point $x_\tau \in \partial \tilde \cD_\tau$ from the probability measure proportional to the boundary quantum length measure, then mark the point $y_\tau\in \partial \tilde \cD_\tau$ which is $v_\tau$ units of quantum length clockwise of $x_\tau$. 
    \item For each $\tilde\cD_\tau$, forest its left and right boundaries 
    to get $\tilde \cD_\tau^f$.
    \item Concatenate the surfaces $\tilde\cD_\tau^f$ according to the ordering induced by $\tau$ at the points $x_\tau,y_\tau$ to get $\tilde \cD^f$.
\end{enumerate}

We first show that $\cD^f \stackrel d= \tilde \cD^f$, then check the law of $\tilde \cD^f$ is $c {\mathcal M}_2^{\mathrm{f.d.}}(\gamma^2-2)$ to complete the proof.

We decompose $\wt X$ in terms of its excursions above its infimum. Let $F^\mathrm{left}$ and $F^\mathrm{right}$ be the forested parts to the left and right of the marked point, so $\cD^f = F^\mathrm{left} \cup \cD \cup F^\mathrm{right}$, and $\cD$ is a doubly marked beaded quantum surface. Then as in Figure~\ref{fig:forestline-poisson}, note that the disks $(\cD_\tau)_{\tau\in\Lambda}$ of $\cD$ correspond to the downward jumps hitting the running infimum of $(\wt X_t)_{t\in[0,a]}$. Each $\cD_\tau$ carries two marked boundary points; the left and right boundary lengths of $\cD_\tau$ are equal to $u_\tau$ and $v_\tau$, and $\cD$ is obtained by concatenating all the $\cD_\tau$'s together. Furthermore,  for each excursion with starting time $\sigma$ and ending time $\tau$, the left boundary side of $\cD_\tau$ is glued to a segment of forested line generated by $(X_s)_{a-\tau<s<a-\sigma}$ (which is determined by $(\wt X_s)_{\sigma<s<\tau})$. Moreover, by the Markov property of $(X_t)_{t\ge0}$, the segment $F^\mathrm{right}$ is independent of $(F^{\mathrm{left}}, \cD)$ given the right boundary length of $\cD$, and has the law of a forested segment of that length.   From this construction we see that the surfaces $\cD^f$ and $\tilde \cD^f$ have the same law.

On the other hand, by Lemma~\ref{lem:excursion-jump}, the joint law of $(u_\tau, v_\tau)$ is given by $c\mathds{1}_{u,v>0}(u+v)^{-1-\beta}dudv$, whereas $|\Md_2(2;\ell,r)| = c\mathds{1}_{\ell,r>0}(\ell+r)^{-1-\beta}$~\cite[Proposition 7.7]{AHS23}.  {Recall that $(\wt X_t)_{0\leq t\leq a}$ can be  generated by 
 {sampling $a \sim \mathrm{Leb}_{\bbR_+}$ and an independent index-$\beta$ stable L\'evy process (still denoted by $(\wt X_t)_{t\geq 0}$) with only downward jumps.} Consider the excursion decomposition $(\wt u, \wt e_{\wt u})$ of $(\wt X_t - \inf_{s\leq t}\wt X_s)$ as in Lemma~\ref{lem:excursion-jump}. Then the collection of quantum disks $\{\tilde D_\tau:\tau\in\Lambda\}$ are generated by the excursions $(\wt u,\wt e_{\wt u})$ before time $a$. On the other hand, using an identical argument as in~\cite[Proposition 4.2]{AHS23}, these excursions can be generated by first sampling $\wt T\sim \mathrm{Leb}_{\bbR_+}$ and then considering a Poisson point process  {with intensity measure} $\mathrm{Leb}_{[0,\wt T]}\times\wt M$. In particular, by Definition~\ref{def-thin-disk}, given $\wt T$, $\{\tilde D_\tau:\tau\in\Lambda\}$ has the same law as $c\Md_2(\gamma^2-2)$  {conditioned on having} cut point measure $\wt T$.} 
It then follows that $\tilde \cD^f$ has the same law as   a sample from  $c {\mathcal M}_2^{\mathrm{f.d.}}(\gamma^2-2)$. This concludes the proof of the proposition.
\end{proof}

\begin{remark}\label{rmk:QD-def}
    Our definition of generalized quantum surfaces is via gluing independent forest lines as in~\cite{DMS14}. In~\cite{MSW22CLE} as well as~\cite[Definition 5.8]{HL22CLE}, the generalized quantum disk is defined in a similar manner to Definition~\ref{def:forested-line} via \emph{excursions} of the stable L\'evy process of index $\frac{4}{\gamma^2}$, where there is one marked point on the boundary corresponding to the starting point of the excursion. 
    The generalized quantum half plane there is defined via gluing  {a Poisson point process} of generalized quantum disks to weight $\gamma^2-2$ quantum wedges. It is proved implicitly in~\cite[Section 2.2]{MSW22CLE} that their definition and our definition 
    are equivalent. By decomposing $(X_t-I_t)_{t\ge0}$ into excursions as in the proof above, it immediately follows from  Proposition~\ref{prop:typical} that the generalized quantum disks defined in~\cite{MSW22CLE} and~\cite[Definition 5.8]{HL22CLE} equals (some constant times) $\GQD_1$.\footnote{  {Note that in \cite{HL22CLE} the generalized quantum disk is represented as an excursion along with a simply connected LQG surface for each jump, but the authors could equivalently have viewed the generalized quantum disk as a beaded quantum surface as we do in this paper.}} 
\end{remark}

The following is {a corollary of Proposition~\ref{prop:typical} via} the  argument of~\cite[Lemma 4.1]{AHS23}.
\begin{lemma}\label{lem:forest-decomposition}
    Fix $\ell>0$. The following three procedures yield the same infinite measure on generalized quantum surfaces.
    \begin{itemize}
        \item Sample a forested line $\cL$ and truncate it to have quantum length $\ell$. Then
take a point from the generalized quantum length measure on the forested boundary arc (this induces a weighting by the forested arc length).
        \item Sample a forested line $\cL$ and truncate it to have quantum length $\ell$. then independently take $(u,\cD^f)$ from $\mathrm{Leb}_{[0,\ell]}\times \Mfd_2(\gamma^2-2)$. Insert $\cD^f$ into $\cL$ at the point with distance $u$ from the root.
        \item  Take $(u,\cD^f)$ from $\mathrm{Leb}_{[0,\ell]}\times \Mfd_2(\gamma^2-2)$, then independently sample two forested lines and truncate them to have quantum lengths $u$ and $\ell-u$. Concatenate the three surfaces.
    \end{itemize}
\end{lemma}

\begin{definition}
    Let $W>0$. First sample  a quantum disk from $\Mfd_2(W)$ and weight its law by the generalized quantum length of its left boundary arc. Then sample a marked point on the left boundary according to the probability measure proportional to the generalized quantum length. We denote the law of the triply marked quantum surface  by $\mathcal{M}_{2, \bullet}^{\textup{f.d.}}(W)$. 
\end{definition}

The following is a quick consequence of Lemma~\ref{lem:forest-decomposition} by recalling that the third marked point for a sample from $\Md_{2,\bullet}(W)$ {(Definition~\ref{def:three-pointed-disk})} is sampled from the quantum length measure.
\begin{lemma}\label{lem:m2fd}
    Let $W>0$. A sample from $\mathcal{M}_{2, \bullet}^{\textup{f.d.}}(W)$ can be produced by
    \begin{enumerate}[(i)]
        \item Sampling $(\cD,\cD^f)$ from $\Md_{2,\bullet}(W)\times \Mfd_2(\gamma^2-2)$, and concatenating $\cD^f$ to $\cD$ at the third marked point of $\cD$;
        \item Foresting the boundary of $\cD$.
    \end{enumerate}
\end{lemma}

Now we are ready to prove the following analog of~\cite[Proposition A.8]{DMS14} in the setting of forested quantum disks.
\begin{proposition}\label{prop:def-gqd}
    Let $(D^f,\phi,x,y)$ be an embedding of a sample from $\Mfd_2(\gamma^2-2)$. If we independently sample $x',y'\in\partial D^f$ from the probability measure proportional to the generalized quantum length measure on $\partial D^f$, then the law of the quantum surface $(D^f,\phi,x',y')$ is still $\Mfd_2(\gamma^2-2)$.
\end{proposition}
\begin{proof}
    Following Lemma~\ref{lem:m2fd} and Definition~\ref{def:m2dot-thin}, if we weight the law of $(D^f,\phi,x,y)$ by the generalized quantum length of $\partial D^f$ and sample $x'\in\partial D^f$ from the probability measure proportional to the generalized quantum length measure, then the law $(D^f,\phi,x,y,x')$ is a constant multiple of $\Mfd_2(\gamma^2-2)\times\Mfd_2(\gamma^2-2)\times\Mfd_2(\gamma^2-2)\times\QD_{3}^f$, i.e., we first sample a surface from $\QD_{3}$ and concatenate three independent forested quantum disks from $\Mfd_2(\gamma^2-2)$, and then glue truncated forested lines to the boundary arcs of $\QD_{3}$. In particular, this implies that if we forget the marked point $x$, the quantum surface $(D^f,\phi,x',y)$ has the law of $\Mfd_2(\gamma^2-2)$ weighted by its total forested boundary length. Applying the previous argument once more yields the proposition.
\end{proof}

 \begin{definition}\label{def:GQD}
     Let $(D^f,\phi,x,y)$ be an embedding of a sample from $\Mfd_2(\gamma^2-2)$, and $L$ be the generalized quantum length of the forested boundary. Let $\mathrm{GQD}$ be the law of $(D^f,\phi)$ under the reweighted measure $L^{-2}\Mfd_2(\gamma^2-2)$. For $n\ge1$, let  $(D^f,\phi)$ be a sample from {$\frac{1}{(n-1)!}L^n\mathrm{GQD}$} and then sample $s_1,...,s_n$ on $\partial D^f$ ordered clockwise according to the probability measure proportional to the generalized quantum length measure. Let $\mathrm{GQD}_{n}$ be the law of $(D^f,\phi,s_1,...,s_n)/{\sim_\gamma}$, and we call a sample from $\GQD_{n}$ \emph{a generalized quantum disk with $n$ boundary marked points}.
 \end{definition}

We have the following law on the boundary lengths of generalized quantum disks.
\begin{proposition}\label{prop:gqd-bdry-length}
    Let $n\ge 1$. For a sample  from $\GQD_{n}$, let $\ell_1,...,\ell_n$ be the generalized quantum lengths of the $n$ boundary segments. Then for some constant $c$, the law of $(\ell_1,...,\ell_n)$ is $$c\,1_{\ell_1,...,\ell_n>0}(\ell_1+...+\ell_n)^{-\frac{\gamma^2}{4}-1}d\ell_1...d\ell_n.$$
\end{proposition}
\begin{proof}
The $n=1$ case is immediate from Lemma~\ref{lem:fd-exponent} and Definition~\ref{def:GQD1}.
    For $n=2$, the claim follows from the same argument as~\cite[Proposition 7.8]{AHS23} via Lemma~\ref{lem:fd-exponent} and Proposition~\ref{prop:def-gqd}. Assume the statement has been proved for $n$. If $(D^f,\phi,s_1,...,s_n)$ is an embedding of a sample $\GQD_{n}$ with $s_1,...,s_n$ ordered clockwise, then as we weight the law of $(D^f,\phi,s_1,...,s_n)$ by the generalized quantum length $\ell_n$ of the forested boundary segment between $s_n$ and $s_1$, and sample $s_{n+1}$ on this segment from the probability measure proportional to the generalized quantum length measure,  the law of $(D^f,\phi,s_1,...,s_{n+1})$ is  $\GQD_{n+1}$. 
    Let $\ell_n'$ be the generalized quantum length of the boundary arc between $s_n$ and $s_{n+1}$. Then the joint law of $(\ell_1,...,\ell_n,\ell_n')$ is $c1_{\ell_1,...,\ell_n>0;0<\ell_n'<\ell_n}(\ell_1+...+\ell_n)^{-\frac{\gamma^2}{4}-1}d\ell_1...d\ell_nd\ell_n'.$ Therefore the claim follows by setting $\ell_{n+1}' = \ell_n-\ell_n'$ and a change of variables.
\end{proof}

\subsection{Welding of forested quantum surfaces}

In this section we prove Theorem~\ref{thm:fd-disk}. The idea is to start with the quantum wedge counterpart in~\cite[Theorem 1.15]{DMS14} and use a pinching argument. We start with the definition of thin quantum wedges.

\begin{definition}[Thin quantum wedge]\label{def:thin-wedge}
    Fix $W\in (0,\frac{\gamma^2}{2})$ and sample a Poisson point process $\{(u,\mathcal{D}_u)\}$ from the measure $\mathrm{Leb}_{\bbR_+}\otimes\Md_2(\gamma^2-W)$.  The weight $W$ quantum wedge is the infinite beaded surface obtained by concatenating the $\mathcal{D}_u$ according to the ordering induced by $u$. We write $\mathcal{M}_2^{\textup{wedge}}(W)$ for the probability measure on weight $W$ quantum wedges.
\end{definition}


The following is from~\cite[Theorem 1.15]{DMS14}.

\begin{theorem}\label{thm:wedge-forest}
    Let $\kappa\in(4,8)$ and $\gamma = \frac{4}{\sqrt{\kappa}}$. Consider a quantum wedge $\mathcal{W}$ of weight $W=2-\frac{\gamma^2}{2}$, and let $\eta$ be the concatenation of an independent $\SLE_\kappa(\frac{\kappa}{2}-4;\frac{\kappa}{2}-4)$ curve on each bead of $\cW$. Then $\eta$
 divides $\cW$ into two independent forested lines $\cL_-,\cL_+$, whose forested boundaries are identified with one
another according to the generalized quantum length. Moreover, $({\mathcal{W}},\eta)$ is measurable  {with respect to} $({\mathcal{L}^-},{\mathcal{L}^+})$.
\end{theorem}
Recall that 
$\SLE_\kappa$ has double points but no triple points for $\kappa\in(4,8)$ (see e.g.~\cite{MW17}). To view $\cL_\pm$ as beaded quantum surfaces, suppose $p$ is a double point and $\eta$ visits $p$ at times $0<t<t'$; pick $\eps \in (0, t'-t)$. If $\eta(t')$ hits $\eta((0,t+\e))$ on the left (resp.\ right) side, then we view $\eta(t)$ and $\eta(t')$ as two different points on the boundary of $\cL_+$ (resp.\ $\cL_-$) staying infinitesimally close to each other. Under this point of view, $\cL_-$ and $\cL_+$ can be embedded on the closure of the union of the bubbles cut out by $\eta$ lying on the left and right side of $\eta$, respectively.

In light of Theorem~\ref{thm:wedge-forest}, we can extend the notion of generalized quantum length to $\SLE_\kappa$ curves on $\gamma$-LQG surfaces with $\kappa=\frac{16}{\gamma^2}\in(4,8)$.  In fact, the generalized quantum length of $\eta$ as above agrees with the \emph{quantum natural parametrization} of $\SLE_\kappa$  up to a constant. The quantum natural parametrization is roughly the quantum version  of natural time parameterization of SLE~\cite{LS11pa,lawlerparamet}; see ~\cite[Section 6.5]{DMS14} for more details. 

Recall the law of the perimeter of a sample from $\GQD_{1}$ has law $c\ell^{-\frac{\gamma^2}{4}-1}d\ell$. For $0<d<2$,  a squared Bessel process $(Z_t)_{t\ge0}$ of dimension $0<d<2$ can be constructed by sampling a Poisson point process $\{(u,e_u)\}$ from $1_{u>0}du\times\mathcal{E}$ and concatenating the $e_u$'s according to the ordering of $u$, where $\mathcal{E}$ is the It\^{o} excursion measure. Moreover, the law of the duration of an  excursion from $\mathcal{E}$ is $c\ell^{\frac{d}{2}-2}d\ell$ {for some constant $c$}.  Using  the Poissonian description of forested lines from Proposition~\ref{prop:typical}, we obtain the following:
\begin{lemma}\label{lem:forest-line-bessel}
    Let $Z$ be a squared Bessel process of dimension $d=2-\frac{\gamma^2}{2}$. Consider its It\^{o} decomposition $\{(u,e_u)\}$ of excursions over 0. Fix a line $\cL$ and parametrize it by quantum length. For each $u$, independently sample a generalized quantum disk from $\GQD_{1}$ conditioned on having {generalized quantum boundary length equal to} the duration of $e_u$, and attach it to $\cL$ at distance $u$ to the root. Then $\cL$ is a forested line.
\end{lemma}

\begin{lemma}\label{lem:length-bead}
In the setting of Theorem~\ref{thm:wedge-forest}, let $\{ (u, \cD_u)\}$ be the decomposition of $\cW$ as in Definition~\ref{def:thin-wedge}, and let $\ell_u$ be the generalized quantum length of $\eta$ in $\cD_u$. Then for some constant $c  {>0}$, $\{(u, \ell_u) \}$ is a Poisson point process on $\bbR_+ \times \bbR_+$ with intensity measure $c 1_{u>0} du \times 1_{\ell > 0} \ell^{-\frac{\gamma^2}2} d\ell$.
\end{lemma}
\begin{proof}
      Let $Z^\pm$ be the associated squared Bessel processes of $\cL_\pm$ as in Lemma~\ref{lem:forest-line-bessel}. Then the cut points of $\cW$ corresponds to common zeros of $Z^-$ and $Z^+$, while $Z:=Z^-+Z^+$ is a squared Bessel process of dimension $4-\gamma^2<2$ (see e.g.~\cite[Section XI.1]{revuzyorbook}). This finishes the proof.
\end{proof}

\begin{proposition}\label{prop:besselcondition}
     Let $d\in(0,2)$ and $\mathcal{E}$ be the It\^{o} excursion measure of the $d$-dimensional squared Bessel process above 0. The laws of $Z$ and $\wt Z$ sampled from the following two procedures agree:
     \begin{enumerate}[(i)]
         \item First sample $T$ from $1_{t>0}dt$. Then sample a Poisson point process $\{(u,e_u)\}$ from $1_{0<u<T}du\times\mathcal{E}$ and concatenate the $e_u$'s according to the ordering of $u$ to generate  {the process $(Z_t)_{0\le t\le T}$}.
         \item First sample $\wt L$ from $1_{\ell>0}\ell^{-\frac{d}{2}}d\ell$. Then sample a $d$-dimensional squared Bessel bridge $(\wt Z_t)_{0\le t\le \wt L}$, i.e., a $d$-dimensional squared Bessel process $(\wt Z_t)_{t\ge0}$ from 0 conditioned on $\wt Z_{\wt L}=0$.
     \end{enumerate}
 \end{proposition}

{
\begin{proof}
Let $\mathcal M$ denote the law of the process $Z$ from (i). 
    By a direct computation similar to~\cite[Lemma 2.18]{AHS23}, the law of the duration of $Z$ is  $1_{\ell>0}\ell^{-\frac{d}{2}}d\ell$. Thus, writing $\mathcal M(\ell)^\#$ to denote the law of a sample from $\mathcal M$ conditioned to have duration $\ell$, we have $\mathcal M = \int_0^\infty \mathcal M(\ell)^\# \ell^{-\frac d2}\, d\ell$. 
    Therefore, it suffices to prove that for each $\ell>0$ the probability measure $\mathcal M(\ell)^\#$ agrees with the law of the duration $\ell$ $d$-dimensional squared Bessel bridge. 

 Let $(\mathfrak B_t)_{t \geq 0}$ be a $d$-dimensional squared Bessel process started from 0. Let $\tau  = \sup \{ t \leq \ell : \mathfrak B_t = 0\}$ and $\sigma = \inf \{ t \geq \ell : \mathfrak B_t = 0\}$, and define the event $E_\eps = \{\sigma < \ell + \e \}$. We will show that the law of $(\mathfrak B_t)_{0 \leq t \leq \ell}$ conditioned on $E_\eps$ converges to $\mathcal M(\ell)^\#$  as $\eps \to 0$, and it also converges in law to the duration $\ell$ squared Bessel bridge. This would complete the proof.

    Since the law of the length of an excursion from $\mathcal{E}$ is $c\ell^{\frac{d}{2}-2}d\ell$, the proofs in~\cite[Section 4]{AHS23} can be carried line by line if we replace the measure $\Md_2(\gamma^2-W)$ there by $\mathcal{E}$, $\mathcal{M}_2^\mathrm{wedge}(W)$ there by the law of $d$-dimensional squared Bessel process from 0, the quantum length measure by $\mathrm{Leb}_{\bbR_+}$ and $\Md_2(W)$ there $\mathcal{M}$, where $W=\frac{\gamma^2d}{4}$. In particular, it follows from the identical proofs that:
    \begin{itemize}
        \item By~\cite[Corollary 4.3]{AHS23}, conditioned on the event $E_\e$, the law of $((\mathfrak B_t)_{0\le t\le \tau}, (\mathfrak B_t)_{\tau\le t\le\sigma})$ agrees with $\mathcal{M}\times\mathcal{E}$ conditioned on $E_\e'$ 
        , where $E_\e'$ is the event that the durations $(x,y)$ of  a pair of processes satisfies $x<\ell<x+y<\ell+\e$. Then the joint law of $(x,y)$ is $$\frac{1}{Z_{\ell,\e}}1_{0<x<\ell<x+y<\ell+\e}x^{-\frac{d}{2}}y^{\frac{d}{2}-2}\,dxdy$$ where $Z_{\ell,\e} = \iint_{0<x<\ell<x+y<\ell+\e} x^{-\frac{d}{2}}y^{\frac{d}{2}-2}\,dxdy$. Thus, conditioned on $\tau$, the conditional law of $(\mathfrak B_t)_{0 \leq t \leq \tau}$ is $\mathcal M(\tau)^\#$. Moreover, conditioned on $E_\eps$, we have $\tau \to \ell$ in probability as $\eps \to 0$, since for any $\ell>\delta>0$ we have  $$ \mathbb P[\tau < \ell - \delta \mid E_\e] = \frac{\int_0^{\ell-\delta}\int_{\ell-x}^{\ell-x+\e} x^{-\frac{d}{2}}y^{\frac{d}{2}-2}\,dydx }{\int_0^{\ell}\int_{\ell-x}^{\ell-x+\e} x^{-\frac{d}{2}}y^{\frac{d}{2}-2}\,dydx} \xrightarrow{\e\to0} 0.$$
         \item We have the weak convergence $\mathcal M(\ell')^\# \to \mathcal M(\ell)^\#$ as $\ell' \to \ell$, with respect to the topology on function space corresponding to uniform convergence on compact subsets of $(0, \ell)$. This follows from \cite[Corollary 4.7]{AHS23}.
    \end{itemize}
    Combining the above two inputs, we conclude that the law of  $(\mathfrak B_t)_{0 \leq t \leq \ell}$ conditioned on $E_\e$ converges in law as $\e \to 0$ to $\mathcal M(\ell)^\#$.

        On the other hand, given $E_\eps$ we have $\mathfrak B_\ell \xrightarrow{p} 0$ as $\e \to 0$. This is immediate from the transition densities and hitting times of Bessel processes given in~\cite[Section 3]{lawler2018notes}. Next, given $E_\eps$ and $\mathfrak B_\ell$, the process $(\mathfrak B_t)_{0\le t\le\ell}$ is a squared Bessel bridge from 0 to $\mathfrak B_\ell$; indeed this is trivially true when conditioning only on $\mathfrak B_\ell$, and since $\mathfrak B$ is a Markov process the further conditioning on $E_\eps$ does not affect the law of $(\mathfrak B_t)_{0 \leq t \leq \ell}$. Finally, the law of the duration $\ell$ squared Bessel bridge from $0$ to $b$ converges as $b \to 0$ to the law of the duration $\ell$ squared Bessel bridge from $0$ to $0$ \cite[Section XI.3]{revuzyorbook}. We conclude that  $(\mathfrak B_t)_{0 \leq t \leq \ell}$ conditioned on $E_\eps$ converges in law as $\eps \to 0$ to the duration $\ell$ squared Bessel bridge.
\end{proof}
}

The following corollary is immediate from Definition~\ref{def:line-segment}, Proposition~\ref{prop:typical} and Proposition~\ref{prop:besselcondition}.
\begin{corollary}\label{cor:segment-bessel}
    For $\ell>0$, a sample from $\mathcal{M}_2^\mathrm{f.l.}(\ell)^\#$ can be generated by a squared Bessel bridge $(Z_t)_{0\le t\le \ell}$ of length $\ell$ from the same method as in Lemma~\ref{lem:forest-line-bessel}.
\end{corollary}

\begin{proposition}\label{prop:weld:segment}
     Let $\kappa\in(4,8)$ and $\gamma = \frac{4}{\sqrt{\kappa}}$. Consider a quantum disk $\mathcal{D}$ of weight $W=2-\frac{\gamma^2}{2}$, and let $\tilde\eta$ be the concatenation of an independent $\SLE_\kappa(\frac{\kappa}{2}-4;\frac{\kappa}{2}-4)$ curve on each bead of $\cD$. Then for some constant $c$, $\tilde\eta$
 divides $\cD$ into two   forested lines segments $\wt\cL_-,\wt\cL_+$, whose law is
 \begin{equation}\label{eq:weld:segment}
      c\int_0^\infty \mathcal{M}_2^\mathrm{f.l.}(\ell)\times\mathcal{M}_2^\mathrm{f.l.}(\ell)d\ell.
 \end{equation}
{Moreover, $\wt\cL_\pm$ a.s.\ uniquely determine $(\cD,\tilde\eta)$ in the sense that $(\cD,\tilde\eta)$ is measurable with respect to the $\sigma$-algebra generated by $\wt\cL_\pm$.}
\end{proposition}
\begin{proof}
   We start with the setting of Theorem~\ref{thm:wedge-forest}, and let $Z^\pm$ be the squared Bessel process of dimension $d=2-\frac{\gamma^2}{2}$ associated with $\cL_\pm$ as in the proof of Lemma~\ref{lem:length-bead}. Then following Lemma~\ref{lem:length-bead} and Definition~\ref{def:thin-wedge}, the curve decorated surface $(\cW,\eta)$ can be generated by 
   \begin{enumerate}[(i)]
       \item Sample a squared Bessel process $(Z_t)_{t\ge0}$ of dimension $4-\gamma^2$, and decompose it into excursions $\{(u,e_u)\}$;
       \item For each excursion $(u,e_u)$, sample a curve-decorated surface $(\cD_u,\eta_u)$ from $\Md_2(\gamma^2-W)\otimes\SLE_\kappa(\frac{\kappa}{2}-4;\frac{\kappa}{2}-4)$ conditioned on the interface length being the excursion length of $e_u$;
       \item Concatenate all the $(\cD_u,\eta_u)$'s together according to the ordering induced by $u$.
   \end{enumerate}
   Moreover, $Z$ is coupled with $Z^\pm$ such that $Z=Z^++Z^-$. By Definition~\ref{def-thin-disk}, if we sample $T$ from $1_{t>0}dt$ and concatenate the $(\cD_u,\eta_u)$'s with $u<T$, then the output surface has law $\Md_2(W)\otimes\SLE_\kappa(\frac{\kappa}{2}-4;\frac{\kappa}{2}-4)$. Therefore here and below we assume $(\cD,\tilde\eta)$ and $(\cW,\eta)$ are coupled as above. On the other hand, for $\ell>0$, by Proposition~\ref{prop:besselcondition}, conditioning on the interface length being $\ell$ is the same as conditioning on $Z_\ell^+=Z^-_\ell=Z_\ell=0$. Then $(Z_t^+)_{0\le t\le \ell}$ and $(Z_t^-)_{0\le t\le \ell}$ are independent squared Bessel bridges of length $\ell$. Indeed, if we let $\sigma_\ell=\inf\{t>\ell:Z_t=0\}$ and condition on $E_\e:=\{\ell<\sigma_\ell<\ell+\e\}$, for any $\delta>0$, as $\e\to0$, the joint law of  $(Z_t^+)_{0\le t\le \ell-\delta}$ and $(Z_t^-)_{0\le t\le \ell-\delta}$ converges to that of independent squared Bessel bridges of lengths $\ell$ truncated at time $\ell-\delta$. This is because, conditioned on $Z_\ell^\pm$, $(Z_t^\pm)_{0\le t\le \ell}$ are squared Bessel bridges from 0 to $Z_\ell^\pm$, while $Z_\ell^-+Z_\ell^+$ converges to 0 in probability.  Since given $(Z_t^+)_{0\le t\le \ell}$ and $(Z_t^-)_{0\le t\le \ell}$, the surfaces to the left and right of $\eta$ are collections of generalized quantum disks with perimeters matching the duration of the excursions of $Z^\pm$ above 0, it then follows from Corollary~\ref{cor:segment-bessel} that $\wt\cL_-$ and $\wt\cL_+$ are independent forested line segments of forested boundary length $\ell$. Moreover, following the same argument as~\cite[Lemma 2.18]{AHS23}, by Lemma~\ref{lem:length-bead}, the law of the generalized quantum length of $\tilde\eta$ is $c\ell^{\frac{\gamma^2}{2}-2}d\ell$. On the other hand, by taking $q=0$ in Lemma~\ref{lem:fs-len}, the law of the generalized quantum length of the interface in~\eqref{eq:weld:segment} is  $c\ell^{\frac{\gamma^2}{2}-2}d\ell$. Therefore the law of  $(\wt\cL_-,\wt\cL_+)$ agrees with~\eqref{eq:weld:segment}.

   To prove the final measurability statement,  for fixed $t>0$, let $(\cD^t,\eta^t)$ be the collection of concatenation of $(\cD_u,\eta_u)$ with $u<t$, and  $(\cL^t_-,\cL_+^t)$ be the part of $(\cL_-,\cL_+)$ in $\cD^t$. Then $(\cD^t,\eta^t)$ is independent of $(\cW\backslash\cD^t,\eta\backslash\eta^t)$. In particular, this means $((\cD^t,\eta^t), (\cL^t_-, \cL^t_+))$ is independent of  $(\cL_-\backslash\cL^t_-,\cL_+\backslash\cL^t_+)$. Moreover, $(\cD^t,\eta^t)$ is measurable with respect to $(\cL_-,\cL_+)$ by Theorem~\ref{thm:wedge-forest}. {Recall that if $A,B,C$ are random variables such that $(A,B)$ determine $C$ and $(A,C)$ is independent of $B$ then $A$ determines $C$. Applying this result in our setting we get that}  $(\cD^t,\eta^t)$ is measurable with respect to  $(\cL^t_-, \cL^t_+)$, and thus $(\tilde \cL_-,\tilde \cL_+)$ on the event $T>t$.  Since $t$ can be  arbitrary, we conclude the proof.
\end{proof}


To prove Theorem~\ref{thm:fd-disk}, we recall the following result on the conformal welding of quantum disks.
\begin{theorem}[Theorem 2.2 of \cite{AHS23}]\label{thm:disk-welding}
Let $\gamma\in(0,2),\tilde\kappa=\gamma^2$ and $W_-,W_+>0$. Then there exists a constant $c:=c_{W_-,W_+}\in(0,\infty)$ such that
$$\Md_2(W_-+W_+)\otimes \SLE_{\tilde\kappa}(W_--2;W_+-2) = c\,\Wd(\Md_2(W_-),\Md_2(W_+)).  $$
\end{theorem}

\begin{proof}[Proof of Theorem~\ref{thm:fd-disk}]
     Let $\tilde\kappa = \frac{16}{\kappa}$. Consider the welding of samples $( \cD_-,\wt\cL_-,\wt\cL_+,\cD_+)$ from
     \begin{equation}\label{eq:pf-weld-disk}
         \int_{\bbR_+^3} \mathcal{M}_2^\mathrm{disk}(W_-;\ell_-)\times \mathcal{M}_2^\mathrm{f.l.}(\ell_-;\ell)\times \mathcal{M}_2^\mathrm{f.l.}(\ell_+;\ell)\times \mathcal{M}_2^\mathrm{disk}(W_+;\ell_+)\, d\ell_-d\ell_+d\ell.
     \end{equation}
     By Proposition~\ref{prop:weld:segment} and a disintegration, for fixed $\ell_\pm$ in~\eqref{eq:pf-weld-disk},  we may first weld $(\wt\cL_-,\wt\cL_+)$ together, yielding a sample $(\cD_0,\eta_0)$ from $\Md_2(2-\frac{\gamma^2}{2};\ell_-,\ell_+)\otimes \SLE_\kappa(\frac{\kappa}{2}-4;\frac{\kappa}{2}-4)$. Then we may glue $\cD_-$ to $\cD_0$, and from Theorem~\ref{thm:disk-welding} we get a weight $W_-+2-\frac{\gamma^2}{2}$ quantum disk decorated by an independent $\SLE_{\tilde\kappa}(W_--2;-\frac{\tilde\kappa}{2})$ process. Finally, we glue $\cD_+$ to the right boundary of $\cD_0$. By Theorem~\ref{thm:disk-welding}, the surface and interfaces $(\cD,\eta_-,\eta_0,\eta_+)$ has law $\Md_2(W)\otimes\mathcal{P}$, where $\mathcal{P}$ is the probability measure with following description. For $W\ge\frac{\gamma^2}{2}$, when $\cD$ is embedded as $(\bbH,\phi,0,\infty)$, then under $\mathcal{P}$, (i) $\eta_+$ is the $\SLE_{\tilde\kappa}(W_--\frac{\tilde \kappa}{2};W_+-2)$ from 0 to $\infty$ (ii) given $\eta_+$, $\eta_-$ is the concatenation of  $\SLE_{\tilde\kappa}(W_--2;-\frac{\tilde\kappa}{2})$ curves in each connected component of $\bbH\backslash\eta_+$ to the left of $\eta_+$ and (iii) given $\eta_\pm$, $\eta_0$ is the concatenation of independent $\SLE_\kappa(\frac{\kappa}{2}-4;\frac{\kappa}{2}-4)$ curves in each pocket of $\bbH\backslash(\eta_-\cup\eta_+)$ between the two curves. On the other hand, following the \emph{SLE duality} argument ~\cite[Theorem 5.1]{zhan2008duality} and~\cite[Theorem 1.4 and Proposition 7.30]{MS16a}, the law of the union $\eta$ of $(\eta_-,\eta_0,\eta_+)$ agrees with  $\SLE_\kappa(\rho_-;\rho_+)$. If $W<\frac{\gamma^2}{2}$, the same argument applies for each bead of $\cD$ and the interface is the concatenation of independent $\SLE_\kappa(\rho_-;\rho_+)$ curves. Therefore we conclude the proof by foresting the boundary arcs of $(\cD,\eta)$.  
\end{proof}

\section{Multiple-SLE via conformal welding}
In this section we prove  Theorems~\ref{thm:existence-uniqueness},~\ref{thm:partition-func} and~\ref{thm:main}. 
The proof is based on an induction. We first prove Theorem~\ref{thm:main} for the $N=2$ case, and Theorem~\ref{thm:existence-uniqueness} and Theorem~\ref{thm:partition-func} for $N=2$ hold by~\cite{miller2018connection}. For the induction step, we apply the probabilistic construction in Section~\ref{subsec:pre-msle} to define the measure $\mSLE_{\kappa,\alpha}$ for $\alpha\in\LP_{N+1}$ and show that Theorem~\ref{thm:main} holds in this setting. Then using  the welding result from Theorem~\ref{thm:main}, we prove that the measure $\mSLE_{\kappa,\alpha}$ is finite, and the $N+1$ case of Theorem~\ref{thm:existence-uniqueness} and Theorem~\ref{thm:partition-func}    
follows  by the same arguments in~\cite{peltola2019toward}. This finishes the entire induction and concludes the proof of Theorems~\ref{thm:existence-uniqueness}-\ref{thm:main}.

\subsection{Multiple-SLE via conformal welding: two-curve case}\label{subsec:N=2}

For $\beta\in\bbR$, $ {\rho_-,\rho_+}> {\max\{-2,\frac{\kappa}{2}-4\}}$, define the measure $\wt{\SLE}_\kappa( {\rho_-;\rho_+};\beta)$ on curves $\eta$ from 0 to $\infty$ on $\mathbb{H}$ as follows. Let $D_\eta$ be the component of $\mathbb{H}\backslash \eta$  {whose boundary contains $1$}, and $\psi_\eta$ the unique conformal map from $D_\eta$ to $\mathbb{H}$ fixing 1 and sending the first (resp. last) point on $\partial D_\eta$ hit by $\eta$ to 0 (resp. $\infty$). Then our $\widetilde{\SLE}_\kappa(\rho_-;\rho_+;\beta)$ on $\bbH$ is defined by
\begin{equation}\label{eqn-sle-CR-1}
\frac{d \widetilde{\SLE}_\kappa(\rho_-;\rho_+;\beta)}{d {\SLE}_\kappa(\rho_-;\rho_+)}(\eta) = |\psi_{\eta}'(1)|^\beta.
\end{equation} 
This definition can be  extended to other domains via conformal transforms.

For $\beta\in\bbR$, recall the 
 {notation} $\Md_{2,\bullet}(W;\beta)$ from Definition~\ref{def:m2dot-alpha} for $W>\frac{\gamma^2}{2}$ and Definition~\ref{def:m2dot-thin} for $W\in(0,\frac{\gamma^2}{2})$. We write $\Mfd_{2,\bullet}(W;\beta)$ for the law of the  {generalized} quantum surface obtained by foresting the three boundary arcs of a sample from $\Md_{2,\bullet}(W;\beta)$.\footnote{Note that $\Md_{2,\bullet}(W;\gamma) = \Md_{2,\bullet}(W)$, while $\Mfd_{2,\bullet}(W;\gamma)$ is different from $\Mfd_{2,\bullet}(W)$.} 

The following is immediate from~\cite[Proposition 4.5]{AHS21}, Proposition~\ref{prop:weld:segment} and Theorem~\ref{thm:disk-welding}.
\begin{proposition}\label{prop:3-pt-disk}
	Let $\kappa\in(4,8)$ and $\gamma = \frac{4}{\sqrt{\kappa}}$. Suppose $W_+,W_->0$ and let $c_{W_-, W_+}\in (0,\infty)$ be the constant in  Theorem~\ref{thm:fd-disk}. Let $\rho_\pm = \frac{4}{\gamma^2}(2-\gamma^2+W_\pm)$, and $W = W_-+W_++2-\frac{\gamma^2}{2}$. Then for all $\beta\in\mathbb{R}$, 
	\begin{equation}\label{eqn-3-pt-disk}
	\begin{split}
	\mathcal{M}_{2,\bullet}^{\textup{f.d.}}(W;\beta)&\otimes \widetilde{\SLE}_\kappa(\rho_-;\rho_+;1-\Delta_{\beta})= c_{W_-, W_+}\Wd(\mathcal{M}_2^{\textup{f.d.}}(W_-), \mathcal{M}_{2,\bullet}^{\textup{f.d.}}(W_+;\beta)),
	\end{split}
	\end{equation} 
	where we are welding along the unmarked boundary arc of  $\mathcal{M}_{2,\bullet}^{\textup{f.d.}}(W_+;\beta)$ and  $\Delta_\beta = \frac{\beta}{2}(Q-\frac{\beta}{2})$.
\end{proposition}
\begin{proof}
    The proof is almost identical to that of Theorem~\ref{thm:fd-disk}. We consider the the welding of samples $( \cD_-,\wt\cL_-,\wt\cL_+,\cD_+)$ from
     \begin{equation}\label{eq:pf-weld-disk-a}
         \int_{\bbR_+^3} \mathcal{M}_2^\mathrm{disk}(W_-;\ell_-)\times \mathcal{M}_2^\mathrm{f.l.}(\ell_-;\ell)\times \mathcal{M}_2^\mathrm{f.l.}(\ell_+;\ell)\times \mathcal{M}_{2,\bullet}^\mathrm{disk}(W_+;\beta;\ell_+)\, d\ell_-d\ell_+d\ell,
     \end{equation}
     where $\mathcal{M}_{2,\bullet}^\mathrm{disk}(W_+;\beta;\ell_+)$ is the disintegration over the quantum  length of the unmarked boundary arc. Then we may first apply Proposition~\ref{prop:weld:segment} to glue $\wt\cL_-$ to $\wt\cL_+$ to get $(\cD_0,\eta_0)$, then apply Theorem~\ref{thm:disk-welding} to glue $\cD_-$ to $\cD_0$ from the left, and finally apply ~\cite[Proposition 4.5]{AHS21} to weld $\cD_+$ to $\cD_0$ on the right to get $(\cD,\eta_-,\eta_0,\eta_+)$. The interface law follows from the same SLE duality argument, and we conclude the proof by foresting the boundary arcs of $\cD$.
\end{proof}

We also need the disintegration of Liouville fields according to quantum lengths and the argument of changing weight of insertions in the conformal welding.
\begin{lemma}[{Lemma 3.2 of~\cite{SY23}}]\label{lm:lf-length}
    Let $m\ge 2$ and $0=s_1<s_2<  {\dots}<s_m=+\infty$. Fix $\beta_1, {\dots},\beta_m<Q$.  Let $C_\bbH^{(\beta_i,s_i)_i}$ and $P_\bbH$ be as in Definition~\ref{def-lf-H-bdry}, and  
    $\tilde{h} = h - 2Q\log|\cdot|_++\frac{1}{2}\sum_{i=1}^m\beta_iG_\bbH(s_i,\cdot)$, and $L=\nu_{\tilde{h}}((-\infty,0))$. For $\ell>0$, let $\LF_{\bbH,\ell}^{(\beta_i,s_i)_i}$ be the law of $\tilde{h}+\frac{2\ell}{\gamma L}$ under the reweighted measure $\frac{2}{\gamma}\frac{\ell^{\frac{1}{\gamma}(\sum_j\beta_j -2Q)-1}}{L^{\frac{1}{\gamma}(\sum_j\beta_j -2Q)}}\cdot\\ C_\bbH^{(\beta_i,s_i)_i}P_\bbH(dh)$. Then $\LF_{\bbH,\ell}^{(\beta_i,s_i)_i}$ is supported on $\{\phi:\nu_\phi((-\infty,0))=\ell\}$, and we have 
    \begin{equation}
       \LF_{\bbH}^{(\beta_i,s_i)_i} = \int_0^\infty\LF_{\bbH,\ell}^{(\beta_i,s_i)_i}d\ell.
    \end{equation}
\end{lemma}

\begin{lemma}[{Lemma 3.3 of~\cite{SY23}}]\label{lm:lf-change-weight}
   In the setting of Lemma~\ref{lm:lf-length}, for fixed $j\in\{2,\dots,m \}$ and $\beta_j'<Q$, we have the vague convergence of measures
    \begin{equation*}
        \lim_{\e\to 0}\e^{\frac{\beta_j'^2-\beta_j^2}{4}}e^{\frac{\beta_j'-\beta_j}{2}\phi_\e(s_j)}\LF_{\bbH,\ell}^{(\beta_i,s_i)_i}(d\phi) = \LF_{\bbH,\ell}^{(\beta_i,s_i)_{i\neq j}, (\beta_j',s_j)}(d\phi).
    \end{equation*}
\end{lemma}

\begin{figure}[H]
    \centering
    \begin{tabular}{ccc} 
		\includegraphics[scale=0.5]{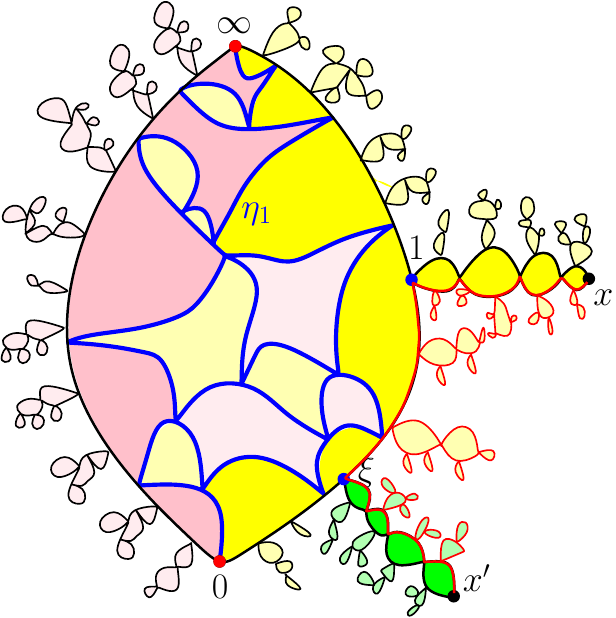}
		&
	   \includegraphics[scale=0.5]{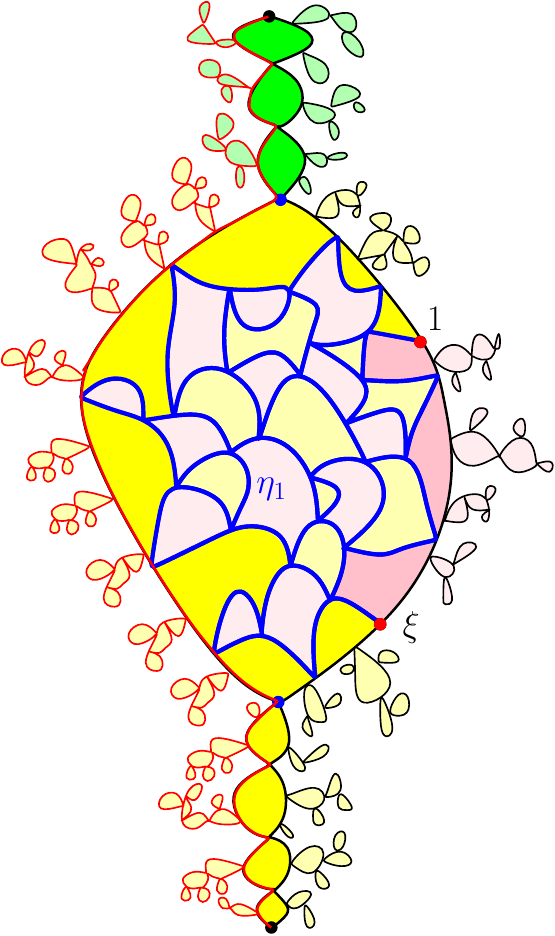}
    &\includegraphics[scale=0.5]{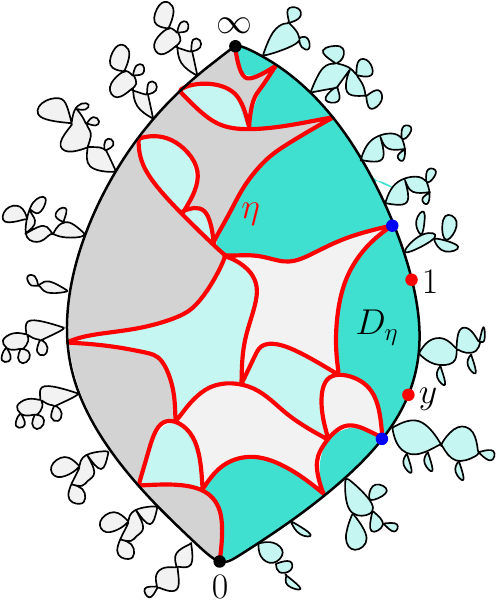}
	\end{tabular}
 \caption{An illustration of the proof of Proposition~\ref{prop:N=2}. \textbf{Left:} The setup in Step 1. The welding $S_+^f$ of a surface (pink) from $\GQD_{2}$  and a surface  (yellow) from $\GQD_{3}$. {In the proof,} the large disk $S_1$ in the picture is embedded as $(\bbH,\phi,0,1,\infty)$, while $S_2$ is the chain of (dark) yellow disks connecting 1 and $x$. We sample a marked point $x'$ on {the} boundary arc of $S_+^f$ from 0 to $x$ from the generalized quantum length measure, which gives the welding of $\GQD_2$ and $\GQD_4$. By further welding a sample $S_-^f$ from $\GQD_2$ to $S^f$ along the red boundary arc, we obtain the welding of $\GQD_2,\GQD_4$ and $\GQD_2$ as in the statement, and restricting to the event that the spine of the output is simply connected is the same as requiring $x'$ falling on the boundary arc between 0 and 1. By Lemma~\ref{lem:forest-decomposition}, this can be done by sampling a point $\xi$ on the spine $S_1$ from the quantum length measure and then concatenating an independent $\GQD_2$ (green) at the point $\xi$. The points $\xi,1$ are marked blue. \textbf{Middle:} The re-embedding of $S_+^f$ from the left panel in Step 2. Applying the conformal map $f_\xi(z)=\frac{\xi(z-1)}{z-\xi}$ to the  embedding of $S_1$,  the points $(0,\xi,1,\infty)$ from the left panel are mapped to $(1,\infty,0,\xi)$. 
 \textbf{Right:} The welding $\tilde S^f$ of a surface $\tilde S_-^f$ (grey) from $\GQD_2$ and a sample $\tilde S_+^f$ (turquoise) from $\Mfd_{2,\bullet}(\gamma^2-2;\beta)$ as in Proposition~\ref{prop:3-pt-disk}. The field $Y$ is  embedded on the spine of $\tilde S^f$, and $X$ is the restriction of $Y$ to the component $D_\eta$ with 1 on the boundary. In Steps 3 and 4, we sample a point $y$ on the boundary of $D_\eta$ from quantum length measure and change the singularity of $Y$ near $y$. By  drawing an $\SLE_\kappa$ from $y$ to 1 in $D_\eta$, $\tilde S_+^f$ would have the same law as the $S_+^f$ in the middle panel. In particular, the turquoise surface would agree with the surface $S_+^f$ in the middle panel, and the conformal maps between the three panels map the marked points to the marked points with the same color. 
 }
 \label{fig:pf-N=2}
 \end{figure}

\begin{proposition}\label{prop:N=2}
    Theorem~\ref{thm:main} holds when $N=2$.
\end{proposition}
At high level the proof is organized as follows. In Step 1 we conformally weld $\GQD_2$ and $\GQD_3$ via Proposition~\ref{prop:3-pt-disk} and add another boundary point from generalized quantum length measure to get $S_+^f$, see the left panel of Figure~\ref{fig:pf-N=2}. The welding $S_+^f$ with another sample $S_-^f$ from $\GQD_2$ would give the desired welding picture. In Step 2 we re-embed the surface $S_+^f$ as in the middle panel of Figure~\ref{fig:pf-N=2}. In Steps 3 and 4, we work on the conformal welding of a sample {$\tilde S_-^f$} from $\Mfd_{2}(\gamma^2-2)$ and a sample {$\tilde S_+^f$} from $\Mfd_{2,\bullet}(\gamma^2-2;\beta)$ as in Proposition~\ref{prop:3-pt-disk}. We modify the surface $\tilde S_+^f$ such that the welding of $\tilde S_-^f$ with $\tilde S_+^f$ agrees in law with the welding of $S_-^f$ with $S_+^f$, as in the right panel of Figure~\ref{fig:pf-N=2}. 

\begin{proof}
	\emph{Step 1: The setup.}
	We start with the conformal welding of two samples from $\GQD_{2}$ and $\GQD_{3}$. Let $c\in(0,\infty)$ be the constant in Proposition~\ref{prop:3-pt-disk} for $W_-=W_+=\gamma^2-2$. By Theorem~\ref{thm:fd-disk} (where we sample a marked point on the boundary from the generalized quantum length measure), we obtain a curve-decorated quantum surface $S_+^f$ whose law can be written as $c^{-1}\Mfd_{2,\bullet}(\frac{3\gamma^2}{2}-2)\otimes\SLE_\kappa$. By Lemma~\ref{lem:m2fd}, $S^f_+$ can be constructed by (i) sampling $((S_1,\eta_1),S_2)$ from $c^{-1}\big(\Md_{2,\bullet}(\frac{3\gamma^2}{2}-2)\otimes \SLE_\kappa)\times\Md_2(\gamma^2-2)$, (ii) concatenating $S_2$ to $S_1$ at its third marked point and (iii) foresting the boundary arcs of $S_1\cup S_2$. Moreover by Proposition~\ref{prop:m2dot}, we may assume $(S_1,\eta_1)$ is embedded as $(\bbH,\phi,\eta_1,0,\infty,1)$ where 
	\begin{equation}\label{eq:phieta1}
		(\phi,\eta_1)\sim \frac{\gamma}{2c(Q-\beta)^2}\LF_\bbH^{(\beta,0),(\beta,\infty),(\gamma,1)}\times \mu_\bbH(0,\infty),
	\end{equation}
	and  $S_2$ is sampled from $\Md_2(\gamma^2-2)$ and embedded as $(D_2, \phi_2, 1, x)$ with $\ol D_2\cap \ol \bbH = \{1\}$, and   $S_+^f$  is obtained by foresting the boundary of $(\bbH\cup D_2, \phi\sqcup\phi_2, 0, x, \infty)${, where $\phi\sqcup\phi_2$ is the distribution given by $\phi$ (resp.\ $\phi_2$) when restricted to the domain of $\phi$ (resp.\ $\phi_2$)}, and $\beta=\frac{4}{\gamma}-\frac{\gamma}{2}$. We sample a marked point $x'$ on the boundary arc of {$S_+^f$} from $0$ to $x$ from the generalized quantum length measure; this induces a weighting by the generalized quantum length of this boundary arc, so we get the conformal welding of $\GQD_2$ and $\GQD_4$. Restrict to the event that $x'$ lies on the boundary arc of $S_+^f$ from $0$ to 1. See the left panel of Figure~\ref{fig:pf-N=2}.

 Note that if one were to further glue a sample $S_-^f$ from $\GQD_{2}$ to $S_+^f$ along the forested boundary arc between $x$ and $x'$, 
 then one would obtain the conformal welding of two samples from $\GQD_{2}$ and one sample from $\GQD_{4}$. The restriction on $x'$ lying between $0$ and $1$ corresponds to the spine of the conformal welding being simply connected, as in the theorem statement. 

    \emph{Step 2: Add the new marked point and re-embed the surface $S_+^f$.}
   Following Lemma~\ref{lem:forest-decomposition}, to construct the surface $S_+^f$ with $x'$ on the boundary, we can (i) weight the law of $(\phi,\eta_1)$ from~\eqref{eq:phieta1} by $\nu_\phi([0,1])$ and sample a point $\xi\in(0,1)$ from the quantum length measure on $[0,1]$ and (ii) sample a surface $S_3^f$ from $\Mfd_2(\gamma^2-2)$ and attach it to $S_+^f$ at the point $\xi$.  {Indeed, once we disintegrate over $\{\nu_\phi([0,1])=\ell\}$,  the procedure of taking the point $x'$ according to generalized quantum length is described in the first procedure in Lemma~\ref{lem:forest-decomposition}, and the procedure discussed just above matches the second procedure in Lemma~\ref{lem:forest-decomposition}. The weighting over $\nu_\phi([0,1])$ comes from the weighting over $\ell$ when taking $u\sim \mathrm{Leb}_{[0,\ell]}$ as in Lemma~\ref{lem:forest-decomposition}.}
    By Lemma~\ref{lm:gamma-insertion}, the law of $(\phi,\xi)$ after procedure (i) is now
    \begin{equation}\label{eq:pf-N=2-1}
       \frac{\gamma}{2c(Q-\beta)^2} \int_0^1 \LF_\bbH^{(\beta,0),(\beta,\infty),(\gamma,1),(\gamma,\xi)}(d\phi) d\xi.
    \end{equation}
    On the other hand, if we perform the coordinate change $z\mapsto f_\xi(z):=\frac{\xi(z-1)}{z-\xi}$, then by Lemma~\ref{lm:lcft-H-conf}, when viewed as quantum surfaces with marked points, ~\eqref{eq:pf-N=2-1} is equal to 
\begin{equation}\label{eq:pf-N=2-2}
\begin{split}
    &\ \ \frac{\gamma}{2c(Q-\beta)^2}\int_0^1 f_\xi'(0)^{\Delta_\beta}f_\xi'(\infty)^{\Delta_\beta}f_\xi'( {\xi})^{\Delta_\gamma}f_\xi'(1)^{\Delta_\gamma}\LF_\bbH^{(\beta,f_\xi(0)),(\beta,f_\xi(\infty)),(\gamma,f_\xi(\xi)),(\gamma,f_\xi(1))}(d\phi)\,d\xi \\&=\frac{\gamma}{2c(Q-\beta)^2}\int_0^1 (1-\xi)^{2\Delta_\beta-2}\, \LF_\bbH^{(\gamma,0),(\beta,\xi),(\beta,1),(\gamma,\infty)}(d\phi)\,d\xi.
    \end{split}
\end{equation}
    In other words, as shown in the middle panel of Figure~\ref{fig:pf-N=2}, the quantum surface $S_+^f$ can be constructed by (i) sampling $S_1:=(\bbH,\phi,0,1,\infty,\xi)$ where $(\phi, {\xi})$ is from the law~\eqref{eq:pf-N=2-2} and drawing an independent $\SLE_\kappa$ curve $\eta_1$ from $\xi$ to 1 (ii) sampling $S_2,S_3$ from $\Md_2(\gamma^2-2)\times\Md_2(\gamma^2-2)$ and concatenating $S_2,S_3$ to $S_1$ at the points $0$ and $\infty$ and (iii) foresting the six boundary arcs of $S_1\cup S_2\cup S_3$.

    \emph{Step 3: Add a typical point to the welding of $\Mfd_{2}(\gamma^2-2)$ and $\Mfd_{2,\bullet}(\gamma^2-2;\beta)$.}  We work on the conformal welding  {$\tilde S^f$} of a sample $\tilde S_-^f$ from $\Mfd_{2}(\gamma^2-2)$ and a sample $\tilde S_+^f$ from $\Mfd_{2,\bullet}(\gamma^2-2;\beta)$. 
    By Proposition~\ref{prop:3-pt-disk} and Definition~\ref{def:m2dot-alpha}, the surface can be constructed by foresting the boundary of $(\bbH,Y,\eta,0,1,\infty)$ with $(Y,\eta)$ sampled from 
    \begin{equation}\label{eq:pf-N=2-0}
        \frac{\gamma}{2c(Q-\beta)^2}\LF_\bbH^{(\beta,0),(\beta,1),(\beta,\infty)}\times \psi_\eta'(1)^{1-\Delta_\beta}\mu_\bbH(0,\infty)(d\eta).
    \end{equation}
   In this step and next step, we shall add a marked point to $\tilde S_+^f$ and change the boundary insertion via Lemma~\ref{lm:lf-change-weight}. The surface $\tilde S_+^f$ will eventually have the same law as $S_+^f$ as in the conclusion of Step 2, and the welding of  $\tilde S_+^f$ with $\tilde S_-^f$ will agree in law with that of $S^f_+$ with $S_-^f$.  Moreover, the law of this conformal welding is~\eqref{eq:pf-N=2-8}, completing the proof.

   Let $D_\eta^-$ (resp.\ $D_\eta^+$) be the union of the connected components of $\bbH\backslash\eta$ whose boundaries contain a segment of $(-\infty,0)$ (resp.\ $(0,\infty)$), and (recall the notion of $\psi_\eta$ and $D_\eta$ in~\eqref{eqn-sle-CR-1}) 
\begin{equation}\label{eq:pf-N=2-3}
   X=Y\circ\psi_\eta^{-1}+Q\log|(\psi_\eta^{-1})'|. 
\end{equation}
     Let $\tilde S_-=(D_\eta^-,Y)/{\sim_\gamma}$,   $\tilde S_+=(D_\eta^+,Y)/{\sim_\gamma}$, $\tilde S_1 = (\bbH,X,0,\infty,1)/{\sim_\gamma}$. Then $\tilde S_-$ and $\tilde S_+$ are the spines of $\tilde S_-^f$ and $\tilde S_+^f$. 
     We sample a marked point $\xi$ on $\tilde S_1$ from the measure $1_{\xi\in(0,1)}(1-\xi)^{2\Delta_\beta-2}\nu_X(d\xi)$.
Then by Lemma~\ref{lm:gamma-insertion}, the surface $\tilde S_+$ is the concatenation of two samples $\tilde S_2$ and $\tilde S_3$ from $\Md_2(\gamma^2-2)$ with a sample $(\bbH,X,0,\xi,1,\infty)$ from
\begin{equation}\label{eq:tildeS1}
    \frac{2}{\gamma}\int_0^1 (1-\xi)^{2\Delta_\beta-2}\, \LF_\bbH^{(\gamma,0),(\gamma, {\xi}),(\beta,1),(\gamma,\infty)}(dX)\,d\xi.
\end{equation}
 {The surfaces $\tilde S_2$ and $\tilde S_3$ are attached to the latter surface}
at the points $0$ and $\infty$. The constant $\frac{2}{\gamma}$ follows from Proposition~\ref{prop:m2dot}, Definition~\ref{def:m2dot-thin} and $(1-\frac{2(\gamma^2-2)}{\gamma^2})^2\cdot\frac{\gamma}{2(Q-\gamma)^2} = \frac{2}{\gamma}$.
On the other hand, for $y=\psi_\eta^{-1}(\xi)$, the law of $(Y,\eta,y)$ is given by
\begin{equation}\label{eq:pf-N=2-4}
    \frac{\gamma}{2c(Q-\beta)^2}\bigg[\int_0^1 \mathds{1}_{E_{\eta,y}} (1-\psi_\eta(y))^{2\Delta_\beta-2}\,\nu_Y(dy)\,\LF_\bbH^{(\beta,0),(\beta,\infty),(\beta,1)}(dY)\bigg]\cdot\psi_\eta'(1)^{1-\Delta_\beta}\mu_\bbH(0,\infty)(d\eta),
\end{equation}
where $E_{\eta,y}$ is the event where $y$ and $1$ lie on the boundary of the same connected component of $\bbH\backslash\eta$. By Lemma~\ref{lm:gamma-insertion}, ~\eqref{eq:pf-N=2-4} is equal to 
\begin{equation}\label{eq:pf-N=2-5}
    \frac{\gamma}{2c(Q-\beta)^2}\bigg[\int_0^1 \mathds{1}_{E{\eta,y}}(1-\psi_\eta(y))^{2\Delta_\beta-2}\,\LF_\bbH^{(\beta,0),(\beta,\infty),(\beta,1),(\gamma,y)}(dY)\, {dy}\bigg]\cdot\psi_\eta'(1)^{1-\Delta_\beta}\mu_\bbH(0,\infty)(d\eta).
\end{equation}

\emph{Step 4: Change the insertion from $\gamma$ to $\beta$.} We weight the law of $(Y,\eta,y)$ from~\eqref{eq:pf-N=2-5} by $\frac{\gamma^2}{4c(Q-\beta)^2}\e^{\frac{\beta^2-\gamma^2}{4}} e^{\frac{\beta-\gamma}{2}X_\e(\xi)}$, where $X$ is given by~\eqref{eq:pf-N=2-3} and $\xi=\psi_\eta(y)$. Then following from the same argument as in~\cite[Proposition 4.5]{AHS21}, we have:
\begin{enumerate}[(i)]
    \item Given  {the LQG generalized quantum length} $\ell$ of $\eta$, the law of $\tilde S_-^f$ is unchanged and is given by $\Mfd_2(\gamma^2-2;\ell)$. By Lemma~\ref{lm:lf-change-weight}, as $\e\to0$, given the quantum length $\ell_0$ of $\nu_X((-\infty,0))$, the law of $(\tilde S_1,\xi)$ from~\eqref{eq:tildeS1} converges in vague topology to
    \begin{equation}\label{eq:pf-N=2-6}
        \frac{\gamma}{2c(Q-\beta)^2}\int_0^1 (1-\xi)^{2\Delta_\beta-2}\, \LF_{\bbH,\ell_0}^{(\gamma,0),(\beta,y),(\beta,1),(\gamma,\infty)}(d\phi)\,d\xi,
    \end{equation}
    In particular, by comparing with~\eqref{eq:pf-N=2-2}, the conformal welding of $\tilde S_-^f$ and $\tilde S_+^f$ {(after changing the singularity from $\gamma$ to $\beta$)} agree in law with the conformal welding of $S_-^f$ and $S_+^f$.
    \item The law of 
    $(Y, {\eta,y})$ is weighted by 
    \begin{equation}\label{eq:pf-N=2-7}
        \frac{\gamma^2}{4c(Q-\beta)^2}\e^{\frac{\beta^2-\gamma^2}{4}} e^{\frac{\beta-\gamma}{2}\big((Y,\theta_\e(y))+Q\log|(\psi_\eta^{-1})'(\xi)| \big)} = \frac{\gamma^2}{4c(Q-\beta)^2}\bigg(\frac{\e}{\psi_\eta'(y)}\bigg)^{\frac{\beta^2-\gamma^2}{4}}e^{\frac{\beta-\gamma}{2}(Y,\theta^\eta_\e(y))}\big|\psi_\eta'(y)\big|^{1-\Delta_\beta} 
    \end{equation}
    where $\theta_\e^\eta$ is  {the} push-forward of the uniform probability measure on $B_\e(y)\cap\bbH$ under $\psi_\eta^{-1}$ and we used the fact that $\log|(\psi_\eta^{-1})(z)|$ is a harmonic function along with Schwartz reflection. As argued in~\cite[Proposition 4.5]{AHS21}, by Girsanov's theorem, under the weighting~\eqref{eq:pf-N=2-7}, as $\e \to0$, the law of $(Y,\eta,y)$ from~\eqref{eq:pf-N=2-5} converges in vague topology to $\frac{\gamma^3}{8c^2(Q-\beta)^4}=c_2$ times 
    \begin{equation}\label{eq:pf-N=2-8}
        \bigg[\int_0^1 \mathds{1}_{E{\eta,y}}\big|\psi_\eta'(y)\big|^{1-\Delta_\beta} (1-\psi_\eta(y))^{2\Delta_\beta-2}\,\LF_\bbH^{(\beta,0),(\beta,\infty),(\beta,1),(\beta,y)}(dY)\,dy\bigg]\cdot\psi_\eta'(1)^{1-\Delta_\beta}\mu_\bbH(0,\infty)(d\eta).
    \end{equation}
    Intuitively, this is because when $\e\to0$,  $\theta_\e^\eta$ is roughly the uniform measure on $B_{\frac{\e}{\psi_\eta'(y)}}\cap\bbH$,  {and $(Y, \theta^\eta_\e(y))$ 
     {is} close to $Y_{\frac{\e}{\psi_\eta'(y)}}(y)$, i.e., the average of $Y$ over $\partial B_{\frac{\e}{\psi_\eta'(y)}}(y)\cap\bbH$. For fixed $y$ and $\eta$, as we weight the law of $Y$ by $\big(\frac{\e}{\psi_\eta'(y)}\big)^{\frac{\beta^2-\gamma^2}{4}}e^{\frac{\beta-\gamma}{2}Y_{\frac{\e}{\psi_\eta'(y)}}(y)}$, then the law of $Y$ converges in vague topology to $\LF_\bbH^{(\beta,0),(\beta,\infty),(\beta,1),(\beta,y)}$ thanks to Lemma~\ref{lm:lf-change-weight} (with $\e$ replaced by $\frac{\e}{\psi_\eta'(y)}$). See also the explanations after~\cite[Eq.(4.12)]{AHS21} for more details.} 
\end{enumerate}
On the other hand, observe that $1-\Delta_\beta = \frac{6-\kappa}{2\kappa} = b$ and $(1-\psi_\eta(y))^{2\Delta_\beta-2}\big|\psi_\eta'(y)\big|^{1-\Delta_\beta}\psi_\eta'(1)^{1-\Delta_\beta} = H_{D_\eta}(y,1)^b$. This implies that if we further draw the interface $\eta_1$ in $D_\eta$ from $y$ to $1$, from the construction of the multiple SLE in Section~\ref{subsec:pre-msle}, the joint law of $(Y,y,\eta,\eta_1)$ in~\eqref{eq:pf-N=2-8} can be described in terms of~\eqref{eq:thm-main}. This concludes the proof for the $N=2$ case.
\end{proof}

\subsection{Multiple-SLE via conformal welding: general case}

In this section we work on the induction step for the proof of Theorem~\ref{thm:main}. To be more precise, we prove the following:
\begin{proposition}\label{prop:induction}
    Suppose Theorems~\ref{thm:partition-func}--\ref{thm:main} holds for $1,2,...,N$. Let $\alpha\in\LP_{N+1}$, and define the measure $\mSLE_{\kappa,\alpha}$ as in Section~\ref{subsec:pre-msle}. Then Theorem~\ref{thm:main} holds for $\alpha$. 
\end{proposition}

We first show that the expression~\eqref{eq:thm-main} has the following rotational invariance. Given two link patterns $\alpha = \{\{i_1,j_1\},...,\{i_N,j_N\}\}$ and $\alpha'=\{\{i_1',j_1'\},...,\{i_N',j_N'\}\}$ in $\mathrm{LP}_N$, we say $\alpha$ and $\alpha'$ are rotationally equivalent if there exists some integer $0\le m\le 2N-1$ such that for every $1\le k\le N$, $i_k'=i_k+m$ and $j_k'=j_k+m$ ($\mathrm{mod}$ $2N$), and we write $\alpha'=\alpha+m$. 
\begin{lemma}\label{lm:rotation}
    In the setting of Proposition~\ref{prop:induction}, the measure $\mSLE_{\kappa,\alpha}$ satisfies the conformal covariance property~\eqref{eq:msle-conformal-covariance}.
     Moreover, for any $0\le m\le 2N+1$, the expression~\eqref{eq:thm-main} for $N+1$ when viewed as  {a} {(non-forested)} curve-decorated quantum 
      {surface} 
     is equal to 
    \begin{equation}
       c_{N+1}\int_{0<y_1<...<y_{2N-1}<1}\LF_\bbH^{(\beta,0),(\beta,1),(\beta,\infty),(\beta,y_1),...,(\beta,y_{2N-1})}\times \mathrm{mSLE}_{\kappa,\alpha+m}(\bbH,0,y_1,...,y_{2N-1},1,\infty)dy_1...dy_{2N-1}.
    \end{equation}
\end{lemma}
\begin{proof}
From the conformal covariance of the measure $\mSLE_{\kappa,\alpha_0}$ for $\alpha_0\in \bigsqcup_{k=1}^N \LP_k$,  
it is straightforward to verify that the measure $\mSLE_{\kappa,\alpha}$ satisfies~\eqref{eq:msle-conformal-covariance} for any conformal map $f:\bbH\to\bbH$. {Combining this with the conformal covariance of Liouville fields (Lemma~\ref{lm:lcft-H-conf}) and using the relation $\Delta_\beta+b=1$ gives the result; see~\cite[Lemma 3.6]{SY23} for a similar computation.} 
\end{proof}

For $\alpha\in\LP_N$ and $(\eta_1,...,\eta_N)\in X_\alpha(D;x_1,...,x_{2N})$, suppose $\eta_k$ links $x_{i_k}$ and $x_{j_k}$. We call $\eta_k$ a \emph{good link} if $j_k=i_k+1 (\mathrm{mod}\, 2N)$, and $x_1,...,x_{i_k-1},x_{i_k+2},...,x_{2N}$ are lying on the boundary of the same connected component of $D\backslash\eta_k$.
\begin{figure}[htb]
    \centering
    \begin{tabular}{ccc} 
		\includegraphics[scale=0.5]{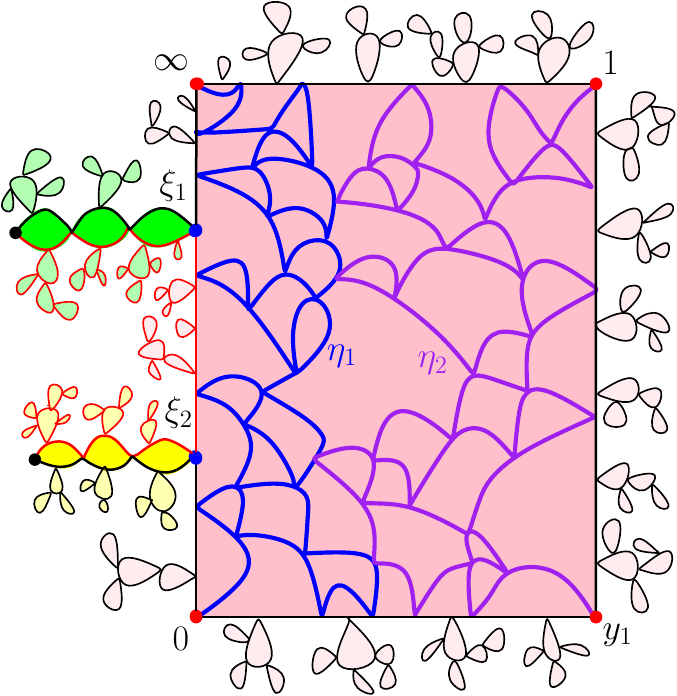}
		&
	   \includegraphics[scale=0.5]{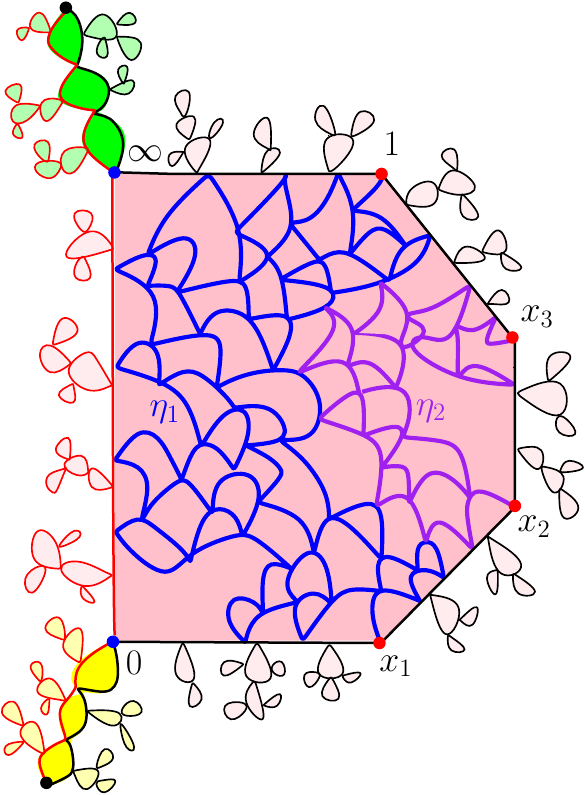}
    &\includegraphics[scale=0.5]{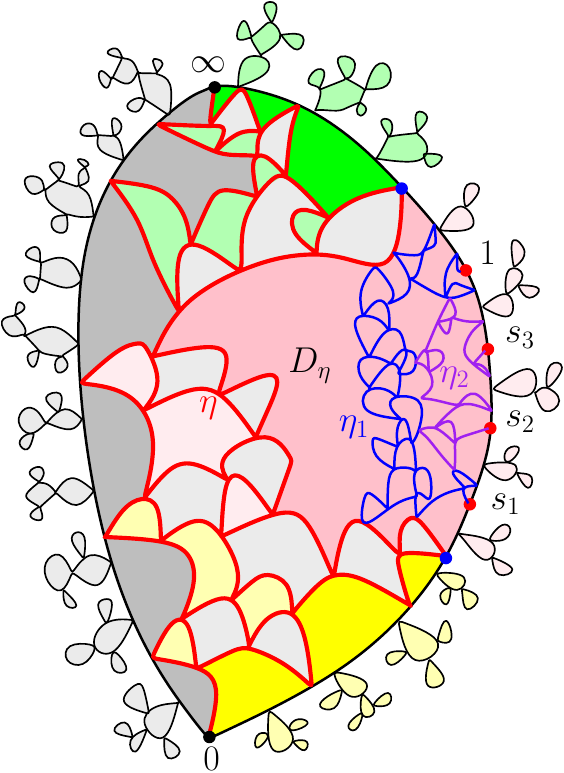}
	\end{tabular}
 \caption{An illustration of the proof of Proposition~\ref{prop:induction} for $N+1 = 3$. \textbf{Left:} The surface $S_+^f$ constructed in Step 1. If we glue a sample $S_-^f$ from $\GQD_2$, we obtain the welding $\Wd_\alpha(\GQD^{4})$ for $\alpha = \{\{1,6\}, \{2,5\},\{3,4\}\}$ restricted to the event $E_1$ that the spine is simply connected and $\{1,6\}$ is a good link. \textbf{Middle:} The re-embedding of $S_+^f$ in Step 2 via the conformal map $f_{\xi_1,\xi_2}(z)=\frac{z-\xi_2}{z-\xi_1}$ and a change of variables. \textbf{Right:} The welding $\tilde S^f$ of a surface $\tilde S_-^f$ (grey) from $\GQD_2$ and a surface $\tilde S_+^f$ (union of pink, green and yellow) from $\Mfd_{2,\bullet}(\gamma^2-2;\beta)$ as in Proposition~\ref{prop:3-pt-disk}. In Steps 3 and 4 we modify the surface $\tilde S_+^f$ as in Proposition~\ref{prop:N=2} such that $\tilde S_+^f$ agrees in law with the surface $S_+^f$ in the middle panel. This will give the expression of $\Wd_\alpha(\GQD^{4})$ restricted to the event $E_1$, and in Step 5 we remove this extra constraint from $E_1$.}\label{fig:pf-N-induc}
 \end{figure}

\begin{proof}[Proof of Proposition~\ref{prop:induction}]
 The proof is similar to Proposition~\ref{prop:N=2} based on induction. By Lemma~\ref{lm:rotation}, without loss of generality we may assume $\{1,2N+2\}\in\alpha$, and let $\hat\alpha\in\mathrm{LP}_N$ be the link pattern induced by $\alpha\backslash\{1,2N+2\}$. On the event $E$, let $E_1$ be the event where the link $\{1,2N+2\}$ in the conformal welding picture $\Wd_\alpha(\GQD^{N+2})$ is a good link. In Steps 1 and 2, we begin with the welding $\Wd_{\hat\alpha}(\GQD^{N+1})$ restricted to the event that the spine is simply connected as from induction hypothesis. We construct and re-embed the surface $S_+^f$ by adding two new marked points and re-embedding, such that the welding of $S_+^f$ with a surface $S_-^f$ from $\GQD_2$ gives $\Wd_{\hat\alpha}(\GQD^{N+1})$ restricted to the event $E_1$. In Steps 3 and 4, we begin with the welding of $\tilde S_-^f$ and $\tilde S_+^f$ as in the proof of  Proposition~\ref{prop:N=2} and modify $\tilde S_+^f$ such that the law of $\tilde S_+^f$ would agree with that of $S_+^f$. Finally, in Step 5 we remove the extra constraint $E_1$. 
 
\emph{Step 1: The setup.}
     We first prove
    Theorem~\ref{thm:main} when restricted to the event $E_1$.
    
    Let $(\bbH,\phi,y_1,...,y_{2N-3},\eta_1,...,\eta_N)$ be a sample from ~\eqref{eq:thm-main} with link pattern $\hat\alpha$ where the multiple $\SLE$ is sampled from $\mathrm{mSLE}_{\kappa,\hat\alpha}(\bbH,0,y_1,...,y_{2N-3},1,\infty)$.  We truncate and glue independent forested lines to the boundary to obtain the curve-decorated quantum surface $S_+^f$. By our induction hypothesis,  we obtain the conformal welding  $\Wd_{\hat\alpha}(\GQD^{N+1})$ restricted to the event where all the marked points are on the same connected component. We weight the law of  $S_+^f$ by the generalized quantum length of the forested boundary segment between $-\infty$ and $0$ and sample two marked points from the generalized quantum length measure on this segment. As in the proof of Proposition~\ref{prop:N=2}, we restrict to the event where the two points are on different trees of disks. From Lemma~\ref{lem:forest-decomposition} and Lemma~\ref{lm:gamma-insertion}, this is the same as sampling $(\bbH,\phi,\xi_1,\xi_2,y_1,...,y_{2N-3},\eta_1,...,\eta_N)$ from 
    \begin{equation}\label{eq:pf-induc-1}
  c_N\int_{\xi_1<\xi_2<0<y_1<...<y_{2N-3}<1}\bigg[\LF_\bbH^{(\beta,y_1),...,(\beta,y_{2N}),(\gamma,\xi_1),(\gamma,\xi_2)}\times \mathrm{mSLE}_{\kappa,\hat\alpha}(\bbH,y_1,...,y_{2N})\bigg]\,dy_1...dy_{2N-3}d\xi_1d\xi_2,
\end{equation}
foresting the boundary, and insert two samples from $\Mfd_{2 }(\gamma^2-2)\times\Mfd_{2}(\gamma^2-2)$ to the points $\xi_1,\xi_2$. Here we used the convention $(y_{2N-2},y_{2N-1},y_{2N})=(1,\infty,0)$.
If we glue a sample $S_-^f$ from $\GQD_{2}$ to $S_+^f$ along the two newly sampled marked points, then the output equals  {$\Wd_\alpha(\GQD^{N+2})|_{E_1\cap E}$.  Indeed, without the constraint that the spine of $\Wd_{\hat\alpha}(\GQD^{N+1})$ is simply connected and the two newly sampled marked points are on different trees of disks, then we get the surface $\Wd_\alpha(\GQD^{N+2})$. The spine for the welding $\Wd_\alpha(\GQD^{N+2})$ is simply connected and the link $\eta_1$ is a good link if and only if the spine for $\Wd_{\hat\alpha}(\GQD^{N+1})$ is simply connected and the two newly sampled marked points are on different trees of disks as described.}

\emph{Step 2: The re-embedding of the surface $S^f$.} For $\xi_1<\xi_2<0$, consider the conformal map $f_{\xi_1,\xi_2}(z)=\frac{z-\xi_2}{z-\xi_1}$ from $\bbH$ to $\bbH$ sending $(\xi_1,\xi_2,\infty)$ to $(\infty,0,1)$. Let $x_1=f_{\xi_1,\xi_2}(0)$, $x_{2N-1}=f_{\xi_1,\xi_2}(1)$, $x_{2N}=1$ and $x_k = f_{\xi_1,\xi_2}(y_{k-1})$ for $2\le k\le 2N-2$.  Then by Lemma~\ref{lm:lcft-H-conf} and~\eqref{eq:msle-conformal-covariance}, when viewed as the law of curve-decorated quantum surfaces, ~\eqref{eq:pf-induc-1} is equal to
    \begin{equation}\label{eq:pf-induc-2}
    \begin{split}
        c_N\int_{\xi_1<\xi_2<0<y_1<...<y_{2N-3}<1}&\bigg[f_{{\xi_1,\xi_2}}'(\xi_1)f_{{\xi_1,\xi_2}}'(\xi_2)\prod_{k=1}^{2N}f_{{\xi_1,\xi_2}}'(y_k) \cdot\LF_\bbH^{(\beta,f_{{\xi_1,\xi_2}}(y_1)),...,(\beta,f_{{\xi_1,\xi_2}}(y_{2N})),(\gamma,0),(\gamma,\infty)}\\&\times \mathrm{mSLE}_{\kappa,\hat\alpha}(\bbH,f_{{\xi_1,\xi_2}}(y_1),\,...,f_{{\xi_1,\xi_2}}(y_{2N}))\bigg]\,dy_1...dy_{2N-3}d\xi_1d\xi_2
    \end{split}
\end{equation}
where we used $\Delta_\beta+b=\Delta_\gamma=1$. Then by a change of variables $x_k = f_{\xi_1,\xi_2}(y_{k-1})$ for $2\le k\le 2N-2$, ~\eqref{eq:pf-induc-2} is equal to
\begin{equation}\label{eq:pf-induc-3}
\begin{split}
     &c_N\int_{0<f_{{\xi_1,\xi_2}}(0)<x_2<...<x_{2N-2}<f_{{\xi_1,\xi_2}}(1) <1}\bigg[f_{{\xi_1,\xi_2}}'(\xi_1)f_{{\xi_1,\xi_2}}'(\xi_2)f_{{\xi_1,\xi_2}}'(0) f_{{\xi_1,\xi_2}}'(1)f_{{\xi_1,\xi_2}}'(\infty)\cdot\\&\LF_\bbH^{(\beta,f_{{\xi_1,\xi_2}}(0)),(\beta,x_2),...,(\beta,x_{2N-2}), (\beta,f_{{\xi_1,\xi_2}}(1)),(\beta,1),(\gamma,0),(\gamma,\infty)}\times \mathrm{mSLE}_{\kappa,\hat\alpha}(\bbH,x_2,...,x_{2N-2})\bigg]dx_2...dx_{2N-2}d\xi_1d\xi_2.
     \end{split}
\end{equation}
Since $x_1=f_{{\xi_1,\xi_2}}(0)=\frac{\xi_2}{\xi_1}$, $x_{2N-1}=f_{{\xi_1,\xi_2}}(1)=\frac{1-\xi_2}{1-\xi_1}$, it is straightforward to check that
$$\frac{\partial(x_1,x_{2N-1})}{\partial(\xi_1,\xi_2) } = \frac{\xi_2-\xi_1}{\xi_1^2(1-\xi_1)^2}. $$
On the other hand, we may compute
$$f_{{\xi_1,\xi_2}}'(\xi_1)f_{{\xi_1,\xi_2}}'(\xi_2)f_{{\xi_1,\xi_2}}'(0) f_{{\xi_1,\xi_2}}'(1)f_{{\xi_1,\xi_2}}'(\infty) =  \frac{\xi_2-\xi_1}{\xi_1^2(1-\xi_1)^2}.$$
Therefore by a change of variables to~\eqref{eq:pf-induc-3}, the law of $(\bbH,\phi,\eta_1,...,\eta_N,0,y_1,...,y_{2N-3},\xi_1,\xi_2,1,\infty)/\sim_\gamma$ agrees with that of $(\bbH,\tilde\phi,\eta_1,...,\eta_N,x_1,...,x_{2N-1},0,1,\infty)/\sim_\gamma$, where $(\tilde\phi,\tilde\eta_1,...,\tilde\eta_N,x_1,...,x_{2N-1})$ is sampled from
\begin{equation}\label{eq:pf-induc-4}
    \begin{split}
        c_N\int_{0<x_1<...<x_{2N-1}<1}&\bigg[\LF_\bbH^{(\beta,x_1),...,(\beta,x_{2N-1}),(\beta,1),(\gamma,0),(\gamma,\infty)}\times \mathrm{mSLE}_{\kappa,\hat\alpha}(\bbH,x_1,...,x_{2N-1},1)\bigg]\,dx_1...dx_{2N-1}.
    \end{split}
\end{equation}
    
\emph{Step 3: Add boundary typical points to the welding of $\Mfd_{2}(\gamma^2-2)$ and $\Mfd_{2,\bullet}(\gamma^2-2;\beta)$.} Parallel to Step 2 the proof of Proposition~\ref{prop:N=2}, consider the conformal welding of a sample $\tilde S_-^f$ from $\Mfd_{2}(\gamma^2-2)$ and a sample $\tilde S_+^f$ from $\Mfd_{2,\bullet}(\gamma^2-2;\beta)$. Define $Y,\eta,X,\tilde S_-,\tilde S_+, \tilde S_1$ 
in the same way as in the proof of Proposition~\ref{prop:N=2}. By Lemma~\ref{lm:gamma-insertion}, as we sample $2N-1$ marked points on $\partial \tilde S_1$ from the measure
$$\mathds{1}_{0<x_1<...<x_{2N-1}<1}\cdot\mathcal{Z}_{\hat\alpha}(\bbH,x_1,...,x_{2N-1},1)\nu_X(dx_1)...\nu_X(dx_{2N-1}), $$
the surface $\tilde S_+$ is the concatenation of two samples $\tilde S_2$, $\tilde S_3$ from $\Md_2(\gamma^2-2)$ with a sample $(\bbH,X,0,x_1,...,\\x_{2N-1},1,\infty)$ from
\begin{equation}\label{eq:pf-induc-5}
    \begin{split}
        \frac{2}{\gamma}\int_{0<x_1<...<x_{2N-1}<1}&\bigg[\LF_\bbH^{(\gamma,x_1),...,(\gamma,x_{2N-1}),(\beta,1),(\gamma,0),(\gamma,\infty)}\times \mathcal{Z}_{\hat\alpha}(\bbH,x_1,...,x_{2N-1},1)\bigg]\,dx_1...dx_{2N-1}.
    \end{split}
\end{equation}
at the points 0 and $\infty$.
On the other hand, for $s_k = \psi_\eta^{-1}(x_k)$ and $k=1,...,2N-1$, the law of $(Y,\eta,s_1,...,s_{2N-1})$ is given by
\begin{equation}\label{eq:pf-induc-6}
    \begin{split}
        \frac{\gamma}{2c(Q-\beta)^2}\int_{0<s_1<...<s_{2N-1}<1}&\bigg[\mathds{1}_{E_{\eta,s_1,...,s_{2N-1}}}\mathcal{Z}_{\hat\alpha}\big(\bbH,\psi_\eta(s_1),...,\psi_\eta(s_{2N-1}),1\big)\nu_Y(ds_1)...\nu_Y(ds_{2N-1}) \bigg]\cdot\\&\LF_\bbH^{(\beta,1),(\beta,0),(\beta,\infty)}(dY)\cdot\psi_\eta'(1)^{1-\Delta_\beta}\,\SLE_\kappa(d\eta),  
    \end{split}
\end{equation}
where $E_{\eta,s_1,...,s_{2N-1}}$ is the event where $s_1,...,s_{2N-1},1$  lie on the boundary of the same connected component of $\bbH\backslash\eta$. By    Lemma~\ref{lm:gamma-insertion}, ~\eqref{eq:pf-induc-6} is equal to 
\begin{equation}\label{eq:pf-induc-7}
    \begin{split}
        \frac{\gamma}{2c(Q-\beta)^2}\int_{0<s_1<...<s_{2N-1}<1}&\bigg[\mathds{1}_{E_{\eta,s_1,...,s_{2N-1}}}\mathcal{Z}_{\hat\alpha}\big(\bbH,\psi_\eta(s_1),...,\psi_\eta(s_{2N-1}),1\big) \cdot\\&\LF_\bbH^{(\gamma,s_1),...,(\gamma,s_{2N-1}),(\beta,1),(\beta,0),(\beta,\infty)}(dY)ds_1...ds_{2N-1}\bigg]\cdot\psi_\eta'(1)^{1-\Delta_\beta}\,\SLE_\kappa(d\eta).  
    \end{split}
\end{equation}

\emph{Step 4: Change the insertion from $\gamma$ to $\beta$.} We weight the law of $(s_1,...,s_{2N-1},Y,\eta)$ from~\eqref{eq:pf-induc-7} by 
$$\frac{\gamma c_N}{2}  \prod_{k=1}^{2N-1}\big(\e^{\frac{\beta^2-\gamma^2}{4}} e^{\frac{\beta-\gamma}{2}X_\e(x_k)}\big),$$ where $X$ is given by~\eqref{eq:pf-N=2-3} and $x_k=\psi_\eta(s_k)$. Following Step 4 of the proof of Proposition~\ref{prop:N=2}, as we send $\e\to 0$,  by comparing with~\eqref{eq:pf-induc-4}, the conformal welding of $\tilde S_-^f$ and $\tilde S_+^f$ agree in law with  the conformal welding of $S_-^f$ and $S_+^f$. Meanwhile, the law of $(s_1,...,s_{2N-1},Y,\eta)$ converges weakly to 
\begin{equation}\label{eq:pf-induc-8}
    \begin{split}
        \frac{c_N\gamma^2}{4c(Q-\beta)^2} &\int_{0<s_1<...<s_{2N-1}<1}\bigg[\mathds{1}_{E_{\eta,s_1,...,s_{2N-1}}}\mathcal{Z}_{\hat\alpha}\big(\bbH,\psi_\eta(s_1),...,\psi_\eta(s_{2N-1}),1\big)\cdot \prod_{k=1}^{2N-1}\psi_\eta'(s_k)^{1-\Delta_\beta} \\&\cdot\psi_\eta'(1)^{1-\Delta_\beta}\cdot\LF_\bbH^{(\beta,s_1),...,(\beta,s_{2N-1}),(\beta,1),(\beta,0),(\beta,\infty)}(dY)\,ds_1...ds_{2N-1}\bigg]\SLE_\kappa(d\eta).  
    \end{split}
\end{equation}
From our construction in Section~\ref{subsec:pre-msle}, if we first sample $\eta$ from $\mu_\bbH(0,\infty)$ and weight its law by 
$$\mathds{1}_{E_{\eta,s_1,...,s_{2N-1}}}\mathcal{Z}_{\hat\alpha}\big(\bbH,\psi_\eta(s_1),...,\psi_\eta(s_{2N-1}),1\big)\cdot \prod_{k=1}^{2N-1}\psi_\eta'(s_k)^{1-\Delta_\beta} \cdot\psi_\eta'(1)^{1-\Delta_\beta},$$
and sample $(\eta_1,...,\eta_N)$ from $\mSLE_{\kappa,\hat\alpha}(D_\eta;s_1,...,s_{2N-1},1)^\#$ (here we used the conformal covariance~\eqref{eq:partition-conformal-conf}), then the joint law of $(\eta,\eta_1,...,\eta_N)$ is the $\mSLE_{\kappa,\alpha}(\bbH;0,s_1,...,s_{2N-1},1,\infty)$ restricted to the event $E_1$ that $\eta$ is a good link. This proves Theorem~\ref{thm:main} 
when restricting to the event $E_1$. By the rotation symmetry in Lemma~\ref{lm:rotation},  Theorem~\ref{thm:main} 
extends to general $\alpha\in\LP_{N+1}$ when restricted to the event $E_k$ where the link starting from $k$ in $\alpha$ is a good link on both sides of the equation. 

\emph{Step 5. Remove the extra constraint $E_1$ and conclusion.} Let $(D,\eta_1,...,\eta_{N+1};x_1,...,x_{2N+2})$ be an embedding of the spine of $\Wd_{\alpha}(\GQD^{N+2})$ restricted to the event $E$. We first show that $(\eta_1,...,\eta_{N+1})\in X_\alpha(D;x_1,...,x_{2N+2})$. For $1\le k\le N+1$, let $\eta_k^L$ (resp.\ $\eta_k^R$) be the left (resp.\ right) boundary of $\eta_k$. If we only work on the conformal welding of the two generalized quantum disks whose interface is $\eta_k$, then by Theorem~\ref{thm:fd-disk}, $\eta_k$ is an $\SLE_\kappa$ process in the spine of a weight $\frac{3\gamma^2}{2}-2>\frac{\gamma^2}{2}$ forested quantum disk (with a number of marked points on the boundary). In particular, since the left boundary of an $\SLE_\kappa$ from $\mu_\bbH(0,\infty)$ would not touch $(0,\infty)$, this implies that $\eta_k^L\cap\eta_k^R\cap\partial D=\emptyset$, from which we further deduce that $(\eta_1,...,\eta_{N+1})\in X_\alpha(D;x_1,...,x_{2N+2})$. On the other hand, from a simple induction, one can show that if  $(\eta_1,...,\eta_{N+1})\in X_\alpha(D;x_1,...,x_{2N+2})$, there must exist some $1\le k\le N+1$ such that $\eta_k$ is a good link. Therefore $E = \cup_{k=1}^{2N+2} E_k$ and we conclude the proof.
\end{proof}

\subsection{Finiteness of multiple-SLE partiton function and consequences}
In this section we work on the induction step of Theorems~\ref{thm:existence-uniqueness} and~\ref{thm:partition-func}. In particular, we aim to prove the following.
\begin{proposition}\label{prop:induc-finite}
    In the setting of Proposition~\ref{prop:induction}, the measure $\mSLE_{\kappa,\alpha}$ is finite for a.e.\ $\underline{x}\in\mathfrak{X}_{2N+2}$. If $\mathcal{Z}_\alpha(\bbH;x_1,...,x_{2N+2})$ is the size of $\mSLE_{\kappa,\alpha}(\bbH;x_1,...,x_{2N+2})$, then $\mathcal{Z}_\alpha$ is lower semicontinuous and is in the space $L_{\mathrm{loc}}^1(\mathfrak{X}_{2N+2})$. Moreover, if $\mathcal{Z}_\alpha(\bbH;x_1,...,x_{2N+2})<\infty$, then $\mSLE_{\kappa,\alpha}(\bbH;x_1,...,x_{2N+2})^\#$ satisfies the resampling property as in Definition~\ref{def:msle}.
\end{proposition}

\begin{lemma}\label{lem:converge-prob}
    Let $\gamma\in(\sqrt 2,2)$, $N\ge2$ and $\beta_0,...,\beta_{2N-1}<Q$. Let $\tilde{\mathfrak{X}}_{2N} = \{(x_1,...,x_{2N-3})\in \bbR^{2N-3}:0<x_1<...<x_{2N-3}<1\}$. Let $\underline x^n,\underline x\in \tilde{\mathfrak{X}}_{2N}$ with $\underline  x^n\to\underline x$. For $\underline x = (x_1,...,x_{2N-3})\in \tilde{\mathfrak{X}}_{2N}$, set
    \begin{equation}\label{eq:converge-prob}
        f_{\underline x}(z) = \frac{1}{2}\sum_{j=0}^{2N-1}\beta_jG_\bbH(x_j,z)-2Q\log|z|_+,
    \end{equation}
    where $(x_0,x_{2N-2},x_{2N-1})=(0,1,\infty)$.
    For $h\sim P_\bbH$, let $\phi_{\underline x}^0 = h+f_{\underline x}$, and $$I_{\underline x} = \big( \nu_{\phi_{\underline x}^0}((-\infty,0)),\nu_{\phi_{\underline x}^0}((0,x_1)),...,\nu_{\phi_{\underline x}^0}((x_{2N-3},1)),\nu_{\phi_{\underline x}^0}((1,\infty))  \big). $$
    Then we have $I_{\underline x^n}\to I_{\underline x}$ in probability.
\end{lemma}
\begin{proof}
Since on $|f_{\underline x^n}(z)-f_{\underline x}(z)|$ converges uniformly to 0 on $(-\infty,0)\cup (1,\infty)$, it is clear that $$(\nu_{\phi_{\underline x^n}^0}((-\infty,0)), \nu_{\phi_{\underline x^n}^0}((1,\infty)))\to (\nu_{\phi_{\underline x}^0}((-\infty,0)), \nu_{\phi_{\underline x}^0}((1,\infty)))$$ almost surely, and this convergence extends to intervals with positive distance from $\underline x$ as well. Therefore it suffices to show that for any $1\le k\le 2N-3$ and $\e>0$, there exists some $N_0,\delta>0$ such that, for all $n>N_0$ we have $\nu_{\phi_{\underline x^n}^0}((x_k-\delta,x_k+\delta))<\e$ with probability at least $1-\e$. To see this, we pick $\delta$ such that $\nu_{\phi_{\underline x}^0}((x_k-3\delta,x_k+3\delta))<\frac{\e}{4}$ with probability at least $1-\frac{\e}{2}$. Define $g_n(z) = \frac{x_k}{x_k^n}z$,  $\tilde f^n = f_{\underline x^n}(z)\circ g_n^{-1}+Q\log|(g_n^{-1})'|$ and $\tilde \phi^n = h\circ g_n^{-1}+\tilde f^n$. Then for sufficiently large $n$,  $\nu_{\phi_{\underline x^n}^0}((x_k-\delta,x_k+\delta)) \le  \nu_{\tilde\phi^n}((x_k-3\delta,x_k+3\delta))$, and $|\tilde f^n(z)-f_{\underline x}(z)|<\frac{1}{100}$ on $(x_k-3\delta,x_k+3\delta)$. Moreover, by our choice of normalization and the conformal invariance of the unnormalized GFF, $h\circ g_n^{-1}$ has the same distribution as $h-h_{\frac{x_k}{x_k^n}}(0)$, where $h_{\frac{x_k}{x_k^n}}(0)$ is the average of $h$ on the semicircle $\{z\in\bbH:|z| = \frac{x_k}{x_k^n}\}$. On the other hand,   $h_{\frac{x_k}{x_k^n}}(0)\to 0$ in probability as $n\to\infty$. Therefore it follows that for sufficiently large $n$, 
$$\bbP(\nu_{\tilde\phi^n}((x_k-3\delta,x_k+3\delta))>\e)\le \bbP(\nu_{h+f_{\underline x}(z)+\frac{1}{100}-h_{\frac{x_k}{x_k^n}}(0)}((x_k-3\delta,x_k+3\delta))>\frac{\e}{4})\le \e $$
 and the claim follows.   
\end{proof}

Now we prove the following lemma on the boundary lengths of Liouville fields.
\begin{lemma}\label{lem:unif-bd}
    Let $\gamma\in(\sqrt 2,2)$, $N\ge2$ and $\beta = \frac{4}{\gamma}-\frac{\gamma}{2}$. Let $\tilde{\mathfrak{X}}_{2N} = \{(x_1,...,x_{2N-3})\in \bbR^{2N-3}:0<x_1<...\\<x_{2N-3}<1\}$ and $K$ be a compact subset of $\tilde{\mathfrak{X}}_{2N}$. Then there exists a constant $p_0>0$ such that for all $(x_1,...,x_{2N-3})\in K$, the following holds.  For $\phi$ sampled from the measure $\LF_\bbH^{(\beta,0),(\beta,1),(\beta,\infty),(\beta,x_1),...,(\beta,x_{2N-3})}$, we forest the boundary of $(\bbH,\phi,0,x_1,...,x_{2N-3},1)$. The event $F_0$ where the generalized quantum lengths of all the forested boundary segments between $(-\infty,0),(0,x_1),...,(x_{2N-3},1),(1,\infty)$ are in $[1,2]$ 
    has measure at least $p_0$.
\end{lemma}

\begin{proof}
By the continuity of the constant $C_{\mathbb{H}}^{(\beta_i, s_i)_i}$ over $s_i$ in Definition~\ref{def-lf-H-bdry}, it suffices to show that for $(h,\mathbf{c})\sim P_\mathbb{H}\times [e^{(N\beta - Q)c}dc]$, if we let $\phi_{\underline x} = h+f_{\underline x}+\mathbf{c} = \phi_{\underline x}^0+\mathbf{c}$, then the event $F_0$ has measure at least $p_0$ for any $\underline x\in K$, where $f_{\underline x}$ is defined in~\eqref{eq:converge-prob}. 

Now let $\underline x^n,\underline x\in \tilde{\mathfrak{X}}_{2N}$ with $\underline  x^n\to\underline x$. Let $G(\underline x)$ be the measure of the event  {$F_0$} for $\phi_{\underline x}$ under $ P_\mathbb{H}\times [e^{(N\beta - Q)c}dc]$. It is easy to check that $G(\underline x)>0$ for any fixed $\underline x$,  {since by Lemma~\ref{lem:law-forested-line}, for each given instance of $\nu_{\phi_{\underline x}}((-\infty,0)), ..., \nu_{\phi_{\underline x}}((1,\infty))$, the law of the generalized quantum lengths of the forested boundary segments are described by independent stable subordinators.} 
Recall that a positive lower semicontinuous function on a compact set has uniform lower bound. Therefore it suffices to prove $G(\underline x)\le\liminf_{n\to\infty} G(\underline x^n).$ Following Lemma~\ref{lem:law-forested-line}, for $\ell>0$, let $g(\ell)$ be the probability that a forested line of segment $\ell$ has generalized quantum length between $[1,2]$. Then  {$g(\ell) = \bbP(Y_\ell\in[1,2]) = \bbP(\ell^{4/\gamma^2}Y_1\in[1,2])$, where $(Y_t)_{t\geq0}$ is a stable subordinator of index $\frac{\gamma^2}{4}$, and clearly} $g(\ell)$ is continuous in $\ell$. By Fatou's lemma and Lemma~\ref{lem:converge-prob} we have
\begin{equation*}
\begin{split}
    0& <G(\underline x) = \int_0^\infty \int \bigg[g\big(e^{\frac{\gamma}{2}c}\nu_{\phi_{\underline x}}((-\infty,0))\big)g\big((e^{\frac{\gamma}{2}c}\nu_{\phi_{\underline x}}((0,x_1))\big)\cdots g\big(e^{\frac{\gamma}{2}c}\nu_{\phi_{\underline x}}((1,\infty))\big)\bigg]P_\bbH(dh)  e^{(N\beta - Q)c}dc\\
    &\le \int_0^\infty \liminf_{n\to\infty}\int \bigg[g\big(e^{\frac{\gamma}{2}c}\nu_{\phi_{\underline x^n}}((-\infty,0))\big)g\big((e^{\frac{\gamma}{2}c}\nu_{\phi_{\underline x^n}}((0,x_1))\big)\cdots g\big(e^{\frac{\gamma}{2}c}\nu_{\phi_{\underline x^n}}((1,\infty))\big)\bigg]P_\bbH(dh)  e^{(N\beta - Q)c}dc\\
    &\le \liminf_{n\to\infty}\int_0^\infty\int \bigg[g\big(e^{\frac{\gamma}{2}c}\nu_{\phi_{\underline x^n}}((-\infty,0))\big)g\big((e^{\frac{\gamma}{2}c}\nu_{\phi_{\underline x^n}}((0,x_1))\big)\cdots g\big(e^{\frac{\gamma}{2}c}\nu_{\phi_{\underline x^n}}((1,\infty))\big)\bigg]P_\bbH(dh)  e^{(N\beta - Q)c}dc\\ &=\liminf_{n\to\infty} G(\underline x^n)
    \end{split}
\end{equation*}
and the claim follows.

\end{proof}

\begin{proof}[Proof of Proposition~\ref{prop:induc-finite}]

Draw a planar partition of $\bbH$ according to the link pattern $\alpha\in\LP_{N+1}$, and let $\mathcal{S}_\alpha$ be the collection of the domains of this partition. We label the links by $1,...,N+1$, and the boundary segments by $N+2,...,3N+3$. For $D\in \mathcal{S}_\alpha$, let $\mathcal{I}_D$ be the set of indices of the links on $\partial D$, and $\mathcal{B}_D$ be the  set of indices of the boundary segments on $\partial D$.

For the conformal welding $\Wd_\alpha(\GQD^{N+1})$, on the event $E$, let $F_0$ be the event where the generalized quantum lengths of all the $2N+2$ boundary segments are in $[1,2]$. Then following the definition and Proposition~\ref{prop:gqd-bdry-length},  we have  for some constant $C>0$ depending on $\kappa$ and $N$,
\begin{equation*}
\begin{split}
    &\Wd_\alpha(\GQD^{N+1})[E\cap F_0]\\&\le C\int_{s_{N+2},...,s_{3N+3}\in[1,2]}\int_{\ell_1,...,\ell_{N+1}\in (0,\infty)} \prod_{D\in\mathcal{S}_\alpha}\big(\sum_{j\in \mathcal{B}_D}s_j+\sum_{i\in\mathcal{I}_D}\ell_i\big)^{-\frac{\gamma^2}{4}-1}d\ell_1...d\ell_{N+1}ds_{N+2}...ds_{3N+3}
    \\& \le C \int_{\ell_1,...,\ell_{N+1}\in (0,\infty)} \prod_{D\in\mathcal{S}_\alpha}\big(1+\sum_{i\in\mathcal{I}_D}\ell_i\big)^{-\frac{\gamma^2}{4}-1}d\ell_1...d\ell_{N+1}.
    \end{split}
\end{equation*}
It is easy to show that there exists an injective map $f$ from $\{1,...,N+1\}$ to $\mathcal{S}_\alpha$ such that for each $1\le i\le N+1$, $i\in\mathcal{I}_{f(i)}$, i.e., one can assign each interface $\eta_i$ to a unique domain with $\eta_i$ on the boundary. Therefore
\begin{equation}\label{pf:finite-2}
\begin{split}
      \int_{\ell_1,...,\ell_{N+1}\in (0,\infty)
      } \prod_{D\in\mathcal{S}_\alpha}\big(1+\sum_{i\in\mathcal{I}_D}\ell_i\big)^{-\frac{\gamma^2}{4}-1}d\ell_1...d\ell_{N+1}\le \int_{\bbR_+^{N+1}}\prod_{i=1}^{N+1}(\ell_i+1)^{-\frac{\gamma^2}{4}-1}d\ell_1...d\ell_{N+1}<\infty.
    \end{split}
\end{equation}
If we apply~\eqref{pf:finite-2} to the expression~\eqref{eq:thm-main}, we observe that the integral
\begin{equation}\label{pf:finite-3}
\begin{split}
     \int_{0<y_1<...<y_{2N-1}<1}\bigg[\LF_\bbH^{(\beta,0),(\beta,1),(\beta,\infty),(\beta,y_1),...,(\beta,y_{2N-1})}[F_0]\times \big|\mathrm{mSLE}_{\kappa,\alpha}(\bbH,0,y_1,...,y_{2N-1},1,\infty)\big|\bigg]dy_1...dy_{2N-1}
    \end{split}
\end{equation}
is finite. By Lemma~\ref{lem:unif-bd}, if we set $\mathcal{Z}_\alpha(\bbH;0,y_1,...,y_{2N-1},1,\infty) = \big|\mathrm{mSLE}_{\kappa,\alpha}(\bbH,0,y_1,...,y_{2N-1},1,\infty)\big|$, then the function $(y_1,...,y_{2N-1})\mapsto \mathcal{Z}_\alpha(\bbH;0,y_1,...,y_{2N-1},1,\infty)$ is in the space ${L}_{\mathrm{loc}}^1(\tilde{\mathfrak{X}}_{2N+2})$. Then from the conformal covariance property in Lemma~\ref{lm:rotation}, the measure $\mSLE_{\kappa,\alpha}$ can be extended to any polygons, and the function $(x_1,...,x_{2N+2})\mapsto \mathcal{Z}_\alpha(\bbH;x_1,...,x_{2N+2})$ is in the space $L^1_{\loc}(\mathfrak{X}_{2N+2})$. 

To prove that $\mathcal{Z}_\alpha$ is lower semicontinuous, we first assume that $\{1,2N+2\}$ is in the link pattern $\alpha$. For an  $\SLE_\kappa$ curve $\eta_1$ in $\bbH$ from 0 to $\infty$, let $\mathcal{Z}_{\hat{\alpha}}(\hat\bbH_{\eta_1};y_1,...,y_{2N-1},1)$ be $\mathds{1}_{\mathcal{E}_{\eta_1}}$ times the expression~\eqref{eq:def-msle-b} for $N+1$ and $(D;x_1,...,x_{2N+2})=(\bbH,0,y_1,...,y_{2N-1},1,\infty)$, where $\mathcal{E}_{\eta_1}$ is the event defined above~\eqref{eq:def-msle-a}. Then from the construction in Section~\ref{subsec:pre-msle}, $\mathcal{Z}_\alpha(\bbH;0,y_1,...,y_{2N-1},1,\infty)$ is equal to the expectation of $\mathcal{Z}_{\hat{\alpha}}(\hat\bbH_{\eta_1};y_1,...,y_{2N-1},1)$. Moreover, since the probability of the $\SLE_\kappa$ curve hitting a given boundary marked point is 0, from the induction hypothesis that $\mathcal{Z}_{\alpha_1}$ is smooth when $\alpha_1\in\bigsqcup_{k=1}^N\LP_k$, one can infer that  $\mathcal{Z}_{\hat{\alpha}}(\hat\bbH_{\eta_1};y_1^m,...,y_{2N-1}^m,1)\to\mathcal{Z}_{\hat{\alpha}}(\hat\bbH_{\eta_1};y_1,...,y_{2N-1},1)$ a.s.\ as $(y_1^m,...,y_{2N-1}^m)\to (y_1,...,y_{2N-1})$. Fatou's lemma thus implies that $(y_1,...,y_{2N-1})\mapsto \mathcal{Z}_\alpha(\bbH;0,y_1,...,y_{2N-1},1,\infty)$ is lower semicontinuous in $\tilde{\mathfrak{X}}_{2N+2}$, and $\mathcal{Z}_\alpha$ is lower semicontinuous by conformal covariance. The other cases follow analogously by the conformal covariance.


Finally, for $(\eta_1,...,\eta_{N+1})$ sampled from $\mSLE_{\kappa,\alpha}(D;x_1,...,x_{2N+2})^\#$,  given any $1\le k\le N+1$ and $\eta_k$, from the construction the law of the rest of the $N$ curves are $\mSLE_{\kappa,\alpha_{\tilde{D}}}^\#$ in the corresponding domain $\tilde D$'s of $D\backslash\eta_k$. Therefore from the induction hypothesis the resampling properties immediately follow and we conclude the proof.
\end{proof}

To show that $\mathcal{Z}_\alpha$ is smooth in $\mathfrak{X}_{2N+2}$, we use a martingale property along with a hypoellipticity argument. This proof strategy is outlined in~\cite[Lemma B.4]{peltola2019toward}, where a brief proof sketch is given. 
For notational simplicity assume $\{1,2\}\in\alpha$; the same proof works if we replace the number 2 by any $3\le k\le 2N+2$. Let $\hat\alpha$ be the link pattern obtained by removing $\{1,2\}$ from $\alpha$. Let $\eta$ be an $\SLE_\kappa$ in $\bbH$ from $x_1$ to $x_2$. Recall the  {notations} $\hat{\bbH}_\eta$ and $\mathcal{Z}_{\hat\alpha}(\hat\bbH_\eta;x_3,...,x_{2N+2})$ from~\eqref{eq:def-msle-a}. We parameterize $\eta$ via the Loewner equation~\eqref{eq:def-sle} and let $ {(W_t)_{t\geq 0}}$ be the driving function. Let $\tau_\e = \inf\{t>0:|g_t(x_2)-W_t|=\e\}$. Then thanks to the domain Markov property of the chordal $\SLE_\kappa$ and the conformal covariance of $\mathcal{Z}_{\hat\alpha}$, for  $\underline x\in\mathfrak{X}_{2N+2}$ with $\mathcal{Z}_\alpha(\bbH;\underline x)<\infty$ (which is a.e.\ by Proposition~\ref{prop:induc-finite}), 
\begin{equation}\label{pf:finite-4}
    \begin{split}
        M_{t\wedge\tau_\e} &:= \bbE\big[\mathcal{Z}_{\hat\alpha}(\hat\bbH_\eta;x_3,...,x_{2N+2})\,|\,\eta([0,t\wedge\tau_\e]) \big]\\
        &= \prod_{k=3}^{2N+2}g_{t\wedge\tau_\e}'(x_k)^{b} \bbE\big[\mathcal{Z}_{\hat\alpha}(g_{t\wedge\tau_\e}(\hat\bbH_\eta);g_{t\wedge\tau_\e}(x_3),...,g_{t\wedge\tau_\e}(x_{2N+2}))\,|\,\eta([0,t\wedge\tau_\e]) \big]\\
        &=\prod_{k=3}^{2N+2}g_{t\wedge\tau_\e}'(x_k)^{b}\times (g_{t\wedge\tau_\e}(x_2)-W_{t\wedge\tau_\e})^{2b}\times\mathcal{Z}_\alpha\big(\bbH;W_{t\wedge\tau_\e},g_{t\wedge\tau_\e}(x_2),...,g_{t\wedge\tau_\e}(x_{2N+2})\big) 
    \end{split}
\end{equation}
defines a martingale.  
Indeed, to sample an $\SLE_\kappa$ curve from $x_1$ to $x_2$ and weight its law by $\mathcal{Z}_{\hat\alpha}(\hat\bbH_\eta;x_3,...,x_{2N+2})$, one may (i) sample $\eta|_{[0,t\wedge\tau_\e]}$ 
, (ii) sample an $\SLE_\kappa$ $\eta'$ from $\eta(t\wedge\tau_\e)$ to $x_2$ in $\bbH\backslash\eta([0,t\wedge\tau_\e])$, and (iii) weight the law of $\eta'$ by $\mathcal{Z}_{\hat\alpha}(\hat\bbH_\eta;x_3,...,x_{2N+2})$. By conformal covariance, (ii) and (iii) can be replaced by (ii') sample an $\SLE_\kappa$ $\tilde\eta'$ in $\bbH$ from $W_{t\wedge\tau_\e}$ to $\infty$ and (iii') weight the law of $\tilde\eta'$ by $\prod_{k=3}^{2N+2}g_{t\wedge\tau_\e}'(x_k)^{b}\mathcal{Z}_{\hat\alpha}(g_{t\wedge\tau_\e}(\hat\bbH_\eta);\\g_{t\wedge\tau_\e}(x_3),...,g_{t\wedge\tau_\e}(x_{2N+2}))$ and set $\eta' = g_{t\wedge\tau_\e}^{-1}\circ\tilde\eta'$. This justifies the second line of~\eqref{pf:finite-4}, while the third line follows from the definition of $\mathcal{Z}_\alpha$.  Let $X_t = (W_t,g_t(x_2),...,g_t(x_{2N+2}),g_t'(x_2),...,g_t'(x_{2N+2}))$. 
 Then following the SLE coordinate changes~\cite{SW05}, $\eta$ evolves as an $\SLE_\kappa(\kappa-6)$ process from $x_1$ to $\infty$ with the force point located at $x_2$, i.e., $(X_t)_{t\ge0}$ solves
 \begin{equation}\label{eq:finite-SDE}
     dW_t = \sqrt{\kappa}dB_t+\frac{\kappa-6}{W_t-g_t(x_2)}dt; \ \ dg_t(x_j) = \frac{2dt}{g_t(x_j)-W_t}; \ \ dg_t'(x_j) = -\frac{2g_t'(x_j)dt}{ {(g_t(x_j)-W_t)^2}}
 \end{equation}
 where $j=2,...,{2N+2}$ and $(B_t)_{t\ge0}$ is a standard Brownian motion.
 The infinitesimal generator  of $(X_t)_{t\ge0}$, when acting on smooth functions, is
 \begin{equation}\label{eq:def-A}
     A = \frac{\kappa}{2}\partial_1^2+\frac{\kappa-6}{x_1-x_2}\partial_1+\sum_{j=2}^{2N+2}\bigg(\frac{2}{x_j-x_1}\partial_j-\frac{2y_j}{(x_j-x_1)^2}\partial_{2N+1+j}  \bigg). 
 \end{equation}
 Consider the function $F$ defined by 
 \begin{equation}\label{eq:martingale-F}
    F(x_1,...,x_{2N+2};y_2,...,y_{2N+2}) = \prod_{k=3}^{2N+2}y_k^{b}\times (x_2-x_1)^{2b}\times\mathcal{Z}_\alpha\big(\bbH;x_1,...,x_{2N+2}\big).
\end{equation}
 {By Proposition \ref{prop:induc-finite},} $F$ is a locally integrable and lower semicontinuous function of  $(\underline{x},\underline{y}):=(x_1,...,x_{2N+2},\\y_2,...,y_{2N+2})\in \mathfrak{X}_{2N+2}\times \bbR_+^{2N+1}$. Let $\tau$ be the first time when $\eta$ hits $  [x_2,\infty)$, i.e., 
 \begin{equation}\label{eq:stop-time-tau}
      \tau=\inf\{t>0:g_t(x_2)=W_t\}.
 \end{equation}
 Then it follows from~\eqref{pf:finite-4} by sending $\e\to0$ that $\{F(X_{t\wedge\tau})\}_{t\ge0}$ is a martingale. Moreover, since $dg_t'(x) = -\frac{2g_t'(x)}{(g_t(x)-W_t)^2}dt$, for $a_2,...,a_{2N+2}>0$, if we let $\wt{X}_t = (W_t,g_t(x_2),...,g_t(x_{2N+2}), a_2g_t'(x_2),...,a_{2N+2}g_t'(x_{2N+2}))$, then $(\wt{X}_t)_{t\ge0}$ solves the~\eqref{eq:finite-SDE} as well and starts from $(x_1,...,x_{2N+2}, a_2,...,a_{2N+2})$. We infer the following.
 \begin{lemma}\label{lm:finite-mg}
     For a.e.\ $(\underline{x}^0,\underline{y}^0)=(x_1^0,...,x_{2N+2}^0,y_2^0,...,y_{2N+2}^0)\in\mathfrak{X}_{2N+2}\times\bbR_+^{2N+1}$ the following holds. Let $(X_t)_{t\ge0}$ be a solution of~\eqref{eq:finite-SDE} starting from $(\underline{x}^0,\underline{y}^0)$, and let $\tau$ be defined as~\eqref{eq:stop-time-tau}. Then $F(X_{t\wedge\tau})_{t\ge0}$ is a martingale. 
 \end{lemma}

  To prove the smoothness of $\mathcal{Z}_\alpha$, the first step is to use the martingale property to prove that $F$ is a distributional solution to the differential equation $AF=0$, then use the hypoellipticity of the differential operator $A$ to prove that $F$ is smooth. Recall that a differential operator $\mathfrak{D}$ is hypoelliptic on domain $U\subset\bbR^n$ if for any open set $\cO\subset U$, for any $f\in (C_c^\infty)^*(\cO)$, $\mathfrak Df\in C^\infty(\cO)$ implies $f\in C^\infty(\cO)$. For  smooth vector fields $X_j:=\sum_{k=1}^na_{jk}(x)\partial_k$ on $U$ where $j=0,...,m$ and $a_{jk}$ are smooth functions in $U$, consider the differential operator 
  \begin{equation}\label{eq:hormander-op}
     \mathfrak D = \sum_{j=1}^m X_j^2+X_0+b
  \end{equation}
   where $b\in C^\infty(U)$.  From~\cite{Hor67}, if the Lie algebra generated by $X_0,...,X_m$ has full rank at every point $x\in U$, then $\mathfrak D$ is hypoelliptic.

\begin{lemma}\label{lm:hypoelliptic}
    The operator $A$ defined in~\eqref{eq:def-A} is hypoelliptic.
\end{lemma}
  
 Let $A^*$ be the dual operator of $A$, i.e.,
 \begin{equation*}
     A^*g = \frac{\kappa}{2}\partial_1^2g-\partial_1\big(\frac{(\kappa-6)}{x_1-x_2}g\big)-\sum_{j=2}^{2N+2}\bigg(\partial_j\big(\frac{2}{x_j-x_1}g\big)-\partial_{2N+1+j}\big(\frac{2y_j}{(x_j-x_1)^2}g\big)  \bigg) 
 \end{equation*}
 for smooth function $g$.

 \begin{proposition}\label{prop:weaksol}
      $F$ is a distributional solution to $AF=0$, i.e., $\inner{F}{A^*g}=0$ for any test function $g\in C_c^\infty(\mathfrak{X}_{2N+2}\times \bbR_+^{2N+1})$.
 \end{proposition}
 {Proposition~\ref{prop:weaksol} is a consequence of Lemma~\ref{lm:finite-mg}.} 
 Lemma~\ref{lm:hypoelliptic} and Proposition~\ref{prop:weaksol} shall be proved in Appendix~\ref{sec:pf-weaksol}. Similar statements are considered in~\cite{dubedat2015sle,peltola2019global}.

\begin{proposition}\label{prop:smooth}
    The function $\mathcal{Z}_\alpha(\bbH;x_1,...,x_{2N+2})$ is smooth in $\mathfrak{X}_{2N+2}$, 
    and solves the PDE~\eqref{eq:msle-pde}.
\end{proposition}
\begin{proof}
    Let $\mathfrak{D}^{(1)} = \frac{\kappa}{2}\partial_1^2+\sum_{j=2}^{2N+2}(\frac{2}{x_j-x_1}\partial_j-\frac{2b}{(x_j-x_1)^2})$ be the differential operator in~\eqref{eq:msle-pde} for $i=1$. By Lemma~\ref{lm:hypoelliptic}, Proposition~\ref{prop:weaksol} and~\cite[Theorem 1.1]{Hor67},  $F$ is smooth in $\mathfrak{X}_{2N+2}\times\bbR_+^{2N+1}$. Therefore $\mathcal{Z}_\alpha(\bbH;x_1,...,x_{2N+2})$ is smooth in $\mathfrak{X}_{2N+2}$. A direct computation shows that $AF = \mathfrak{D}^{(1)}\mathcal{Z}_\alpha=0$, i.e.,  $\mathcal{Z}_\alpha\big(\bbH;x_1,...,x_{2N+2})$ solves the PDE~\eqref{eq:msle-pde} for $N+1$ and $i=1$.  {The equation for $i=2$ follows from the reversibility of $\SLE_\kappa$~\cite{IGIII} and an indentical argument as we swap the role of $x_1$ and $x_2$. The equation for other $i\geq 3$ follows from the identical argument as the $i=1,2$ case.} 
\end{proof}

\begin{proof}[Proof of Theorems~\ref{thm:existence-uniqueness}, ~\ref{thm:partition-func} and~\ref{thm:main}]
    For $N=2$, Theorem~\ref{thm:partition-func} and Theorem~\ref{thm:existence-uniqueness} hold by~\cite{miller2018connection}, where Theorem~\ref{thm:main} follows from Proposition~\ref{prop:N=2}. Suppose the theorems hold for $1,...,N$. Then for $N+1$ Theorem ~\ref{thm:main} follows by  Proposition~\ref{prop:induction}. Theorem~\ref{thm:existence-uniqueness} and and the first part of Theorem~\ref{thm:partition-func} follows from Proposition~\ref{prop:induc-finite} and Proposition~\ref{prop:smooth}. The second half of Theorem~\ref{thm:partition-func} follows from Proposition~\ref{prop:global-local} in Section~\ref{subsec:local-msle}.  This concludes the induction step and the whole proof.
\end{proof}

\subsection{Relation with local multiple SLE}\label{subsec:local-msle}
In this section, we show that the global multiple SLE $\mSLE_{\kappa,\alpha}$ agrees with the \emph{local $N$-$\SLE_\kappa$} driven by the partition function $\mathcal{Z}_\alpha$  {as studied in}~\cite{dubedat2007commutation,graham2007multiple,KP16partitionfunc}. To begin with, we recall the definition of local multiple $\SLE_\kappa$ from~\cite{dubedat2007commutation} and~\cite[Appendix A]{KP16partitionfunc}.

Let $(D;x_1,...,x_{2N})$ be a polygon, and $(U_1,...,U_{2N})$ be  localization neighborhoods, in the sense that $x_k\in \ol U_k\subset \ol D$, $D\backslash U_k$ is simply connected and $U_j\cap U_k = \emptyset$ for $1\le j,k\le 2N$ with $j\neq k$. Consider 2N-tuples of oriented unparameterized curves $(\eta_1,...,\eta_{2N})$, where each $\eta_k$ is contained in $\ol U_k$ and connects $x_k$ and a point $x_k'\in \partial U_k$. Choose a parametrization such that $\eta_k:[0,1]\to \ol U_k$ such that $\eta_k(0)=x_k$ and $\eta_k(1) = x_k'$.  A \emph{local $N$-$\SLE_\kappa$} in  $(D;x_1,...,x_{2N})$ localized in $(U_1,...,U_{2N})$, is a probability measure $\mathsf{P}_{(U_1,...,U_{2N})}^{(D;x_1,...,x_{2N})}$ on $(\eta_1,...,\eta_{2N})$ with conformal invariance (CI), domain Markov property (DMP), and absolute continuity of marginals with respect to the chordal $\SLE_\kappa$ (MARG) as follows:

\begin{itemize}
\item[(CI)] If  $(\eta_1,...,\eta_{2N})\sim \mathsf{P}_{(U_1,...,U_{2N})}^{(D;x_1,...,x_{2N})}$, then for any conformal map $f:D\to f(D)$, 
$(f\circ\eta_1,...,f\circ\eta_{2N})\sim \mathsf{P}_{(f(U_1),...,f(U_{2N}))}^{(f(D);f(x_1),...,f(x_{2N}))}$;
    \item[(DMP)] Fix stopping times $(\tau_1,...,\tau_{2N})$ for  $(\eta_1,...,\eta_{2N})$. Given initial segments $(\eta_1|_{[0,\tau_1]},...,\eta_{2N}|_{[0,\tau_{2N}]})$ the conditional law of the remaining parts $(\eta_1|_{[\tau_1,1]},...,\eta_{2N}|_{[\tau_{2N},1]})$ is $\mathsf{P}_{(\wt U_1,...,\wt U_{2N})}^{(\wt D; \wt x_1,...,\wt x_{2N})}$ where for each $1\le k\le 2N$, $\wt x_k$ is the tip $\eta_k(\tau_k)$, $\wt D$ is the connected component of $D\backslash \big(\cup_{k=1}^{2N}\eta_k([0,\tau_k])\big)$ with $\wt x_1,...,\wt x_{2N}$ on the boundary and $\wt U_k = \wt D\cap U_k$. 
    \item[(MARG)] There exist smooth functions $F_j$: $\mathfrak{X}_{2N}\to \bbR$, for $j=1,...,2N$, such that for the domain
$D=\bbH$, boundary points $x_1 <...<x_{2N}$, and their localization neighborhoods $U_1,...,U_{2N}$, the marginal law of $\eta_j$ under $\mathsf{P}_{(U_1,...,U_{2N})}^{(\bbH;x_1,...,x_{2N})}$ is the Loewner evolution driven by $W_t$ which solves 
\begin{equation}\label{eq:local-msle}
    \begin{split}
        &dW_t = \sqrt{\kappa}dB_t + F_j(V_t^1,..., V_t^{j-1},W_t,V_t^{j+1},...,V_t^{2N}); \ \ W_0 = x_j;
        \\& dV_t^k = \frac{2}{V_t^k-W_t}; \ \ V_0^k = x_k\ \ \text{for } k\neq j. 
    \end{split}
\end{equation}
\end{itemize}
 Dub\'{e}dat~\cite{dubedat2007commutation} proved that the local $N$-$\SLE_\kappa$ processes are classified by partition functions as below. We use the following version stated in \cite[Proposition 4.7]{peltola2019global}. 
 \begin{proposition}\label{prop:local-msle}
    Let $\kappa>0$.
    \begin{enumerate}[(i)]
        \item     Suppose $\mathsf P$ is a local $N$-$\SLE_\kappa$. Then there exists a function $\mathcal{Z}:\mathfrak {X}_{2N}\to \bbR_+$ satisfying (PDE)~\eqref{eq:msle-pde}
and (COV)~\eqref{eq:msle-conformal-conf}, such that for all $j=1,...,2N$, the drift functions in (MARG) take the form
$F_j = \kappa\partial_j\log \mathcal{Z}$.  Such a function $\mathcal{Z}$ is determined up to a multiplicative constant.
\item Suppose $\mathcal{Z}:\mathfrak {X}_{2N}\to \bbR_+$ satisfies (PDE)~\eqref{eq:msle-pde}
and (COV)~\eqref{eq:msle-conformal-conf}. Then, the random collection of
curves obtained by the Loewner chain in (MARG) with $F_j = \partial_j\log\mathcal{Z}$, for all $j=1,...,2N$,
is a local $N$-$\SLE_\kappa$. Two functions $\mathcal{Z}$ and $\wt{\mathcal{Z}}$ give rise to the same local $N$-$\SLE_\kappa$ if and only if $\mathcal{Z} = const\times \wt{\mathcal{Z}}$.
    \end{enumerate}
 \end{proposition}
Now we show that the initial segments of our global multiple SLE agrees with the local multiple SLE driven by the parition function $\mathcal{Z}_\alpha$. The argument is almost the same as~\cite[Lemma 4.8]{peltola2019global}, except that we perform truncations since a priori we do not have the strong power law bounds on $\mathcal{Z}_\alpha$ as in the setting there.
\begin{proposition}\label{prop:global-sde}
    Let $\kappa\in(4,8)$ and $\alpha\in\LP_N$. Fix $x_1<...<x_{2N}$. Assume $\{j,k\}\in\alpha$, and suppose $W_t^j$ solves~\eqref{eq:local-msle} with $F_j = \kappa\partial_j\log\mathcal{Z}_\alpha$. Let $$T_j = \inf\{t>0:\min_{i\neq j}|g_t(x_i) - W_t^j| = 0\}.$$
    Then the Loewner equation driven by $W_t^j$ is well-defined up to time $T_j$. Moreover, it
is almost surely generated by a continuous curve up to and including time $T_j$, which has the same law as the curve $\eta_j$ in $\mSLE_{\kappa,\alpha}(\bbH;x_1,...,x_{2N})^\#$ connecting $x_j$ and $x_k$ stopped at the time $\sigma_j$ that it separates any of $\{x_i:i\neq j\}$ from $\infty$. 
\end{proposition}

\begin{proof}
    Let  $(\wt W_t)_{t\ge0}$ be the Loewner driving function for $\eta_j$, and $(\wt g_t)_{t\ge0}$ be the associated Loewner maps. For $\e,M>0$, let $$\wt\tau_{\e,M} = \inf\{t>0:\min_{i\neq j}|\wt g_t(x_i) -\wt W_t| = \e \text{ or } \max_{i\neq j}|\wt g_t(x_i) -\wt W_t| = M\}.$$   Thanks to the domain Markov property and conformal invariance of $\SLE_\kappa$, by~\eqref{pf:finite-4} (with $N+1$ replaced by $N$), $\eta_j|_{[0,\wt\tau_{\e,M}]}$ can be produced by 
    \begin{enumerate}[(i)]
        \item Sample an $\SLE_\kappa$ $\eta$  in $\bbH$ from $x_j$ to $x_k$ parameterized via the Loewner equation~\eqref{eq:def-sle}. Let ${(W_t)_{t\geq 0}}$ be the driving function and $(g_t)_{t\ge0}$ be the Lowener maps, then
        \begin{equation}
            dW_t = \sqrt{\kappa}dB_t+\frac{\kappa-6}{W_t - g_t(x_k)}dt.
        \end{equation}
        Let $\tau_{\e,M} = \inf\{t>0:\min_{i\neq j}| g_t(x_i) -W_t| = \e \text{ or } \max_{i\neq j}| g_t(x_i) -W_t| = M \}$, and $\tau_j = \inf\{t>0:\min_{i\neq j}| g_t(x_i) -W_t| = 0\}$. 
        \item Weight the law of $\eta|_{[0,\tau_{\e,M}]}$ by $\mathsf M_{\tau_{\e,M}}$, where $$\mathsf M_t=\frac{1}{Z}\prod_{i\neq j,k}g_t'(x_i)^{b}\times (g_t(x_k)-W_t)^{2b}\times\mathcal{Z}_\alpha\big(\bbH; g_t(x_1),...,g_t(x_{j-1}),W_t,g_t(x_{j+1}),...,g_t(x_{2N})\big) $$
        where $Z = \mathcal{Z}_\alpha(\bbH;x_1,...,x_{2N})$.
    \end{enumerate}
    Then $(\mathsf M_{t\wedge\tau_{\e,M}})_{t\ge0}$ is a martingale for $\eta$. For fixed $T>0$ and $0<t<T\wedge\tau_{\e,M}$, since for $i\neq k$, $d(g_t(x_i)-g_t(x_k)) = -2\frac{g_t(x_i)-g_t(x_k)}{(g_t(x_i)-W_t)(g_t(x_k)-W_t)}$ and $dg_t'(x) = -\frac{2g_t'(x)}{(g_t(x)-W_t)^2}$, one can check that $|x_i-x_k|\ge |g_t(x_i)-g_t(x_k)|\ge |x_i-x_k|e^{-T\e^{-2}}$ and $1\ge|g_t'(x_k)|\ge e^{-2T\e^{-2}}$. This implies that $(W_t,g_t(x_2),...,g_t(x_{2N}))_{0\le t\le \tau_{\e,M}}$ is contained in some fixed compact subset of $\mathfrak{X}_{2N}$ and thus $(\mathsf M_{t\wedge\tau_{\e,M}})_{0\le t\le T}$ is a martingale bounded from both above and below.  Moreover, since $\mathcal{Z}_\alpha$ solves the (PDE)~\eqref{eq:msle-pde}, it follows that
    $$ \mathsf M_t^{-1}d\mathsf M_t = \sqrt{\kappa}\bigg(\partial_j\log\mathcal{Z}_\alpha - \frac{2b}{g_t(x_k)-W_t} \bigg)dB_t.$$
    Therefore by the Girsanov theorem, if we weight the law of $(B_t)_{t\in[0,T\wedge\tau_{\e,M}]}$ by $\mathsf M_{\tau_{\e,M}}$, then $(W_t)_{t\ge0}$ solves~\eqref{eq:local-msle} up until $T\wedge\tau_{\e,M}$. This proves the statement up until the time $T\wedge\tau_{\e,M}$. Since $\kappa>4$, $\tau_j<\infty$ a.s.. Therefore if we send $M,T\uparrow\infty$ and $\e\downarrow0$, we have $(T\wedge\tau_{\e,M})\uparrow \tau_j$ and the claim follows.
\end{proof}
 
By Proposition~\ref{prop:global-sde}, using the domain Markov property, conformal invariance and the reversibility of $\SLE_\kappa$ for $\kappa\in(4,8)$, we have the following.
\begin{proposition}\label{prop:global-local}
     Let $N\ge1$, $\alpha\in\LP_N$, and $(D;x_1,...,x_{2N})$ be a polygon. Suppose $(\eta_1,...,\eta_{2N})$ is a sample from $\mSLE_{\kappa,\alpha}(D;x_1,...,x_{2N})^\#$. Then for any  localization neighborhoods $(U_1,...,U_{2N})$, the law of $(\eta_1,...,\eta_{2N})$ when restricted to $(U_1,...,U_{2N})$ agrees with $\mathsf{P}_{(U_1,...,U_{2N})}^{(D;x_1,...,x_{2N})}$ driven by partition function $\mathcal{Z}_\alpha$ in the sense of Proposition~\ref{prop:local-msle}.
\end{proposition}

\appendix
\section{Proof of Lemma~\ref{lm:hypoelliptic} and Proposition~\ref{prop:weaksol}}\label{sec:pf-weaksol}

In this section we prove Lemma~\ref{lm:hypoelliptic} and Proposition~\ref{prop:weaksol}. The main idea at high level, which is the same as~\cite[Lemma 4.4]{peltola2019global}, is to consider the semigroup generated by the process $(X_t)_{t\ge0}$ and extend the domain of its infinitesimal generator $A$. Then the martingale property would imply that $F$ is a weak solution as in Proposition~\ref{prop:weaksol}. However, there are several obstacles to directly apply the proof of~\cite[Lemma 4.4]{peltola2019global}: (i) as pointed out by Dapeng Zhan, the extension of the $A$ to the space of generalized functions is not clear in the original proof; (ii) the martingale property in~\eqref{pf:finite-4} is only valid up to the stopping time $\tau$; (iii) a priori we only know by Proposition~\ref{prop:induc-finite} that $F$ is lower semicontinuous while in~\cite{peltola2019global} $F$ is assumed to be continuous; (iv) a further restriction in the proof in~\cite{peltola2019global} is that $g_t'(x_j) = 1$ for $t=0$ from the definition of the Loewner evolution and thus the starting point of $(X_t)_{t\ge0}$ is not arbitrary. 

We write down a complete proof to deal these issues. Issue (iv) is already treated in Lemma~\ref{lm:finite-mg} using the homogeneity of the equations for $g_t'(x)$, and we prove in Lemma~\ref{lm:hypoelliptic} that the operator $A$ (rather than the operator $\mathfrak{D}^{(1)}$, as proved in~\cite{peltola2019global}) is hypoelliptic. For (ii), we apply truncation using the stopping time $\sigma$ as below such that Lemma~\ref{lm:finite-mg} is applicable and the terms in~\eqref{eq:finite-SDE} are smooth. Then for (i), in Lemma~\ref{lm:feller} we show that the truncated proecss $(X_{t\wedge\sigma})_{t\ge0}$ is Feller, and use Bony's theorem\footnote{We thank Eveliina Peltola for introducing us to this theorem.}~\cite{bony} along with properties of the infinitesimal generator to rigorously justify the integration by parts in the proof of~\cite[Lemma 4.4]{peltola2019global}. Finally for (iii), we establish the integral equation in Lemma~\ref{lm:inte-by-part-2} and apply the monotone convergence theorem for general lower semicontinuous functions.


\begin{proof}[Proof of Lemma~\ref{lm:hypoelliptic}]
To check the H\"{o}rmander condition, we set $X_0 = \frac{\kappa-6}{x_1-x_2}\partial_1+\sum_{j=2}^{2N+2}\big(\frac{2}{x_j-x_1}\partial_j-\frac{2y_j}{(x_j-x_1)^2}\partial_{2N+1+j}\big)$, and $X_1 = \sqrt{\frac{\kappa}{2}}\partial_1$. We write $X_0^{[0]} = X_0$. Then for $n\ge1$, by induction we have
\begin{equation}\label{eq:liealgebra}
    X_0^{[n]}:=\frac{1}{n}\big[\partial_1,X_0^{[n-1]}  \big] = \frac{6-\kappa}{(x_2-x_1)^{n+1}}\partial_1+\sum_{j=2}^{2N+2}\big(\frac{2}{(x_j-x_1)^{n+1}}\partial_j -\frac{2(n+1)y_j}{(x_j-x_1)^{n+2}}\partial_{2N+1+j}\big)
\end{equation}
Consider the matrix $\mathbf{A} = (a_{ij})_{1\le i,j\le 4N+2}$, where for $1\le j\le 2N+1$, $a_{ij} = \frac{2}{(x_{j+1}-x_1)^i}$, and for $2N+2\le j\le 4N+2$, $a_{ij} = -\frac{2iy_{j-2N}}{(x_{j-2N}-x_1)^{i+1}}$. Indeed, to prove that the linear space spanned by $X_1,X_0^{[(0)]},...,X_0^{[4N+1]}$ has dimension $4N+3$, it suffices to show that $\det\mathbf{A}\neq0$ for every $(x_1,...,x_{2N+2},y_2,...,y_{2N+2})\in\mathfrak{X}_{2N+2}\times\bbR_+^{2N+1}$. Let $\wt{\mathbf{A}} = (\wt a_{ij})_{1\le i,j\le 4N+2}$ where $\wt a_{ij} = \frac{1}{(x_{j+1}-x_1)^i}$ for $1\le j\le 2N+1$ and $\wt a_{ij} = \frac{i-1}{(x_{j-2N}-x_1)^{i-1}}$ for $2N+2\le j\le 4N+2$. Then $\det\mathbf{A}\neq0$ if and only if $\det\wt{\mathbf{A}}\neq0$, and $\wt{\mathbf{A}}$ is a confluent Vandermonde matrix and hence invertible (see e.g.\ Eq.\,(1.3) in~\cite{Walter62}). This concludes the proof.  
\end{proof}

To prove Proposition~\ref{prop:weaksol}, it suffices to show that for any given $(\underline x^0,\underline y^0)\in\mathfrak X_{2N+2}\times\bbR_+^{2N+1}$, there exists some neighborhood $\cO$ such that $\inner{F}{A^*g} = 0$ for any $g\in C_c^\infty(\cO)$. We let $\cO$ be the interior of the convex hull in $\bbR_+^{4N+3}$ generated by $\{(x_1^0+\e_1\delta_0, ..., x_{2N+2}^0+\e_{2N+2}\delta_0, y_2^0+\e_{2N+3}\delta_0,..., y_{2N+2}^0+\e_{4N+3}\delta_0):\e_{1},...,\e_{4N+3} = \pm1 \}$, and choose $\delta_0>0$ sufficiently small such that $\ol\cO\subset\mathfrak X_{2N+2}\times\bbR_+^{2N+1}$. For the It\^{o} process $(X_t)_{t\ge0}$ described by~\eqref{eq:finite-SDE}, let $\sigma = \inf\{t>0:X_t\notin \cO\}$. 
\begin{lemma}\label{lm:feller}
    $(X_{t\wedge\sigma})_{t\ge0}$ is a Feller process in $\ol\cO$.
\end{lemma}
\begin{proof}
    For  $f\in C(\ol\cO)$ and $\mathbf x\in \ol\cO$, define $(P_tf)(\mathbf x) = \bbE\big[f(X_{t\wedge\sigma})|X_0 = \mathbf x \big]$.   We need to show that 
    \begin{enumerate}[(i)]
        \item $P_tf\in C(\ol\cO)$;
        \item As $t\to0$, $P_tf$ converges to $f$ uniformly.
    \end{enumerate}
    By~\cite[Proposition III.2.4]{revuzyorbook}, to prove (ii) it suffices to show that $P_tf(\mathbf x)\to f(\mathbf x)$ for any $\mathbf x\in\ol\cO$, which readily follows from the continuity of the paths of $X$ and dominant convergence theorem. To prove (i), for any given $\e>0$, we pick $\wt f\in   C^\infty(\ol\cO)$ such that $\|f-\wt f\|_{C(\ol\cO)}<\e/3$.  {Since the coefficients for $X_t$ as defined in~\eqref{eq:finite-SDE} are smooth within a neighborhood of $\cO$,  by Dynkin's formula~\cite[Lemma 17.21]{Kal21book}}, for $\mathbf x\in \cO$,
    \begin{equation}\label{eq:pf-feller-1}
        P_t\wt f(\mathbf x) = \wt f(\mathbf x) + \bbE^{\mathbf x} \int_0^{t\wedge\sigma} A\wt f(X_s)ds =    \wt f(\mathbf x) + \int_0^t \bbE^{\mathbf x} 1_{\sigma>s}A\wt f(X_s)ds.
    \end{equation}
    Note that we have implicitly used the fact that $\bbP(\sigma = s) = 0$ for any $s>0$. This is because,  for any $\vec\lambda = (\lambda_1,...,\lambda_{4N+3})$ with $\lambda_1\neq0$, $\bbP(s<\tau;\vec\lambda \cdot X_{s} = \lambda_0) = 0$ for any given $\lambda_0\in\bbR$. One can check this by applying the Girsanov theorem and comparing with the Brownian motion. Moreover, using 
     {the smoothness of the coefficients in $\ol \cO$, as explained in~\cite[Theorem 8.5(ii)]{legallbook}, we have the continuity of solutions with respect to the initial value, in the sense that if $(X_t')_{t\ge0}$ solves~\eqref{eq:finite-SDE} with $X_0' = \mathbf x'$, there exists a coupling such that $\mathds{1}_{\sigma>s}(X_r')_{0\le r\le s}$ converges uniformly in $r$ in probability as $\mathbf x'\to\mathbf x$ to $\mathds{1}_{\sigma>s}(X_r)_{0\le r\le s}$}. In particular, by applying  {the dominated} 
    convergence theorem to~\eqref{eq:pf-feller-1}, $|P_t\wt f(\mathbf x)-P_t\wt f(\mathbf x')|<\e/3$ and thus $|P_tf(\mathbf x)-P_t f(\mathbf x')|<\e$ for $\mathbf x'$ sufficiently close to $\mathbf x$. This proves $P_tf\in C(\cO)$. Now for $\mathbf x\in\partial\cO$, take $\delta\in(0,\e)$
 such that $|f(\mathbf x')-f(\mathbf x)|<\e/3$ when $|\mathbf x'-\mathbf x|<\delta$,  {and without loss of generality assume that the size of the coefficients in~\eqref{eq:finite-SDE} are bounded by some constant $K$ not depending on $\delta$ within the $\delta$-neighborhood of $\ol\cO$}. Without loss of generality, assume that $\mathbf x$ is on the right half of $\partial \cO$, i.e., $\mathbf x+\lambda\mathbf{e}_1\notin\cO$ for any $\lambda>0$ where $\mathbf{e}_1 = (1,0,...,0)$. Let $(B_t)_{t\geq 0}$ be the standard Brownian motion in~\eqref{eq:finite-SDE}. Then we have  {that}  {the event $E_\delta = \{\sup_{0\le s\le \delta^{3}}B_s>\delta^{7/4}\}\cap\{\sup_{0\le s\le \delta^{3}}|B_s|<\delta^{5/4}\}$   has probability $1-o_\e(1)$}. For $|\mathbf x'-\mathbf x|<\delta^2$,  {we claim}  that for $X_t' = ((X_t^1)',...,(X_{t}^{4N+3})')$ starting from $\mathbf x'$, with probability $1-o_\e(1)$, $\sup_{0<s<\delta^3}(X_s^1)'-(X_0^1)'>\delta^{9/5}$ while $\sup_{0<s<\delta^3}|(X_s^j)'-(X_0^j)'|<\delta^{8/3}$ for $j=2,...,4N+3$.  {Indeed, let $\sigma_\delta$ be the first time $X_t'$ exits the $\delta$-neighborhood of $\ol\cO$. Then for $s\in (0,\delta^3\wedge\sigma_\delta)$, from~\eqref{eq:finite-SDE} one has $|(X_s^1)'-(X_0^1)'-\sqrt\kappa B_s|\leq K\delta^3$ and $|(X_s^j)'-(X_0^j)'|<K\delta^3$ for $j\geq2$. On the event $E_\delta$, $\|(X_s^j)'-(X_0^j)'\|\leq c(K\delta^3+\delta^{5/4})$ for some constant $c$, and for sufficiently small $\delta$ we must additionally have $\sigma_\delta>\delta^3$. This verifies our claim.} In particular, by our choice of $\cO$, with probability $1-o_\e(1)$, $X'$ exits $\cO$ before time $\delta^3$ and the exit location is within $\delta$-neighborhood of $\mathbf x$. Therefore $P_tf(\mathbf x')\to f(\mathbf x) = P_tf(\mathbf x)$ for $\mathbf x\in\partial \cO$. This finishes the proof of (i).
 \end{proof}
Let $\wt A$ be the infinitesimal generator of $(X_{t\wedge\sigma})_{t\ge0}$ on $C(\ol\cO)$. Write $D(\wt A)$ for the domain of $\wt A$, i.e., $$D(\wt A) = \{ f\in C(\ol\cO): \lim_{t\to0}\frac{P_tf-f}{t}\ \text{exists in } C(\ol\cO) \}.$$
 \begin{lemma}\label{lm:generator}
     Suppose $f\in C^\infty(\cO)\cap C(\ol\cO)$ such that $Af(\mathbf x_n)\to 0$ as $\mathbf x_n \to\mathbf x\in \partial \cO$. Then $f\in D(\wt A)$, and if we define $Af(\mathbf x) = 0$ for $\mathbf x\in\partial \cO$, then $\wt Af = Af$.
 \end{lemma}
\begin{proof}
    For $\e>0$, let $\sigma_\e = \inf\{t>0: \dist(X_t,\partial \cO)<\e\}$. Then $(t\wedge\sigma_\e)\uparrow (t\wedge\sigma)$ as $\e\to0$. By Dynkin's formula, for $\mathbf x\in \cO$, 
    \begin{equation}
        \bbE^{\mathbf{x}}f(X_{t\wedge\sigma_\e}) = f(\mathbf x)+\bbE^{\mathbf{x}}\big[\int_0^{t\wedge\sigma_\e}Af(X_s)ds \big].
    \end{equation}
    Sending $\e\to0$, by dominant convergence theorem,
    \begin{equation}\label{eq:generator-1}
        P_tf(\mathbf{x}) = f(\mathbf x)+\bbE\big[\int_0^{t\wedge\sigma}Af(X_s)ds \big] = f(\mathbf x)+\bbE^{\mathbf{x}}\big[\int_0^{t}Af(X_{s\wedge\sigma})ds \big]
    \end{equation}
    since $Af = 0$ on $\partial\cO$. By definition, ~\eqref{eq:generator-1} continues to hold when $\mathbf x\in\partial\cO$. Meanwhile, by definition, 
    \begin{equation}\label{eq:generator-2}
        P_tf(\mathbf{x}) = f(\mathbf x)+\bbE^{\mathbf{x}}\big[\int_0^{t}Af(X_{s\wedge\sigma})ds \big] =  f(\mathbf x) + \int_0^t P_sAf(\mathbf{x})ds.
    \end{equation}
     Therefore, since $Af\in C(\ol\cO)$, 
     $$\frac{P_tf - f}{t} = \frac{1}{t}\int_0^t P_sAfds \to Af $$
    in $C(\ol\cO)$ as $t\to 0$, which implies that $f\in D(\wt A)$ and  $Df = Af$.  
\end{proof}

A vector $v$ is called an exterior normal to a closed set $F$ at a point $x_0\in F$, if there exists an open ball contained in $\bbR^n\backslash F$ centered at a point $x_1$, such that $x_0$ belongs to the closure of this ball and  $v = \lambda(x_1 - x_0)$ for some $\lambda>0$. Let $\mathfrak D$ be a differential operator which can be written by~\eqref{eq:hormander-op} and satisfies the H\"{o}rmander condition.  Write $\mathfrak D$ as $\mathfrak D = \sum_{i,j=1}^n a_{ij}(x)\partial_{ij}+ \sum_{j=1}^n b_i(x)\partial_i+c(x)$. Suppose $ O\subset \bbR^n$ is a bounded domain, such that for any $x\in\partial  O$, there exists an exterior normal $v$ vector to $\ol O$ at $x$ with $\sum_{i,j=1}^na_{ij}(x)v_iv_j>0$. Then in~\cite[Theorem 5.2]{bony}, Bony has proved that if $c(x)<c_0<0$ for some constant $c_0$ in $O$, then for continuous functions $f$ and $g$, the equation
\begin{equation*}
    \mathfrak Du = f \ \text{in } O; \ \ u = g \ \text{on } \partial O
\end{equation*}
{has a unique solution $C(\ol O)$. If $f$ is smooth, then $u$ is also smooth by~\cite[Theorem 1.1]{Hor67}.  {Since in our choice of $\cO$, for each point $x_0$ on $\partial\cO$, we are able to find an exterior normal vector at $x_0$ whose first coordinate is nonzero, and for the operator $A$ in Lemma~\ref{lm:hypoelliptic}, $a_{ij}\neq0$ only when $i=j=1$, the domain $\cO$ satisfies the constraints for Bony's theorem.} In particular,   combined with Lemma~\ref{lm:hypoelliptic} we have}
\begin{theorem}\label{thm:bony}
    Let $\lambda>0$. Then for any $f\in C^\infty(\cO)\cap C(\ol\cO)$ and $g\in C(\ol\cO)$, the equation
    \begin{equation}\label{eq:bony}
    (A-\lambda)u = f \ \text{in } \cO; \ \ u = g \ \text{on } \partial \cO
\end{equation}
admits a unique solution $u\in C^\infty(\cO)\cap C(\ol\cO)$.
\end{theorem}
\begin{lemma}\label{lm:inte-by-part}
    For any $\varphi\in D(\wt A)$ and $g\in C_c^\infty(\cO)$, $\inner{\wt A\varphi}{g} = \inner{\varphi}{A^*g}$.
\end{lemma}
\begin{proof}
    By definition, $\varphi,\wt A\varphi\in C(\ol\cO)$. Fix $\lambda>0$. Take $\psi_m\in C^\infty(\ol\cO)$ such that $\|\psi_m-(\wt A-\lambda)\varphi\|_{C(\ol\cO)}\to 0$ as $m\to\infty$. By Theorem~\ref{thm:bony}, the equation 
    \begin{equation}\label{eq:bony-1}
    (A-\lambda)u_m = \psi_m \ \text{in } \cO; \ \ -\lambda u_m = \psi_m \ \text{on } \partial \cO
\end{equation}
has a unique solution $u_m\in  C^\infty(\cO)\cap C(\ol\cO)$. Then the boundary condition in~\eqref{eq:bony-1} implies that $Au_m(\mathbf x)$ converges to 0 as $\mathbf{x}$ goes to the boundary. Therefore by Lemma~\ref{lm:generator}, $u_m\in D(\wt A)$ and $\wt A u_m = Au_m$. In particular,
\begin{equation}\label{eq:closed-operator}
    (\wt A-\lambda)u_m = (A-\lambda)u_m = \psi_m\to (\wt A-\lambda)\varphi \ \text{in } C(\ol\cO).
\end{equation}
On the other hand, by the maximal principle, $\|u_m-u_n\|_{C(\ol\cO)}\le \lambda^{-1}\|\psi_m-\psi_n\|_{C(\ol\cO)}$, which implies that there exists some $u\in C(\ol\cO)$ such that $u_m\to u$ in $C(\ol\cO)$. Since the infinitesimal generator of a Feller semigroup is closed (see e.g.~\cite[Theorem 34.4]{laxbook}), ~\eqref{eq:closed-operator} together with $u_m\to u$ implies that $u\in D(\wt A)$ and $(\wt A-\lambda)u = (\wt A-\lambda)\varphi$. Since $\wt A-\lambda$ is the inverse of the resolvent operator (see e.g.~\cite[Proposition 6.12]{legallbook}), this further implies that $u = \varphi$.

Now by the dominant convergence theorem, we have
\begin{equation*}
    \inner{(\wt A-\lambda)\varphi}{g} = \lim_{m\to\infty}\inner{(A-\lambda)u_m}{g} = \lim_{m\to\infty}\inner{u_m}{(A^*-\lambda)g} = \inner{u}{(A^*-\lambda)g} = \inner{\varphi}{(A^*-\lambda)g},
\end{equation*}
and we conclude the proof by subtracting $\lambda\inner{\varphi}{g}$.
\end{proof}

\begin{lemma}\label{lm:inte-by-part-2}
    For any $\varphi\in C(\ol\cO)$ and $g\in C_c^\infty(\cO)$, $\inner{P_t\varphi}{g} = \inner{\varphi}{g} + \int_0^t\inner{P_s\varphi}{A^*g}\,ds$.
\end{lemma}
\begin{proof}
    Since the domain of $\wt A$ is dense (see e.g.~\cite[Theorem 34.4]{laxbook}, we may pick $\varphi_m\in D(\wt A)$ such that $\varphi_m\to\varphi$ in $C(\ol\cO)$. Then for each $m$, by~\cite[Proposition 6.11]{legallbook}, $P_s\varphi_m\in D(\wt A)$ and $P_t\varphi_m = \varphi_m + \int_0^t \wt A(P_s\varphi_m)ds$. Then  by Lemma~\ref{lm:inte-by-part}, 
\begin{equation*}
    \begin{split}
       \inner{P_t\varphi}{g} &= \lim_{m\to\infty} \inner{P_t\varphi_m}{g} = \lim_{m\to\infty}\bigg(\inner{\varphi_m}{g} + \int_0^t  \inner{\wt A(P_s\varphi_m)}{g}\,ds \bigg)\\&= \lim_{m\to\infty}\bigg(\inner{\varphi_m}{g} + \int_0^t \inner{P_s\varphi_m}{A^*g}\,ds\bigg) \\&= \inner{\varphi}{g} + \int_0^t\inner{P_s\varphi}{A^*g}ds.
    \end{split}
\end{equation*}
\end{proof}

\begin{proof}[Proof of Proposition~\ref{prop:weaksol}]
    By Lemma~\ref{lm:finite-mg}, for $t>0$ and a.e.\ $\mathbf{x}\in\ol\cO$, $F(\mathbf{x}) = \bbE^{\mathbf{x}}F(X_{t\wedge\sigma})$.  {By Proposition~\ref{prop:induc-finite}, $\mathcal{Z}_\alpha(\bbH;x_1,...,x_{2N+2})$  is lower semicontinuous and locally integrable, which implies that} $F$ is lower semicontinuous and locally integrable.  {This enables us to pick} $f_m\in C(\ol\cO)$ with $f_n\uparrow F$. Then by Lemma~\ref{lm:inte-by-part-2} and the dominant convergence theorem, for any $g\in C_c^\infty(\cO)$,
    \begin{equation}
    \begin{split}
        \inner{F}{g} &= \int_{\cO}\bbE^{\mathbf{x}}F(X_{t\wedge\sigma})g(\mathbf{x})\,d\mathbf x=\lim_{m\to\infty} \int_{\cO}\bbE^{\mathbf{x}}f_m(X_{t\wedge\sigma})g(\mathbf{x})\,d\mathbf x \\
        &=\lim_{m\to\infty}\inner{P_tf_m}{g} = \lim_{m\to\infty}\bigg(\inner{f_m}{g}+\int_0^t\inner{P_sf_m}{A^*g}\,ds  \bigg)\\
        &=\lim_{m\to\infty}\bigg(\inner{f_m}{g}+\int_0^t\int_\cO\bbE^{\mathbf{x}}f_m(X_{s\wedge\sigma}){A^*g}(\mathbf x)\,d\mathbf x\, ds  \bigg)\\
        &= \inner{F}{g}+\int_0^t\int_\cO\bbE^{\mathbf{x}}F(X_{s\wedge\sigma}){A^*g}(\mathbf x)\,d\mathbf x\, ds  = \inner{F}{g} + t\inner{F}{A^*g}.
        \end{split}
    \end{equation}
    Therefore $\inner{F}{A^*g} = 0$ for any $g\in C_c^\infty(\cO)$, which concludes the proof.
\end{proof}

\bibliographystyle{alpha}
\bibliography{theta}

\end{document}